%% file: actcongr.tex
\magnification=\magstephalf
\input amstex
\documentstyle{amsppt}
\pagewidth{6.5truein}
\pageheight{8.9truein}
\ifx\refstyle\undefinedZQA\else\refstyle{C}\fi

\input plot
\plottrue
\plotPStype1
\input plotsupp

\define\X{{\Cal X}}
\define\Z{{\bold Z}}
\define\N{{\bold N}}
\define\nullset{\varnothing}
\define\QED{ \null\nobreak\hfill$\blacksquare$}
\define\QEd{ \null\nobreak\hfill$\square$}
\define\procl#1{\smallskip{\it {#1}. }}
\define\pp{{\bar p}}
\define\NN{{\bar N}}
\define\axis{\ell}
\define\orbit{{\Cal O}}
\define\isom{\phi}
\define\isomgroup{{O_3}}
\define\eps{\varepsilon}
\define\ri{\therosteritem}
\define\scong{\preceq}
\define\CP#1{{CP_{#1}}}
\define\UNC#1{{UNC_{#1}}}
\define\setgraph{{\Cal G}}
\define\Fgraph{{\Cal F}}

\define\Adams{1}
\define\Dougherty{2}
\define\DoughertyForeman{3}
\define\Hausdorff{4}
\define\Magnus{5}
\define\Nickolas{6}
\define\Robinson{7}
\define\Wagon{8}

\topmatter
\title Open sets satisfying systems of congruences\endtitle
\author Randall Dougherty\endauthor
\affil Ohio State University\endaffil
\date December 21, 1999 \enddate
\address Department of Mathematics, Ohio State University,
Columbus, OH 43210\endaddress
\email rld\@math.ohio-state.edu \endemail
\subjclass Primary: 52B45\endsubjclass
\keywords
Banach-Tarski paradox, congruences, free groups
\endkeywords

\abstract
A famous result of Hausdorff states that a sphere with countably many
points removed can be partitioned into three pieces $A,B,C$ such that
$A$~is congruent to~$B$ (i.e., there is an isometry of the sphere which
sends $A$ to~$B$), $B$~is congruent to~$C$, and $A$~is congruent to $B
\cup C$; this result was the precursor of the Banach-Tarski paradox.
Later, R.~Robinson characterized the systems of congruences like this
which could be realized by partitions of the (entire) sphere with rotations
witnessing the congruences.  The pieces involved were nonmeasurable.

In the present paper, we consider the problem of which systems of
congruences can be satisfied using {\sl open} subsets of the sphere (or
related spaces); of course, these open sets cannot form a partition of
the sphere, but they can be required to cover `most of' the sphere in the
sense that their union is dense.  Various versions of the problem arise,
depending on whether one uses all isometries of the sphere or restricts
oneself to a free group of rotations (the latter version generalizes
to many other suitable spaces), or whether one omits the requirement
that the open sets have dense union, and so on.  While some cases of
these problems are solved by simple geometrical dissections, others
involve complicated iterative constructions and/or results from the theory
of free groups.  Many interesting questions remain open. \endabstract

\endtopmatter
\document

\head 1. Introduction \endhead

Can one find four nonempty, pairwise disjoint open subsets of a sphere
such that the union of any two is congruent to the union of any
other two?  What about five such sets?  Six?  Seven?

This quite concrete geometrical question, and many similar questions,
arose as an offshoot of a study of questions related to the Banach-Tarski
paradox.  More directly, they are related to the following theorem of
Hausdorff~\cite{\Hausdorff, p.~469}, which led to the Banach-Tarski
result:  There is a countable subset~$D$ of the sphere~$S^2$ such that
$S^2\setminus D$ can be partitioned into three sets $A,B,C$ such that $A$
is congruent to~$B$ (i.e., there is an isometry $\rho$ such that $\rho(A)
= B$), $B$~is congruent to~$C$, and $C$ is congruent to $A \cup B$.
(It is easy to see that the sets~$A,B,C$ cannot be measurable with
respect to the standard isometry-invariant probability measure on~$S^2$.)

We will consider various systems of congruences like the one given
above; it will help to fix some notation and terminology now.

Fix a positive integer~$r$.  A {\it congruence\/} is specified by two
subsets $L$ and~$R$ of $\{1,2,\dots,r\}$, and is written formally as
$\bigcup_{k\in L} A_k \cong \bigcup_{k\in R} A_k$, where $A_1,A_2,\dots,A_r$
are variables.  The congruence is {\it proper\/} if both $L$ and~$R$
are nonempty proper subsets of $\{1,\dots,r\}$.  Now suppose $G$~is a
group acting on a set~$X$, and a system of congruences is given by pairs
$L_i,R_i\subseteq\{1,\dots,r\}$ for $i\le m$. 
Then a given sequence of sets $A_k\subseteq X$
($k\le r$) is said to {\it satisfy} the system of congruences
if the sets~$A_k$ are pairwise disjoint and, for each
$i\le m$, there is $\sigma_i\in G$ such that $\sigma_i(\bigcup_{k\in
L_i} A_k) = \bigcup_{k\in R_i} A_k$ (i.e., $\sigma_i$~witnesses
congruence number~$i$).

Of course, if the sets~$A_k$ are all empty, then they trivially
satisfy any system of congruences.  The opposite extreme case is when
the sets~$A_k$ form a partition of~$X$; in this case, if they
satisfy the system of congruences, they are said to be a {\it solution}
to the system.

The argument of Hausdorff generalizes to show that, for any system of
proper congruences, there are subsets of~$S^2$ (which is acted on by its
rotation group) satisfying the congruences and having union~$S^2\setminus
D$, where $D$~is countable.  (See chapter~4 of Wagon~\cite{\Wagon}.)
One cannot always eliminate the countable exceptional set here; Raphael
Robinson~\cite{\Robinson} characterized the systems of congruences which
actually have solutions on $S^2$~with its rotation group.

The above constructions produce extremely wild sets; in the case
of the Banach-Tarski paradox, it is easy to see that the
construction cannot be performed using measurable sets.
Marczewski asked whether a Banach-Tarski decomposition could
be produced using sets with the property of Baire; this question
was answered affirmatively by Dougherty and Foreman~\cite{\DoughertyForeman}.
A characterization of which systems of congruences have solutions
in~$S^2$ with its rotation group (or related spaces)
using sets with the property of Baire is given
in Dougherty~\cite{\Dougherty}.  In both cases, the results for
sets with the property of Baire are obtained from constructions
of {\it open} sets which `almost' satisfy the decomposition
equations or congruences, in the sense that a meager exceptional
set is allowed for each equation or congruence, and also in the
partition(s) of~$S^2$.

This naturally leads to the question of whether one can find
open sets which actually satisfy a system of congruences,
without exceptional sets.  Of course, one cannot require such
sets to form a partition of the space (the sphere~$S^2$, being
connected, cannot be partitioned nontrivially into open sets),
but one can, if one chooses, require the sets to fill `almost all'
of the space in the sense that their union is dense (so the leftover
set is nowhere dense and hence meager).  Such questions are the
focus of this paper.

The reason that the sphere is a good space to study such congruences on
is the same reason that the Hausdorff and Banach-Tarski paradoxes
apply to it --- the rotation group of the sphere has a subgroup which
is a free group on two generators.  Most of the open-set results here
and in the previous papers above apply in a more general context:

\definition{Definition 1.1} A {\it suitable space} is a pair $(\X,G)$
where $\X$~is a complete separable metric space and $G$~is a countable
group acting on~$\X$ by homeomorphisms such that $G$~is a free group
on more than one generator and $G$~acts freely on a comeager subset
of~$\X$.  (Equivalently, for each $g\in G$ other than the identity, the
set of fixed points of~$g$ has empty interior.) \enddefinition

Because of this generalization, we will pay more attention to the
case of the sphere acted on by a free group of rotations
than to the case of the sphere with its entire group of isometries.
Results for the free group case will often generalize to a wide
variety of other suitable spaces (for instance, the group of
bi-Lipschitz homeomorphisms from the Cantor space to itself has a
subgroup which is free on two generators and acts freely on the
Cantor space; we will see other examples of suitable spaces later);
results for the all-isometry case are more isolated.
So it will be of interest to show that a system of congruences cannot
be satisfied on the sphere with elements of a free group of rotations
witnessing the congruences, even when it is easy to get open subsets
of the sphere satisfying the congruences via other isometries.

For a few examples (especially in the case of the sphere with all
isometries), the open sets satisfying certain congruences will be given
by simple dissections.  In other cases, though, the open sets will be
produced by iterative constructions and will be quite complicated, with
infinitely many connected components and often having boundaries of
positive measure.

We will use the symbol~$\circ$ or simple juxtaposition to denote
a group operation, interchangeably.  All group actions will
be written on the left.
For standard basic facts about free groups, such as the unique
expression of any element as a reduced word in the generators and the
fact that any nonidentity element has infinite order, see any text on
combinatorial group theory, such as Magnus, Karrass, and
Solitar~\cite{\Magnus}.  More advanced facts will be referred to
specifically as needed. 
For instance, every subgroup of a free group is free~\cite{\Magnus,
Cor.~2.9}.  Also, a free group on two generators has subgroups
which are free on $n$~generators for any given natural number~$n$
\cite{\Magnus, Prob.~1.4.12}; hence, the group acting on a suitable
space has such subgroups.

\head 2. Basic restrictions \endhead

We start here by describing some properties that a system of congruences
must have in order to be satisfied nontrivially by open sets in the
contexts we are studying.

From the congruences in a given system, one can deduce other congruences.
The fact that the mappings witnessing congruences form a group
means that congruence is an equivalence relation --- the identity
mapping is used to show that $\cong$~is reflexive, inverses give
symmetry of~$\cong$, and composition gives transitivity.  If we are
considering sets which form a partition of the space (i.e., solutions
to the system of congruences), then there is a complementation rule:
from $\bigcup_{k\in L} A_k \cong \bigcup_{k\in R} A_k$
we can deduce $\bigcup_{k\in L^c} A_k \cong \bigcup_{k\in R^c} A_k$
(where $S^c = \{1,\dots,r\} \setminus S$), because the mapping witnessing
the congruence of two sets also witnesses that their complements
are congruent.  A system of congruences is called {\it weak}
if one cannot deduce any self-complementary congruence
$\bigcup_{k\in L} A_k \cong \bigcup_{k\in L^c} A_k$ from it
by the equivalence relation rules and the complementation rule.

It is easy to see that, if a system of congruences has a solution in~$S^2$
with rotations witnessing the congruences, then the system must be weak:
any rotation has fixed points, and hence cannot witness that a set is
congruent to its complement.  Robinson showed that the converse is true:
any weak system of congruences has a solution in~$S^2$ with rotations
witnessing the congruences (using unrestricted pieces in the partition).

If we are not requiring the sets to form a partition of the space
(as noted earlier, we cannot require this for open subsets of the sphere),
then the complementation rule need not hold,
and a system of congruences need not be weak in order to be satisfied.
For instance, the simple system
with $r=2$ and the single congruence $A_1 \cong A_2$ is clearly not
weak, but it is satisfied on~$S^2$ by two complementary open hemispheres.
Or one can just use two smaller disks; these will not have dense union,
but one can use a rotation from a free group to witness the congruence.
(For the hemispheres one would have to use a rotation of order~$2$.)

However, if we want to get open subsets of the sphere with dense union
to satisfy a system of congruences, with rotations from a free group
witnessing the congruences, then the system must be weak.
This was proved in Dougherty~\cite{\Dougherty} (such sets would
form a `quasi-solution' to the system in the sense of that paper).

The proof referred to above uses the following easy result which will also
be needed here:
 
\proclaim{Lemma 2.1 \rm\cite{\Dougherty, Lemma 3.2}}
If an open subset~$A$ of~$S^2$ is invariant under a
rotation of infinite order around an axis~$\axis$, then $A$~is invariant
under {\it all} rotations around~$\axis$.  The same is true
if `invariant' is replaced by `quasi-invariant'
(where $A$~is quasi-invariant under~$\rho$ iff $A$~differs
from~$\rho(A)$ by a meager set).\QED\endproclaim

Hence, if the open set is invariant under rotations of infinite order
around two different axes, then the set must be either empty or the entire
sphere.

Next, say that $B$ is {\it subcongruent} to $C$ ($B \scong C$) if $B$ is
congruent to a subset of~$C$.  From a given system of congruences, one
can deduce subcongruences by the following rules:
if $L \subseteq R$,
then $\bigcup_{k\in L} A_k \scong \bigcup_{k\in R} A_k$;
if $B \scong C$ and $C \scong D$, then $B \scong D$;
and, if $B \cong C$ is in the given system, then $B \scong C$
and $C \scong B$.  Again, there is a complementation rule in the
case where the sets form a partition of the space:
if $\bigcup_{k\in L} A_k \scong \bigcup_{k\in R} A_k$,
then $\bigcup_{k\in R^c} A_k \scong \bigcup_{k\in L^c} A_k$.
Call the system {\it
consistent\/} if there do not exist sets $L,R \subseteq \{1,2,\dots,r\}$
with $R$~a proper subset of~$L$ such that one can deduce $\bigcup_{k \in
L} A_k \scong \bigcup_{k \in R} A_k$ from the system by the above rules.
(For example, the
Hausdorff system $A_1 \cong A_2 \cong A_3 \cong A_2 \cup A_3$
is not consistent.)  Note that any consistent system must consist entirely
of proper congruences, if one ignores trivial identity congruences
such as $\nullset \cong \nullset$.

The main result of Dougherty~\cite{\Dougherty} states that a
system of congruences has a solution on the sphere under its rotation
group using nonmeager sets with the property of Baire if and only if
the system is weak and consistent.

If we consider sets which do not form a partition, then again
the complementation rule no longer applies.  Nonetheless, if
there are nonempty open subsets of the sphere satisfying a system
of congruences (even using arbitrary isometries), then
the system must be consistent.  In fact, even more must hold
in this case.

The reason is the standard isometry-invariant
probability measure on the sphere, which gives every nonempty
open set positive measure.  This measure gives a necessary condition
for there to be nonempty open subsets $A_1,\dots,A_r$ of the sphere
satisfying a system of congruences: there must exist positive
numbers $\mu_1,\dots,\mu_r$ such that, if
$\bigcup_{k\in L} A_k \cong \bigcup_{k\in R} A_k$ is in the system,
then $\sum_{k \in L} \mu_k = \sum_{k \in R} \mu_k$.
(If one wants to allow some of the sets~$A_k$ to be empty, then
one can allow some of the numbers~$\mu_k$ to be~$0$.)

Call a system for which there exist positive numbers~$\mu_k$ as above
{\it numerically consistent}.  A numerically consistent system must
be consistent, because, for each
subcongruence $\bigcup_{k\in L} A_k \scong \bigcup_{k\in R} A_k$
deducible from the system, we have
$\sum_{k \in L} \mu_k \le \sum_{k \in R} \mu_k$.  (Even the complementation
rule preserves this, because we have
$\sum_{k \in L^c} \mu_k = M - \sum_{k \in L} \mu_k$, where
$M = \sum_{k=1}^r \mu_k$.)  This inequality cannot hold
if $R$~is a proper subset of~$L$, so no such subcongruence is deducible.

However, numerical consistency is strictly stronger than consistency.
For example, consider the system $A_1 \cong A_2 \cong A_3 \cong A_4
\cong A_5$, $A_1 \cup A_2 \cong A_1 \cup A_3 \cup A_4$. It is easy to
show that this system is weak and consistent. But there do not exist
positive numbers~$\mu_k$ as above; they would all have to be the same
number~$\mu$, and then the last congruence would give $2\mu = 3\mu$ and
hence $\mu = 0$.

Among the numerically consistent systems of congruences, the following
systems (one for each $s \ge 1$) can be singled out:
$$\UNC s: \qquad \bigcup_{j\in L} B_j \cong \bigcup_{j\in R} B_j,
\quad L,R \subseteq \{1,2,\dots,s\},\quad |L|=|R|$$
System $\UNC s$ states that the sets $B_1,\dots,B_s$ are such that,
for each $m \le s$, any two unions of $m$~of the sets are congruent
to each other.  The system is clearly numerically consistent,
with $\mu_j = 1$ for each $j \le s$.  We will now see that the
systems~$\UNC s$ form a `universal' family of numerically consistent
systems of congruences.

Suppose we have a system of congruences on sets $A_1,\dots,A_r$ and
another system of congruences on sets $B_1,\dots,B_s$. We say that the
first system is {\it reducible} to the second system if there is a
function~$\pi$ from $\{1,2,\dots,s\}$ to $\{1,2,\dots,r\}$ such that,
for each $L,R\subseteq \{1,2,\dots,r\}$, if
$$\bigcup_{k\in L} A_k \cong \bigcup_{k\in R} A_k$$
is in the first system, then
$$\bigcup_{j:\;\pi(j)\in L} B_j \cong \bigcup_{j:\;\pi(j)\in R} B_j$$
is in the second system.  So, if we have sets~$B_j$ satisfying the
second system, we can get sets $A_k$ satisfying the first system by
letting $A_k = \bigcup_{j:\;\pi(j)=k}B_j$.
If $\pi$ maps $\{1,2,\dots,s\}$ {\sl onto} $\{1,2,\dots,r\}$,
then the reduction preserves nonemptiness: if the sets~$B_j$
are all nonempty, then the resulting sets~$A_k$ will also be nonempty.

\proclaim{Proposition 2.2}
A system of congruences is numerically consistent if and
only if it is reducible to~$\UNC s$ for some~$s$ by some function~$\pi$
from $\{1,2,\dots,s\}$ onto $\{1,2,\dots,r\}$.
\endproclaim

\demo{Proof}
For the `if' part, it suffices to show that reducibility via an onto
function preserves numerical consistency. Suppose that a system of
congruences on $A_1,\dots,A_r$ is reducible to a system of congruences
on $B_1,\dots,B_s$ via the onto function~$\pi$. Suppose we have positive
numbers~$\lambda_j$ for $j \le s$ witnessing that the second system is
numerically consistent.  Then we can get positive numbers $\mu_k$ for $k
\le r$ by letting $\mu_k = \sum_{j:\;\pi(j) = k} \lambda_j$, and these
numbers will witness the numerical consistency of the first system.

For the `only if' part, suppose we are given a system of congruences
on $A_1,\dots,A_r$ and positive numbers $\mu_k$, $k \le r$, witnessing
that the system is numerically consistent.  This means that the
numbers~$\mu_k$ satisfy certain linear equations with integer
coefficients.  It now follows
from standard linear algebra results that we can get positive {\sl
rational} numbers $\mu_k$ satisfying these equations.  Then, since we
can multiply through by a common denominator, we may assume that the
numbers~$\mu_k$ are actually positive integers.

Let $s = \sum_{k = 1}^r \mu_k$, and let $\pi$ be a function from
$\{1,2,\dots,s\}$ to $\{1,2,\dots,r\}$ such that, for each $k \le r$,
$k$~has exactly $\mu_k$ preimages under~$\pi$ in $\{1,2,\dots,s\}$.
Since the numbers $\mu_k$ are all nonzero, $\pi$ is surjective.
For each $L,R\subseteq \{1,2,\dots,r\}$, if
$$\bigcup_{k\in L} A_k \cong \bigcup_{k\in R} A_k$$
is in the given system, then $\sum_{k \in L} \mu_k
= \sum_{k \in R} \mu_k$; it follows that
$$\bigcup_{j:\;\pi(j)\in L} B_j \cong \bigcup_{j:\;\pi(j)\in R} B_j$$
is in the system~$\UNC s$.  So $\pi$~reduces the given system
to~$\UNC s$.
\QED\enddemo

So, to show that all numerically consistent systems
are satisfiable (by nonempty sets) in a certain space, it suffices to show
that the systems~$\UNC s$ are all satisfiable (by nonempty sets).  Note that
$\UNC s$~is weak for odd~$s$ but not for even~$s$.
(But not every numerically consistent weak system
is reducible to~$\UNC s$ for an odd~$s$ --- for instance,
look at the system $A_1 \cong A_2 \cong A_3 \cong A_4$.)

\head 3. Initial results \endhead

We are now ready to consider the satisfiability of some particular
systems of congruences using open subsets of the sphere.
Let us start with the systems~$\UNC r$ from the preceding section.
We also consider the following natural subsystem of~$\UNC r$:

$$\CP r: \qquad \bigcup_{k\in L} A_k \cong \bigcup_{k\in R} A_k,
\quad L,R \subseteq \{1,2,\dots,r\},\quad |L|=|R|=2$$

This is just the ``$r$ sets, with the union of any two congruent to
the union of any other two'' system mentioned at the beginning of
section~1.  (It is weak if $r \ne 4$.)

For small enough~$r$ it is easy to produce open subsets of the sphere with
dense union satisfying $\CP r$ and~$\UNC r$.  For $r=1$, let $A_1$ be the
whole sphere; for $r=2$, let $A_1$ and~$A_2$ be complementary hemispheres.
For $r=3$, one can divide a sphere into three $120^\circ$ lunes by
three equally-spaced meridians, and these sets will satisfy~$\UNC3$.
For $r=4$, one can get the desired sets by radially projecting the faces
of a regular tetrahedron to its circumscribing sphere.

A slightly more complicated construction yields sets satisfying $\CP5$
and~$\UNC5$. The faces of a regular icosahedron can be partitioned into
five sets of four such that two faces in the same set do not touch, even
at a vertex; in fact, there are exactly two such partitions, one a
mirror image of the other.  (The arrangement of triangles in one such
set is unique up to rotation and reflection; given one such set, the
other four in the partition can be obtained by rotating the first set
around a vertex of the icosahedron.) The five sets in such a partition
can be projected to the circumscribing sphere to yield five open sets
$A_1,\dots,A_5$ (each with four components); these open sets
satisfy~$\UNC5$ (and, in particular, $\CP5$).

Whether $\UNC6$, or even~$\CP6$, is satisfied on the sphere by nonempty
open sets is not yet known.  One possible way to prove that such
open sets do not exist would be to show that the isometries witnessing
the congruences would have to satisfy enough group-theoretic relations
that the group generated by them could not be a subgroup of the
isometry group of the sphere; Michael Larsen (personal communication)
has suggested an approach along these lines.

The constructions above all make use of finite-order rotations
to witness congruences.  If one wants to restrict the isometries used
to a free group of rotations, then the problem becomes quite different.
It is easy to satisfy~$\UNC2$, which is just the congruence
$A_1 \cong A_2$: let $\sigma$ be a non-identity rotation in
the group, let $x$ be a point of the sphere not fixed by $\sigma$,
let $A_1$ be a neighborhood of~$x$ so small that $A_1$ and~$\sigma(A_1)$
are disjoint, and let $A_2 = \sigma(A_1)$.  (One cannot arrange
for the sets $A_1$ and~$A_2$ to have union dense in the sphere,
because $\UNC2$~is not weak; see section~2.)
But the system~$\CP3$ cannot be satisfied:

\proclaim{Theorem 3.1} Suppose $A_1$, $A_2$, and~$A_3$ are disjoint
open subsets of~$S^2$ such that $A_1 \cup A_2 \cong A_1 \cup A_3 \cong
A_2 \cup A_3$, and these congruences are witnessed by elements of a
free group of rotations of~$S^2$.  Then $A_1 = A_2 = A_3 = \nullset$.
\endproclaim

\demo{Proof}  We will use the following group-theoretic facts: If two
elements $g$ and~$g'$ of a free group~$G$ commute, then there are integers
$a$ and~$b$ and an element~$h$ of~$G$ such that $g = h^a$ and $g' = h^b$,
and hence $g^b = {g'}^a$~\cite{\Magnus,~Prob.~1.4.6}.  If $g$ and~$g'$
do not commute, then they are free generators for a free subgroup of~$G$
of rank~$2$ \cite{\Magnus, Cors. 2.11 and~2.13.1}.  (The rank of a free
group is the number of generators in a free generating set for the group;
this is well-defined~\cite{\Magnus,~Thm.~2.4}.)

There are two cases to consider.  First, suppose the rotations
witnessing $A_1 \cup A_2 \cong A_1 \cup A_3$ and $A_1 \cup A_2 \cong A_2
\cup A_3$ commute; then they are both powers of some rotation~$\sigma$
of infinite order.  By permuting the indices $1,2,3$, we may arrange to
have $\sigma^m(A_1 \cup A_3) = A_1 \cup A_2$ and $\sigma^n(A_1 \cup A_2)
= A_2 \cup A_3$ where $m$ and~$n$ are nonnegative integers.  In fact, we
may assume that $m$ and~$n$ are positive; if one congruence were
witnessed by the identity, this would force two of the three sets to be
empty, and the other congruence would force the third set to be empty as
well.  Now, suppose $C$~is a connected component of~$A_3$.  Then $C$~is
a component of $A_1 \cup A_3$.  (Since $\{A_1,A_3\}$ is a partition of
$A_1 \cup A_3$ into open sets, the components of $A_1 \cup A_3$ are just
the components of~$A_1$ and the components of~$A_3$.) So
$\sigma^m(C)$~is a component of $A_1 \cup A_2$ and hence a component of
either~$A_1$ or~$A_2$. Similarly, if $C$~is a component of either~$A_1$
or~$A_2$, then $\sigma^n(C)$ is a component of either~$A_2$ or~$A_3$.
Applying these two facts repeatedly, we find that, if $C$~is a component
of any one of the three sets, then there is an infinite increasing
sequence of positive integers~$a$ such that $\sigma^a(C)$ is also a
component of one of the three sets. All of the sets $\sigma^a(C)$ have
the same positive measure (using the standard measure on~$S^2$); since
$S^2$~has finite measure, the relevant sets $\sigma^a(C)$ cannot all be
disjoint, so there are positive integers $a < b$ such that the sets
$\sigma^a(C)$ and $\sigma^b(C)$ overlap and are each a component of one
of the sets~$A_i$.  Since the sets $A_i$~are disjoint, $\sigma^a(C)$ and
$\sigma^b(C)$ must be components of the same set~$A_i$; since they
overlap, we must actually have $\sigma^a(C) = \sigma^b(C)$.
Applying~$\sigma^{-a}$ gives $C = \sigma^{b-a}(C)$.  Since
$\sigma^{b-a}$~is a rotation of infinite order, $C$~must be invariant
under all rotations around the axis of~$\sigma^{b-a}$, by Lemma~2.1.  In
particular, $\sigma(C) = C$.

We have now seen that all components of all of the sets~$A_i$ are
invariant under~$\sigma$, so the sets themselves are invariant
under~$\sigma$ and hence under $\sigma^m$ and~$\sigma^n$.  Therefore,
$A_1 \cup A_2 = A_1 \cup A_3 = A_2 \cup A_3$; since the sets~$A_i$ are
disjoint, they must be empty. This completes the first case.

For the remaining case, suppose $\sigma(A_1 \cup A_3) = A_1 \cup A_2$
and $\tau(A_2 \cup A_3) = A_1 \cup A_2$, where $\sigma$ and~$\tau$
do not commute and hence are free generators for their subgroup.  We will
show that $A_3 = \nullset$.  By permuting the indices, one can use the
same proof to get $A_1 = A_2 = \nullset$ (if $\sigma$ and~$\tau$ do not
commute, then $\sigma^{-1}$ and $\sigma^{-1} \circ \tau$ do not commute);
alternatively, one can use Proposition~3.1 of Dougherty~\cite{\Dougherty}
to complete the proof.

Suppose $A_3 \ne \nullset$; then $A_3$~has a connected component~$C$.
As before, we see that both $\sigma(C)$ and $\tau(C)$ are components of
either $A_1$ or~$A_2$.  Furthermore, if $C'$~is a component of~$A_1$,
then $\sigma(C')$ is a component of either~$A_1$ or~$A_2$; if $C'$~is a
component of~$A_2$, then $\tau(C')$ is a component of either~$A_1$
or~$A_2$.  These latter facts can be applied repeatedly starting at
$\sigma(C)$ to get increasingly long compositions of $\sigma$ and~$\tau$
which, when applied to $\sigma(C)$, give components of
either~$A_1$ or~$A_2$.  Eventually two such components must overlap and
hence coincide.  Therefore, there exist non-identity words $u$ and~$v$
in the generators $\sigma$ and~$\tau$, with no inverse powers of
$\sigma$ or~$\tau$ occurring, such that the rightmost term in~$u$
is~$\sigma$, and $v(u(C)) = u(C)$.  Let $w = u^{-1} \circ v \circ u$; then
$C$~is fixed under~$w$, and the reduced form of~$w$ consists of
negative terms ($\sigma^{-1}$ and~$\tau^{-1}$) followed by positive
terms ($\sigma$ and~$\tau$), with more positive terms than negative
terms, and with rightmost term~$\sigma$.

The same procedure starting with $\tau(C)$ leads to a word~$w'$ such
that $C$~is fixed under~$w'$, the reduced form of~$w'$ consists of
negative terms followed by positive terms, with more positive terms than
negative terms, and the rightmost term of~$w'$ is~$\tau$. It is now easy
to see that the reduced form of~$ww'$ ends in~$\tau$, while that
of~$w'w$ ends in~$\sigma$. So $w$ and~$w'$ do not commute; hence, as
rotations, they must have different axes.  But they both have infinite
order, so $C$~is fixed under any rotation around either of these axes,
by Lemma~2.1.  As noted after that lemma, it now follows that $C$~must
be all of~$S^2$.  So $A_3 = S^2$ and $A_1 = A_2 = \nullset$. This
clearly does not satisfy the congruences, so we have a contradiction.
Therefore, $C$~cannot exist, so $A_3$~is empty, as desired. \QED\enddemo

Attempts to generalize Theorem~3.1 lead to the study of solutions
to congruences in terms of {\sl finite} sets, in~$S^2$ or in free
groups themselves.

\proclaim{Theorem~3.2} For any system of congruences, the following
are equivalent:
\roster
\item"(I)" There are open subsets of~$S^2$, not all empty, which satisfy
the congruences, with the witnesses coming from a free group of rotations.
\item"(II)" There are finite subsets of~$S^2$, not all empty, which satisfy
the congruences, with the witnesses coming from a free group of rotations.
\item"(III)" For some free group~$F$ and some element~$w$ of~$F$ which
is not a proper power of another element, there are finite subsets of
$F/\langle w \rangle$ (the set of left cosets of the subgroup\/
$\langle w \rangle$ generated by~$w$), not all empty, which satisfy the
congruences (with witnesses from~$F$).
\endroster
\endproclaim

\demo{Proof} To see that (II) implies (I), suppose the sets~$A_k$ are
disjoint finite sets which satisfy the congruences, and choose $\eps >
0$ so small that the distance between any two points in $\bigcup_k A_k$
is greater than~$2\eps$.  Then the sets $B_k = \{x\colon d(x,A_k) <\nobreak
\eps\}$ are disjoint open sets which satisfy the congruences.

For the proof that (I) implies (II), we first eliminate a trivial case.
Suppose that one of the sets~$A_k$ occurs on both sides of any
congruence of the system in which it appears at all.  Then (II) clearly
holds, because we can let this set~$A_k$ be a single point of~$S^2$ and
all other sets be empty; all congruences would then be witnessed by the
identity map.  So from now on, assume that each set~$A_k$ occurs on only
one side of some congruence; it follows immediately that any nontrivial
solution to the system must have at least two sets nonempty.

As we saw in the proof of Theorem~3.1, the following fact follows
easily from the definition of connectedness:  if $E \subseteq S^2$ is the
union of disjoint open sets~$E_i$, then any connected component of~$E$
is included in one of the sets~$E_i$, and is a component of
that~$E_i$.  Also, since $S^2$~is locally connected, any
component of an open subset of~$S^2$ is open, and therefore has positive
measure under the standard isometry-invariant probability measure on~$S^2$.

Suppose that the open sets~$A_k$ satisfy the congruences, as specified
by~(I), and let $G$ be the free group of rotations which includes
witnesses to these congruences.  Choose a component~$C_0$ of one of the
sets~$A_k$, and let $W$ be the collection of all components of the
sets~$A_k$ which are congruent to~$C_0$ as witnessed by a rotation
in~$G$.  Then $W$~is a collection of pairwise disjoint open sets which
all have the same positive measure, so $W$~is finite (but nonempty,
since $C_0 \in W$).

We now choose a point~$x_0$ as follows: if $C_0$~is not fixed under any
nonidentity element of~$G$, let $x_0$ be any point in~$C_0$; if $C_0$~is
fixed under some nonidentity rotation $g \in G$, let $x_0$ be one of
the two points of~$S^2$ fixed under~$g$.  (The set~$C_0$ cannot be
fixed under rotations in~$G$ around two different axes. If it were,
the remark after Lemma~2.1
would imply that $C_0$~is all
of~$S^2$.  This would make one set~$A_k$ all of~$S^2$ and the rest empty,
the case we eliminated earlier.) In the latter case, $x_0$~need not be
an element of~$C_0$.

For each $C \in W$, define $f(C)$ to be $h(x_0)$ for any $h \in G$ such
that $h(C_0) = C$.  If $h'$~is another element of~$G$ sending~$C_0$
to~$C$, then $h^{-1} \circ h'$ fixes~$C_0$, so it fixes~$x_0$, so
$h(x_0) = h'(x_0)$; hence, $f(C)$ is well defined.  We now verify that
$f$~is one-to-one.  If $C_0$~is not fixed under any nonidentity member
of~$G$, then $f(C) \in C$ for each $C$ in~$W$, and the members of~$W$
are disjoint, so $f$~must be one-to-one.  Now assume $C_0$~is fixed
under $g \in G$, $g$~not the identity.  Suppose $C$ and~$C'$ are in~$W$,
and $f(C) = f(C')$. Choose $h,h' \in G$ such that $h(C_0) = C$ and
$h'(C_0) = C'$.  Since $f(C) = f(C')$, we have $h(x_0) = h'(x_0)$, so
$(h^{-1} \circ h')(x_0) = x_0$, so $h^{-1} \circ h'$ must be a rotation
around the same axis as~$g$.  By Lemma~2.1, we have $h^{-1}(h'(C_0)) =
C_0$, so $h'(C_0) = h(C_0)$, so $C = C'$.  Therefore, $f$~is one-to-one.
It is easy to see that $f$~preserves the action of~$G$ (that is,
$f(h(C)) = h(f(C))$ for $h \in G$).

Now define~$B_k$ to be $\{f(C)\colon C \in\nobreak W,\,\,C
\subseteq\nobreak A_k\}$ for each~$k$.  The sets~$B_k$ are finite,
disjoint (since $f$~is one-to-one), and not all empty.  It remains to
see that the sets~$B_k$ satisfy the given congruences.  Suppose one of
the congruences is $\bigcup_{k \in L}A_k \cong \bigcup_{k \in R}A_k$,
and let $\rho \in G$ witness this congruence for the sets~$A_k$.  If
$x$~is a point in one of the sets~$B_k$ for $k \in L$, then $x = f(C)$
for some $C \in W$ which is a component of one of the sets~$A_k$ for $k
\in L$.  Hence, $C$~is a component of $\bigcup_{k \in L}A_k$, so
$\rho(C)$ is a component of $\rho(\bigcup_{k \in L}A_k) = \bigcup_{k \in
R}A_k$, so $\rho(C)$ is a component of $A_k$ for some $k \in R$.  It
follows that $\rho(C) \in W$, and $f(\rho(C)) \in B_k$ for some $k \in
R$.  Since $f$~preserves the action of~$G$, $f(\rho(C)) = \rho(f(C)) =
\rho(x)$.  We have therefore shown that $\rho(\bigcup_{k \in L}B_k)
\subseteq \bigcup_{k \in R}B_k$; the reverse inclusion is proved the
same way, so the sets~$B_k$ satisfy this congruence.  This completes the
proof that (I) implies (II).

Now, suppose (II) holds; let the sets~$A_k$ be as in (II), and let $G$
be the free group of rotations.  We will show that (III) holds. Choose
any point~$x_0$ in one of the sets~$A_k$.  The subgroup of~$G$
consisting of those elements which fix~$x_0$ is abelian (since any two
rotations around the same axis commute) and therefore cyclic (the subgroup
must be free, because $G$~is a
free group).  Let $w$ be a generator of this subgroup.

If $w$~is the identity, then it is not a proper power of another element
of~$G$, since $G$~has no nonidentity elements of finite order.  If $w$~is
not the identity, then $w$~still cannot be a proper power of another
element~$v$ of~$G$, since then $v$~would have to be a rotation around
the same axis as~$w$ and would therefore also fix~$x_0$, contradicting
the fact that only powers of~$w$ fix~$x_0$.

If $g$ and~$g'$ are elements of~$G$, then $g(x_0) = g'(x_0)$ iff
$g^{-1}\circ g'$ fixes~$x_0$, iff $g^{-1}\circ g' \in \langle w
\rangle$, iff $g$ and~$g'$ are in the same left coset of $\langle w
\rangle$.  Therefore, we can define a one-to-one map~$\isom$ from
$G/\langle w \rangle$ to the $G$-orbit of $x_0$ by $\isom(g\langle w
\rangle) = g(x_0)$.  Clearly the map~$\isom$ preserves the action of~$G$
(that is, if $C \in G/\langle w\rangle$ and $g \in G$, then $\isom(gC)
= g\isom(C)$).  Now let $B_k = \isom^{-1}(A_k)$ for each~$k$; the
sets~$B_k$ are finite, disjoint, and not all empty (one of them contains
$\isom^{-1}(x_0)$), and any elements of~$G$ which witness congruences
between sets~$A_k$ will witness the same congruences between sets~$B_k$.
Therefore, (III) holds.

Finally, suppose (III) holds; we will prove (II).  Since the given
subsets of $F/\langle w \rangle$ are finite, and there are only finitely
many congruences to be witnessed, we may assume that $F$~is a free group
on finitely many generators.  Hence, there is a group~$G$ of rotations
of~$S^2$ which is isomorphic to~$F$; we may assume $F = G$.

We will now find a point~$x_0$ of~$S^2$ which is a fixed point of~$w$,
but is not a fixed point of any element of~$G$ which is not a power
of~$w$.  If $w$~is the identity of~$G$, then we can take $x_0$ to be any
point other that the fixed points of the nonidentity elements of~$G$ (of
which there are only countably many).  If $w$~is not the identity of~$G$,
let $x_0$ be one of the two fixed points of~$w$.  In this latter
case, since $w$~is not a proper power in~$G$, the facts at the beginning
of the proof of Theorem~3.1 imply that the only elements of~$G$
which commute with~$w$ are in $\langle w \rangle$.
So any other element of~$G$ must not have the same axis
as~$w$; in other words, no elements of~$G$ other than the powers of~$w$
fix~$x_0$.

Now define a map from~$G$ to the orbit of~$x_0$ by mapping $g \in G$ to
$g(x_0)$.  Then $g\in G$ and $g' \in G$ are mapped to the same
point in the orbit if and only if $g$ and~$g'$ are in the same left
coset of $\langle w \rangle$, since only elements of $\langle w
\rangle$ fix~$x_0$.  Therefore, we get an induced bijection from
$G/\langle w \rangle$ to the orbit; call this bijection~$\isom$.  Again
we easily see that $\isom$~preserves the action of~$G$.  If sets~$B_k$
are the given finite subsets of $G/\langle w\rangle$ satisfying the
congruences, and $A_k = \isom(B_k)$ for each~$k$, then the sets~$A_k$
satisfy the congruences, as witnessed by the same elements of~$G$ which
witness the congruences for the sets~$B_k$.  Therefore, (II) holds, as
desired. \QED\enddemo

Note that the restriction in (III) that $w$~is not a proper power is
necessary.  Without it, one could let $F=\Z$ (a free group on
one generator under addition) and $w = 3$, so that $F/\langle w\rangle$
has $3$~elements.  Then, letting $A_1,A_2,A_3$ be the three singleton
subsets of $F/\langle w\rangle$, one would get a solution to
the system $A_1 \cup A_2 \cong A_1 \cup A_3 \cong
A_2 \cup A_3$, while Theorem~3.1 states that (I) cannot hold
for this system.

\head 4. Finite subsets of free groups \endhead

Part~(III) of Theorem~3.2 suggests that it is useful to consider
satisfaction of system of congruences by finite sets in certain countable
spaces.  A particular case of special interest is when the word~$w$
is the identity element; here we are talking about finite subsets of
the free group~$F$ itself, under the canonical action of~$F$ on~$F$.
We may assume that $F$~is a free group on countably many, but at least
two, generators.  In this case, if we give~$F$ the discrete topology,
then $F$ acting on itself is actually a suitable space.

This space turns out to be universal for the problem of satisfying
congruences by finite nonempty sets, in the following sense:

\proclaim{Proposition 4.1} A system of congruences can be satisfied by finite
nonempty subsets of a free group if and only if it can be satisfied
by finite nonempty sets in every suitable space. \endproclaim

\demo{Proof}
The right-to-left implication is trivial because the free group on two
generators is itself a suitable space as above.

For the other direction, suppose we have finite nonempty subsets of the
free group~$F$ which satisfy the congruences.  Since only finitely many
generators of~$F$ are used for the elements of the nonempty subsets
and for the witnesses to the congruences, we may assume $F$~is finitely generated.  Hence, for any suitable space $(\X,G)$, $G$ has a subgroup~$F'$
isomorphic to~$F$; this means that there are finite nonempty
sets $A'_k \subseteq F'$ satisfying the congruences.

Since $G$~acts freely on a comeager subset of~$\X$, we can find
a point $x \in \X$ such that $G$~acts freely on the orbit of~$x$.
Let $A_k = \{g(x)\colon g \in A'_k\}$; then the sets~$A_k$
satisfy the congruences in~$\X$.
\QED\enddemo

The next result is quite easy for the case of the sphere (or any
other suitable space with a compatible metric which is invariant
under the group action), but requires a little more care in
the general case:

\proclaim{Proposition 4.2} If a system of congruences can be satisfied
by finite nonempty subsets of a free group, then it can be satisfied
by open nonempty sets in every suitable space. \endproclaim

\demo{Proof}
As in the preceding proof, find a point~$x$ in a free orbit
of the suitable space $(\X,G)$
such that there are nonempty finite subsets of the orbit of~$x$
satisfying the congruences.

If we have a metric for~$\X$ which is invariant under the group
action, then we can just replace the finitely many points with
open balls of the same radius, chosen so small that the balls
do not overlap.

Without assuming such a metric, we can proceed as follows.
Since $\X$~is Hausdorff and $G$~acts by continuous maps,
we have that, for any $g,g' \in G$ such that $g(x) \ne g'(x)$
(this holds for any distinct $g$ and $g'$, because the action is
free on this orbit), there is $\eps > 0$ so small that,
if $U$~is the open ball of radius~$\eps$ centered
at~$x$, then $g(U) \cap g'(U) = \nullset$.
Find such an~$\eps$ so small that it works for any two
of the finitely many group elements~$g$ for which $g(x)$ was
used in the above finite sets.  Then define $U$ as above;
if we replace $g(x)$ with~$g(U)$ for each of these group elements~$g$,
then we get open nonempty subsets satisfying the congruences.
\QED\enddemo

One could hope at this point that a system of congruences satisfied
by finite nonempty subsets of a free group would be satisfied
by open nonempty subsets with dense union in any suitable space.
Unfortunately, this is not the case; the trivial system $A_1 \cong A_2$
is a counterexample.  A less trivial counterexample (one which is
weak) is given in section~5.

However, it turns out that, if one can get finite subsets of a free
group satisfying a system of congruences and meeting two minor extra
restrictions given below, then one can get open subsets of~$S^2$ or any
other suitable space satisfying the system and having dense union.

Suppose the group $F$~is freely generated by $f_i$, $1 \le i \le m$.
For any $g,g' \in F$, there is a unique shortest path from~$g$ to~$g'$ via
the generators and their inverses (i.e., a sequence $g_0,g_1,\dots,g_k$
where $g_0 = g$, $g_k = g'$, and each $g_{j+1}$ is obtained from~$g_j$
by applying a single $f_i$ or~$f_i^{-1}$ on the left).  Call a subset~$S$
of~$F$ {\it connected\/} if, for all $g,g' \in S$, all of the group
elements along the shortest path from~$g$ to~$g'$ are also in~$S$.

If $S$~is a finite subset of~$F$, $g \in S$, and $i \le m$, then there
is a greatest $k \ge 0$ such that $f_i^j \circ g \in S$ for $0 \le j \le
k$, and there is a smallest $k' \le 0$ such that $f_i^j \circ g \in S$
for $k' \le j \le 0$.  The subset $\{f_i^j\circ g \colon k' \le\nobreak j
\le\nobreak
k\}$ of~$S$ is a maximal `line in the $i$\snug-direction' within~$S$; these
lines form a partition of~$S$.  We will say that such a set~$S$ is {\it
prime} if, for each~$i$, the cardinalities of the lines in the
$i$\snug-direction for~$S$ have no common factor greater than~$1$.  An
equivalent form of this definition can be stated as follows: $S$~is
prime if there do not exist $k \ge 2$, $i \le m$, and a set~$T$ such
that $S$~is the disjoint union of the sets $f_i^j(T)$ for $0 \le j < k$.

\proclaim{Theorem 4.3} Suppose a system of congruences has the following
property: there is a free group~$F$ on~$m$ generators such that there
are disjoint nonempty finite subsets of~$F$ satisfying the congruences,
and the union of these finite sets is connected and prime.  Then, for
any suitable space $(\X,G)$, there are nonempty pairwise disjoint open
subsets of~$\X$ with dense union which satisfy the congruences. \endproclaim

\demo{Proof} Let $f_i$ ($1 \le i \le m$) be free generators for a free
subgroup of~$G$; we may assume that $G$~is the group generated by the
elements~$f_i$, and that $F = G$.  Let $P$ be a finite nonempty prime
subset of~$G$.  We will show that there exist nonempty
pairwise disjoint open
subsets $A_p$ ($p \in P$) of~$\X$ with union dense in~$\X$ such that,
for any $p,q \in P$, if $q = f_i \circ p$, then $A_q = f_i(A_p)$.

To see that this suffices to prove the theorem, proceed as follows.
Suppose we have nonempty disjoint finite subsets~$C_k$ of~$F$
satisfying the congruences, and their union~$P$ is connected and
prime.  Construct open sets~$A_p$ as above.  Now let $A'_k = \bigcup_{p
\in C_k} A_p$ for each~$k$; we will see that the sets~$A'_k$ satisfy
the given congruences. (They are clearly nonempty
pairwise disjoint open sets
with union dense in~$\X$.) Suppose that $\bigcup_{k \in L} A'_k \cong
\bigcup_{k \in R} A'_k$ is one of the congruences in the system, and
let $h$ be an element of~$F$ such that $h(\bigcup_{k \in L} C_k) =
\bigcup_{k \in R} C_k$.  Write~$h$ as a reduced word in the generators~$f_i$,
say $h = \rho_n\rho_{n-1}\dots\rho_1$ where each~$\rho_j$ is
either~$f_i$ or~$f_i^{-1}$ for some~$i$.  Now, let $p \in \bigcup_{k
\in L} C_k$ be arbitrary.  We have $hp \in \bigcup_{k \in R} C_k$;
since $P$~is connected, all of the intermediate points $p_j =
\rho_j\rho_{j-1}\dots\rho_1p$ are in~$P$.  By the construction of the
sets~$A_q$ for $q \in P$, we have $A_{p_j} = \rho_j(A_{p_{j-1}})$ for
$1 \le j \le n$; putting these together gives $A_{hp} = h(A_p)$.  Since
$h$~was arbitrary, this gives $h(\bigcup_{k \in L} A'_k) \subseteq
\bigcup_{k \in R} A'_k$.  The same argument with $h^{-1}$ instead of~$h$
gives the reverse inclusion, so $h(\bigcup_{k \in L} A'_k) =
\bigcup_{k \in R} A'_k$.  Therefore, the sets~$A'_k$ satisfy the given
congruences.

The construction of the sets~$A_p$ will be a modification of the
construction in Theorem~2.1 of Dougherty~\cite{\Dougherty}.  Let $r
= |P|$.  We will use~$P$ as an index set instead of $\{1,2,\dots,r\}$.
(Let us fix a listing $P = \{p_1,p_2,\dots,p_r\}$, although it will not
be used much.) In order to construct the sets~$A_p$, we will construct
open sets~$B_p$ for $p \in P$ such that: $\bigcap_{p \in P} B_p =
\nullset$; the sets $\bigcap_{p'\in P:\;p' \ne p}B_{p'}$ for $p \in P$
are all nonempty, and their union is dense in $\X$; and, if $p,q \in P$
and $q = f_i \circ p$, then $f_i(\bigcap_{p'\in P:\;p' \ne p}B_{p'}) =
\bigcap_{p'\in P:\;p' \ne q}B_{p'}$.  Given such sets~$B_p$, the sets
$A_p = \bigcap_{p'\in P:\;p' \ne p}B_{p'}$ have the desired properties.

If $\rho$~is one of the generators~$f_i$ or one of the inverse
generators~$f_i^{-1}$, let
$E(\rho) = \{p \in P \colon \rho \circ p \notin\nobreak P\}$.  If
one views~$P$ as a graph (with an edge joining $p$ to~$q$ if $p = f_i \circ
q$ or $q = f_i \circ p$ for some $i$), then $E(\rho)$ can be thought of as
the `ends of~$P$ in the direction of~$\rho$.'

The sets~$B_p$ will be constructed as increasing unions of sets~$B_p^n$,
$n = 0,1,2,\dotsc$.  The sets~$B_p^n$ will satisfy the following
properties, which will be maintained as induction hypotheses:
\roster
\item[2] $\bigcap_{p \in P} B_p^n = \nullset$.
\item If $p,q \in P$, $i \le m$, and $q = f_i \circ p$, then
$f_i(B_p^n) = B_q^n$; also, for each $i \le m$, $f_i(\bigcap_{p\in
E(f_i)} B_p^n) = \bigcap_{p\in E(f_i^{-1})} B_p^n$.
\item For any $x\in \X$, the set of $y\in \X$ which are connected to~$x$
by a chain of active links is finite.
\endroster
(There is no property~\ri1.)
The definitions of `link' and `active link' are the same as they were
in the proof of Theorem~2.1 of Dougherty~\cite{\Dougherty}:

\procl{Definition} Two points $x$ and~$x'$ are {\it linked}, or there is
a {\it link} from~$x$ to~$x'$, if $x' = f_i(x)$ or $x = f_i(x')$ for
some $i \le m$.  Points $x$ and~$x'$ are {\it connected by a chain of
links} if there are points $x_0,x_1,\dots,x_J$ with $x_0 = x$ and $x_J =
x'$ such that there is a link from $x_{j-1}$ to $x_j$ for each $j \le
J$. A link from~$x$ to~$x'$ is {\it active\/} (for the sets~$B_p^n$) if
there is a point in one or more of the sets~$B_p^n$ which is connected
to~$x$ or to~$x'$ by a chain of at most $2^r$~links.

\smallskip
Note that adding one new point to a set~$B_p^{n+1}$ activates only a
finite number of new links, although the finite number is very large.

We will construct sets~$B_p^n$ (increasing with~$n$) with the above
properties so that, if $B_p = \bigcup_{n=0}^\infty B_p^n$, then the sets
$\bigcap_{p'\in P:\;p' \ne p}B_{p'}$ are nonempty and have dense union.
Given this,
we clearly have $\bigcap_{p \in P} B_p =\nobreak \nullset$, by~\ri2.  Now
suppose $p,q \in P$ and $q = f_i \circ p$; then
$$\alignat2
f_i\biggl(\bigcap_{p'\in P:\;p' \ne p}B_{p'}\biggr)
&= f_i\biggl(\bigcap_{p' \in E(f_i)}B_{p'} \,\,\cap\,\,
  \bigcap_{p' \in P \setminus E(f_i):\; p' \ne p}B_{p'}\biggr)
&&\\
&= f_i\biggl(\bigcap_{p' \in E(f_i)}B_{p'}\biggr) \,\,\cap\,\,
  \bigcap_{p' \in P \setminus E(f_i):\; p' \ne p}f_i(B_{p'})
&&\qquad\text{since $f_i$ is one-to-one}\\\allowdisplaybreak
&= \bigcap_{p' \in E(f_i^{-1})}B_{p'} \,\,\cap\,\,
  \bigcap_{p' \in P \setminus E(f_i):\; p' \ne p}B_{f_i \circ p'}
&&\qquad\text{by \ri3}\\
&= \bigcap_{p' \in E(f_i^{-1})}B_{p'} \,\,\cap\,\,
  \bigcap_{p' \in P \setminus E(f_i^{-1}):\; p' \ne q}B_{p'}
&&\\
&= \bigcap_{p'\in P:\;p' \ne q}B_{p'}.
&&
\endalignat$$
Therefore, the sets~$B_p$ have all of the required properties.

Let $B_p^0 = \nullset$ for all~$p$.  Fix a listing
$\langle Z_n \colon n = 0,1,2,\dotsc \rangle$
of the nonempty sets in some
countable base for~$\X$, making sure that $\X$~itself is listed at least
$r$~times; the $t$\snug'th time we reach the set~$\X$ in the list ($t \le r$),
we will ensure that $\bigcap_{p'\in P:\;p' \ne p_t} B_{p'}$ is nonempty.

So suppose we are given $B_p^n$ ($p \in P$) and $Z = Z_n$.  Let $Z'$ be~$Z$
unless $Z$ is~$\X$ for the $t$\snug'th time ($t \le r$), in which case
let $Z'$ be the interior of the complement of~$B_{p_t}^n$.
(This~$Z'$ must be nonempty, because $B_{p_t}^n$~cannot be dense.
If $B_{p_t}^n$~were open dense, then \ri3~would imply that all of the
sets~$B_p^n$ were open dense, since $P$~is connected; this would
contradict~\ri2.)
Let $D$ be the complement of a ($G$\snug-invariant) comeager set on which
$G$~acts freely, and let $D'$ be the union of the images under the elements
of~$G$ of the boundaries of the sets~$B_p^n$; then $D \cup D'$ is
meager.  Let $x_0$ be any point in $Z' \setminus (D \cup D')$.  By~\ri2,
we can
choose $\pp \in
P$ such that $x_0 \notin B_\pp^n$ (making sure to set $\pp = p_t$ if
$Z$ is~$\X$ for the $t$'th time).  We will enlarge the sets~$B_p^n$ to
sets~$B_p^{n+1}$ so that $x_0 \in B_p^{n+1}$ for all~$p$ other than~$\pp$.

First, we will define~$\hat B_p$ to be $B_p^n \cup \{g(x_0) \colon g
\in\nobreak
T_p\}$ for some $T_p \subseteq G$.  To define~$T_p$, we will give a
recursive definition (based on the reduced form of elements of~$G$) of
a set $M_g \subseteq P$ for each $g \in G$, and then let
$T_p = \{g \in G \colon p \in\nobreak M_g\}$.

If $g$~is the identity of~$G$,
let $M_g = \{p \in P\colon p \ne\nobreak \pp\}$.  Otherwise, we can write~$g$
uniquely as $\rho \circ g'$ where $\rho = f_i^{\pm 1}$ and $g'$~has
a shorter reduced form than $g$~does, and hence $M_{g'}$~is already
defined.  Let $M_{g'}^+ = M_{g'} \cup \{p\colon g'(x_0) \in\nobreak B_p^n\}$.
If $M_{g'} = \nullset$, let $M_g = \nullset$.  Otherwise, let $M_g =
\{\rho p \colon p \in\nobreak
M_{g'}^+ ,\, \rho p \in\nobreak P\} \cup E'$,
where $E'$ is $E(\rho^{-1})$ if $E(\rho) \subseteq M_{g'}^+$,
$\nullset$ otherwise.

The first task is to show by induction on $g \in G$ that
$M_g^+ \ne P$ for all~$g$.  If $g$~is the identity, then $\pp \notin
M_g^+$ by the definition of~$x_0$.  Otherwise, write~$g$ as $\rho \circ
g'$ as above.  If $M_{g'} = \nullset$, then $M_g^+ = \{p\colon g(x_0)
\in\nobreak
B_p^n\} \ne P$ by~\ri2.  Now suppose $M_{g'} \ne \nullset$.  By the
induction hypothesis, choose $q \in P$ such that $q \notin M_{g'}^+$;
in particular, $g'(x_0) \notin B_q^n$.  If $q \in E(\rho)$, then the
definition of~$M_g$ gives $E(\rho^{-1}) \cap M_g = \nullset$, and we
cannot have $g(x_0) \in \bigcap_{p \in E(\rho^{-1})} B_p^n$ because
this and~\ri3 would give $g'(x_0) \in \bigcap_{p \in E(\rho)} B_p^n$,
contradicting $g'(x_0) \notin B_q^n$.  Hence, we cannot have
$E(\rho^{-1}) \subseteq M_g^+$, so $M_g^+ \ne P$.  Finally, suppose $q
\notin E(\rho)$.  Then $g'(x_0) \notin B_q^n$ and \ri3 imply $g(x_0)
\notin B_{\rho q}^n$, while the definition of~$M_g$ gives $\rho
q \notin M_g$, so $\rho q \notin M_g^+$, so $M_g^+ \ne P$.

We now check that, if $g$ and~$g'$ are in~$G$, $\rho$ is~$f_i^{\pm 1}$
for some~$i$, and $g = \rho \circ g'$, then $E(\rho) \subseteq
M_{g'}^+$ if and only if $E(\rho^{-1}) \subseteq M_g^+$; also, for any
$p \in P$ which is not in $E(\rho)$, $p \in M_{g'}^+$ if and only if
$\rho p \in M_g^+$.  We may assume that the reduced form of~$g'$
does not have~$\rho^{-1}$ as its leftmost component (otherwise,
interchange $g$ and~$g'$ and replace~$\rho$ with~$\rho^{-1}$); hence,
$M_g$~is defined from~$M_{g'}$ as above.  If $M_{g'} = \nullset$ and
hence $M_g = \nullset$, then the desired equivalences follow
immediately from~\ri3, so suppose $M_{g'} \ne \nullset$.  The
left-to-right implications are now immediate from the definition of~$M_g$.
For the first right-to-left implication, if $E(\rho)
\not\subseteq M_{g'}^+$, then $E(\rho^{-1}) \cap M_g = \nullset$ by
definition of~$M_g$, while $E(\rho^{-1})\not\subseteq \{p\colon g(x_0)
\in\nobreak B_p^n\}$ because otherwise \ri3~would give $E(\rho) \subseteq
\{p\colon g'(x_0) \in\nobreak B_p^n\} \subseteq M_{g'}^+$, so
$E(\rho^{-1})\not\subseteq M_g^+$.  The second right-to-left implication
is proved in the same way.

We are now ready to prove \ri2--\ri4 for the sets~$\hat B_p$.  The
definitions of $T_p$ and~$\hat B_p$ (and the fact that $G$~acts freely
on the orbit of~$x_0$) easily imply that $\{p \colon g(x_0) \in\nobreak \hat
B_p\} = M^+_g$ for all $g \in G$, while $\{p \colon x \in\nobreak \hat B_p\} =
\{p \colon x \in\nobreak B_p^n\}$ if $x$~is not in the $G$\snug-orbit of~$x_0$.
Therefore, properties \ri2 and~\ri3 for~$\hat B_p$ follow from the same
properties for~$B_p^n$ and the above facts about~$M^+_g$.
For~\ri4, we need some additional Claims.

We first note some useful facts about the sets $E(\rho)$.  We have
$|E(\rho)| = |E(\rho^{-1})|$, because $\rho$~gives a bijection between
$P \setminus E(\rho)$ and $P \setminus E(\rho^{-1})$.  If we view~$P$
as a graph as explained earlier (put an edge between $p$ and $f_i \circ
p$ if these are both in~$P$), then the number of edges `in the
$i$\snug-direction' (i.e., coming from generator~$f_i$ as above) is precisely
$|P \setminus E(f_i)|$.  This graph on~$P$ cannot have any cycles,
because $G$~is a free group (a cycle in the graph would give a
nontrivial reduced word~$w$ and an element~$p$ of~$P$ such that $w
\circ p = p$, so $w$~would be the identity in~$G$).  Therefore, by
standard results in graph theory, the graph must have fewer edges than
vertices; that is, $\sum_{i = 1}^m |P \setminus E(f_i)| < r$.  In
particular, if $i,i' \le m$ are distinct, then $|P \setminus E(f_i)| +
|P \setminus E(f_{i'})| < r$, so $|E(f_i)| + |E(f_{i'})| > r$.

Now, define the labeled directed graph~$\setgraph$
as follows.  The vertices of~$\setgraph$ are the nonempty proper subsets
of~$P$.  Let $\rho$ be $f_i$ or~$f_i^{-1}$ for some~$i$,
and let $S$ be a proper subset of~$P$.
Let~$S'$ be $\{\rho p\colon p \in\nobreak S \setminus E(\rho)\} \cup E'$,
where $E'$ is $E(\rho^{-1})$ if $E(\rho) \subseteq S$,
$\nullset$~otherwise.  If $S' \ne \nullset$, then $\setgraph$~has
an edge from~$S$ to~$S'$ labeled~$\rho$.  This edge is called
{\sl good} if $E(\rho) \subseteq S$ or $E(\rho) \cap S = \nullset$
(in which case there is a corresponding edge from~$S'$ to~$S$
labeled~$\rho^{-1}$), {\sl bad} otherwise.

\procl{Claim 1} No cycle in~$\setgraph$ contains a bad edge.

\procl{Proof}
If there is a good edge from~$S$ to~$S'$, then $|S| = |S'|$;
if there is a bad edge from~$S$ to~$S'$, then $|S| > |S'|$.
Hence, if $S$~were a vertex in a cycle containing a bad edge,
we would get $|S|>|S|$.
\QEd\par\smallskip

Now construct the undirected graph~$\setgraph_0$ by treating each pair
of oppositely-directed good edges in~$\setgraph$ as a single
undirected edge.

\procl{Claim 2} The undirected graph~$\setgraph_0$ is acyclic.

\procl{Proof}
Suppose we have a nontrivial cycle in~$\setgraph_0$; by taking
a minimal such cycle, we may ensure that there are no repeated
edges in the cycle.  This cycle corresponds to a cycle~$c$
in~$\setgraph$ (using good edges only)
which does not use both edges of a good pair
consecutively.  The vertices of~$c$ are subsets of~$P$ of the same size.
Call these vertices $N_0,N_1,\dots,N_{J-1}$, and let $e_j$~be
the edge from~$N_j$ to~$N_{j+1}$ (letting $N_J=N_0$), with
label~$\rho_j$, where $\rho_j$ is $f_{i_j}$ or~$f_{i_j}^{-1}$.

We now show that there must be some~$j$, $0 \le j < J$, such that
$E(\rho_j) \subseteq N_j$.  Suppose this is not the case; then we
simply have $N_{j+1} = \{\rho_j p\colon p \in\nobreak N_j\}$ for each
such~$j$.  Now start with some $p \in N_0$, and get $\rho_0p \in N_1$,
$\rho_1\rho_0p \in N_2$, and so on; eventually we get
$gp \in N_J = N_0$, where $g = \rho_{J-1}\rho_{J-2}\dots\rho_0$.  Note
that this expression for~$g$ is in reduced form, since the assumptions
above forbid $\rho_{j+1} = \rho_j^{-1}$; hence, $g$~is not the
identity.  We can now repeat this process to get $g^2 p \in N_0$,
$g^3p \in N_0$, and so on forever; this gives infinitely many elements
of~$N_0$, contradicting the finiteness of~$P$.

The same argument can be applied to the sets $P \setminus N_j$ instead
of~$N_j$; hence, there must be some~$j$, $0 \le j < J$, such that
$E(\rho_j)$ and~$N_j$ are disjoint.

The next step is to show that the numbers~$i_j$ for $0 \le j < J$ must
all be the same.  Suppose this is not so.  Then there must be numbers
$j,j'$ such that $0 \le j,j' < J$, $E(\rho_j) \subseteq N_j$,
$E(\rho_{j'}) \cap N_{j'} = \nullset$, and $i_j \ne i_{j'}$.  (Choose
$j$ and~$j'$ such that $E(\rho_j) \subseteq N_j$ and $E(\rho_{j'}) \cap
N_{j'} = \nullset$.  If $i_j \ne i_{j'}$, we are done; otherwise, we
can find~$j''$ such that $i_j \ne i_{j''}$, and one of the two pairs
$j,j''$ or $j'',j'$ will work.)  Then previous results give
$|E(\rho_j)| + |E(\rho_{j'})| > r$.  But $E(\rho_j) \subseteq N_j$ and
$E(\rho_{j'}) \subseteq P \setminus N_{j'}$, so $|E(\rho_j)| \le |N_j|$
and $|E(\rho_{j'})| \le r - |N_{j'}|$, so $|N_j| + r - |N_{j'}| > r$,
so $|N_j| > |N_{j'}|$, which is impossible because we established
earlier that the sets $N_0,N_1,\dots,N_J$ have the same size.

So all of the numbers~$i_j$ for $0 \le j < J$ are the same; since
both edges of a good pair cannot appear consecutively,
the values~$\rho_j$ must be identical. From now on, we will just
write~$\rho$ for this common value, and $i$~for the common
value of~$i_j$.

Define the infinite sequence $\NN_j$, $j = 0,1,2,\dotsc$, by letting
$\NN_j = N_{j\,\bmod{\,J}}$.  This sequence is periodic with
period a divisor of $J$.  Also, for each~$j$, we have either
$E(\rho) \subseteq \NN_j$, in which case $\NN_{j+1} =
\{\rho p\colon p \in\nobreak \NN_j\setminus E(\rho)\} \cup E(\rho^{-1})$,
or $E(\rho) \cap \NN_j = \nullset$, in which case $\NN_{j+1} =
\{\rho p\colon p \in\nobreak \NN_j\}$.
Futhermore, $E(\rho) \subseteq \NN_j$ for infinitely many~$j$, and
$E(\rho) \cap \NN_j = \nullset$ for infinitely many~$j$.

Fix an element~$p$ of $E(\rho^{-1})$.  The sequence $p,\rho p,\rho^2p,
\dotsc$ cannot lie entirely within~$P$, so there is a least $k > 0$
such that $\rho^k p \notin P$.  Then $\rho^{k-1}p \in E(\rho)$, and the
set $\{p,\rho p,\dots,\rho^{k-1}p\}$ (of size~$k$) is a `line in
the $i$\snug-direction' for~$P$, as defined in the paragraph preceding
this theorem.

Now, suppose $E(\rho) \subseteq \NN_j$; then $p \in \NN_{j+1}$,
$\rho p \in \NN_{j+2}$, $\rho^2 p \in \NN_{j+3}$, and so on, until
eventually we get $\rho^{k-1}p \in \NN_{j+k}$.  This means that $E(\rho)$
cannot be disjoint from $\NN_{j+k}$, so we must have
$E(\rho) \subseteq N(j+k)$.

We have just shown that, if $E(\rho) \subseteq \NN_j$ and $k$~is the
size of some line in the $i$\snug-direction for~$P$, then
$E(\rho) \subseteq \NN_{j+k}$.  Then, if $k'$~is also the size of
a line in the $i$\snug-direction for~$P$ (possibly the same line), then
$E(\rho) \subseteq \NN_{j+k+k'}$, and so on.  In fact, if $K$~is any
sum of nonnegative multiples of sizes of lines in the
$i$\snug-direction for~$P$,
then $E(\rho) \subseteq \NN_j$ implies $E(\rho) \subseteq \NN_{j+K}$.

Since $P$~is prime, the sizes of the lines in the $i_0$\snug-direction for~$P$
have no common divisor greater than~$1$.  Therefore, by standard number
theory, $1$~is a sum of multiples (not necessarily nonnegative) of these
sizes.  For each negative multiple~$ck$ occurring in this sum ($k$~a line
size, $c < 0$), replace~$ck$ with $(c - Jc)k$, which is
a nonnegative multiple of~$k$; this replacement will increase the sum by
a multiple of $J$.  The result is that we get a number
$K \equiv 1 \pmod{J}$ which is a sum of
nonnegative multiples of line sizes.
Therefore, for any~$j$ such that $E(\rho) \subseteq \NN_j$, we get
$E(\rho) \subseteq \NN_{j+K}$.  But the periodicity of
$\NN_0,\NN_1,\dotsc$ implies that $\NN_{j+K} = \NN_{j+1}$.
Therefore, if $E(\rho) \subseteq \NN_j$, then $E(\rho) \subseteq \NN_{j+1}$;
repeated application of this gives $E(\rho) \subseteq \NN_{j'}$ for
all $j' > j$.  This is the final contradiction, because there are
infinitely many~$j'$ such that $E(\rho) \cap \NN_{j'} = \nullset$;
hence, the claim is proved. \QEd\par\smallskip

%
%

The rest of the proof is just like the last part of the
proof of Theorem~2.1 of Dougherty~\cite{\Dougherty}.
Using the above two claims, we get:

\procl{Claim 3} Every path of length~$2^r$ in the
digraph~$\setgraph$ contains a pair of consecutive edges
with labels $f_i$ and~$f_i^{-1}$, or vice versa, for some~$i$.

\procl{Proof}  Suppose we have a path of length~$2^r$ in~$\setgraph$.
Since there are fewer than~$2^r$ vertices in~$\setgraph$, some
vertex must be visited more than once, so we get a nontrivial subpath
which starts and ends at the same vertex (i.e., a cycle).  By
Claim~1, this subpath consists entirely of good edges, so
it induces a corresponding path in the graph~$\setgraph_0$
which also starts and ends at the same place.  By Claim~2,
this latter path cannot be a nontrivial cycle, so it must
double back on itself (use the same edge twice in succession);
hence, the original path uses both edges of a pair of oppositely-directed
good edges successively, which gives the desired conclusion.
\QEd\par\smallskip

Now, for any $g \in G$, $x_0$~is connected to~$g(x_0)$ by a chain
of links, and this chain can be read off from the reduced form of~$g$.
In order to prove~\ri4 for the sets~$\hat B_p$, it will suffice to show
that, if $M_g \ne \nullset$, then either all of the links in this chain are
active for the sets~$B_p^n$, or the chain has fewer than~$2^r$ links; once
we know this, \ri4 for~$B_p^n$ implies that there are only finitely
many points~$g(x_0)$ such that $M_g \ne \nullset$ (equivalently, since
$G$~acts freely on the orbit of~$x_0$, the set of~$g$ such that
$M_g \ne \nullset$ is finite), so only finitely many new links are
activated when $B_p^n$~is enlarged to~$\hat B_p$, so \ri4 for~$B_p^n$
implies \ri4 for~$\hat B_p$.

So suppose $M_g \ne \nullset$ and the above chain has at least $2^r$~links.
Then $M_h \ne \nullset$ for all of the intermediate points~$h(x_0)$
on the chain. It must now be true that, given any $2^r$~consecutive
links in the chain, at least one of the $2^r+1$ endpoints of
these links is in one of the sets~$B_p^n$, because otherwise the sets~$M_h$
at these $2^r+1$ endpoints would give a counterexample to
Claim~3. (If none of these points~$h(x_0)$ is in any of the sets~$B_p^n$,
then we
always have $M^+_h = M_h$.  Now, if $h$ and $h' = \rho \circ h$ are
final subwords of the reduced word for~$g$,
where $\rho$~is $f_i$ or~$f_i^{-1}$,
then the way in which $M_{h'}$~is computed from~$M_h$ shows that
there is an edge in~$\setgraph$ from~$M_h$ to~$M_{h'}$ labeled~$\rho$.
The resulting path of length~$2^r$ cannot
include consecutive edges labeled $f_i$ and~$f_i^{-1}$ or vice versa
because we are working with the reduced form of~$g$.)
It follows that all~$2^r$
of the links are active for~$B_p^n$; since this was an arbitrary
subchain of the chain, all of the links in the chain are active for~$B_p^n$.
This completes the proof of~\ri4 for~$\hat B_p$.

Now that we have \ri2--\ri4 for~$\hat B_p$, we can enlarge these sets
to get open sets.  Let $S$ be the set of $g \in G$ such that $x_0$~is
connected to~$g(x_0)$ by a chain of links which are active for the
sets~$\hat B_p$, and let $S'$ be the set of $g'\in G$ such that $g'(x_0)$~is
connected to~$g(x_0)$ for some $g \in S$ by a chain of at most
$2^r+1$ links.  Then $T_p \subseteq S$ for all~$p$, $S \subseteq S'$,
and $S$ and~$S'$ are finite by~\ri4.  Let $U_0$ be an
open neighborhood of~$x_0$
so small that the images~$g(U_0)$ for $g \in S'$ are pairwise disjoint and
each of them is either included in or disjoint from each of the sets~$B_k^n$.
(This is possible because, by the choice of~$x_0$, no point
in~$S'$ is on the boundary of any of the sets~$B_p^n$.)  Now let
$B_p^{n+1} = B_p^n \cup \bigcup \{g(U_0) \colon g \in\nobreak T_p \}$
for each~$p$;
we must see that these sets satisfy properties \ri2--\ri4.

From the definition of~$B_p^{n+1}$ and the disjointness of the sets~$g(U_0)$
for $g \in S'$, the following two statements follow easily: If
$x\in g(U_0)$ for some $g\in S'$, then $x\in B_p^{n+1}$ if and only if
$g(x_0)\in \hat B_p$. If $x\in\X$ is not in any of the sets~$g(U_0)$ for
$g\in S$, then $x\in B_p^{n+1}$ if and only if $x\in B_p^n$.

We can now prove \ri2--\ri4 for~$B_k^{n+1}$.

\ri2: If a point~$x$ is in one of the neighborhoods~$g(U_0)$ where $g\in
S$, then $g(x_0)\notin \bigcap_{p \in P} \hat B_p$ by \ri2 for~$\hat B_p$,
so $x\notin \bigcap_{p \in P} B_p^{n+1}$; if $x$~is not in one of these
neighborhoods, then $x\notin \bigcap_{p \in P} B_p^n$ by \ri2 for~$B_p^n$,
so $x\notin \bigcap_{p \in P} B_p^{n+1}$.

\ri3: We prove $f_i(B_p^{n+1}) \subseteq B_q^{n+1}$
where $p,q \in P$ and $q = f_i \circ p$;
the other parts are similar.  Suppose $x \in
B_p^{n+1}$.  If $x \in g(U_0)$ for some $g \in S$,
then $g(x_0) \in \hat B_p$, so $f_i(g(x_0)) \in
\hat B_q$ by \ri3 for~$\hat B_p$; but $f_i \circ g
\in S'$ and $f_i(x) \in f_i(g(U_0))$, so $f_i(x) \in B_q^{n+1}$.
If $x$~is not in~$g(U_0)$ for any $g \in S$, then $x
\in B_p^n$, so $f_i(x) \in B_q^n$ by \ri3 for~$B_p^n$.

\ri4: Let $w$ be any point of~$\X$, and consider the set of all points
connected to~$w$ by a path of links which are active for the
sets~$B_p^{n+1}$. If this set contains no point which is in~$g(U_0)$ for any
$g\in S$, then all of the links connecting the set were in fact active
for $B_p^n$. (Note: If the link from~$x$ to~$x'$ is
activated by~$x''$, because there is a chain of at most~$2^r$ links
connecting~$x''$ to~$x$ or to~$x'$, then all of the links in this chain
are also activated by~$x''$.) Hence, the set is finite by \ri4
for~$B_p^n$.  So suppose $y\in g(U_0)$ is connected by active links to~$w$,
and $g\in S$.  A point is connected to~$w$ if and only if it is
connected to~$y$, so it will suffice to show that only finitely many
points are connected to~$y$.

Suppose $y$~is actively linked to~$y'$, say $y' = f_i(y)$ (the case
$y' = f_i^{-1}(y)$ is similar).  Let $y''$ be a point in one of
the sets $B_p^{n+1}$ such that $y''$~is connected to either
$y$ or~$y'$ by a chain of at most~$2^r$ links.  Then there is an
element~$h$ of~$G$ such that $h(y) = y''$, and the reduced form
of~$h$ in terms of the generators~$f_I$ has length at most $2^r+1$
(and, if it has length $2^r+1$, then the rightmost component is~$f_i$).
Therefore, $h \circ g \in S'$.  We now have $y'' \in h(g(U_0))$,
so, since $y'' \in B_p^{n+1}$, we must have $h(g(x_0)) \in \hat B_p$.
This means that the link from~$g(x_0)$ to $f_i(g(x_0))$ is active for
the sets~$\hat B_p$, so $f_i(y) \in f_i(g(x_0))$ and $f_i \circ g \in S$.

Now this argument can be repeated starting at~$y'$, and so on; the
result is that, for any chain of active (for the sets~$B_p^{n+1}$) links
starting at~$y$, all of the links in the corresponding chain starting
at~$g(x_0)$ are also active (for the sets~$\hat B_p$).  Furthermore, if
$y$~is connected to two different points $y'$ and~$y''$ by such chains of
links, this will give $y'=h'(y)$ and $y''=h''(y)$ for some distinct
elements $h',h''$ of~$G$, and the corresponding points reached
from~$g(x_0)$ will be $h'(g(x_0))$ and $h''(g(x_0))$; since $G$~acts freely
on the orbit of~$x_0$, these two points will also be different.
Therefore, since $g(x_0)$~is connected to only finitely many points, $y$
(and hence~$w$) must be connected to only finitely many points.  This
completes the proof of~\ri4 for the sets~$B_p^{n+1}$.

This completes the recursive construction. \QED\enddemo

Note that both `connected' and `prime' are needed here; neither one
suffices by itself.  The trivial non-weak system $A_1 \cong A_2$
is satisfied by nonempty finite subsets of a free group (singletons,
in fact), which can be placed next to each other so that their union is
connected; or they could be made non-adjacent, in which case their union
would be prime but not connected.  But the system cannot be satisfied
by open subsets of~$S^2$ with dense union using free rotations.

One could strengthen the definition of `prime' by considering `lines
in the $g$\snug-direction' for any group element~$g$, not just the
generators; call this version `strongly prime.'  Then it is a consequence
of Theorem~4.3 (and its proof) that any connected and prime finite subset
of a free group is strongly prime.  (Is there a simple direct proof?)
It might be that `strongly prime' would suffice for Theorem~4.3, without
connectedness being needed; but the proof would need substantial revision.

\head 5. More examples on the sphere \endhead

In this section, we give two examples which make use of the special
properties of the sphere~$S^2$.  The first example uses these properties
to show that certain sets do not exist, while the second uses these
properties to show that certain sets do exist.

First, look at the system
$$A_1 \cong A_2 \cong A_3,\qquad A_1 \cup A_2 \cong A_1 \cup A_3.$$
It is easy to
get finite subsets of the sphere satisfying these congruences via
free rotations: let $\sigma$ and~$\rho$ be two such rotations around
different axes, let $x$ be a fixed point of~$\sigma$, and let
$$A_1 = \{x\},\qquad A_2 = \{\sigma(\rho(x))\}, \qquad A_3 = \{\rho(x)\};$$
then $\rho(A_1) = A_3$, $\sigma(A_3) = A_2$, and $\sigma(A_1 \cup A_3)
= A_1 \cup A_2$.
As in Theorem~3.2, we can enlarge these points to open disks
to get open sets satisfying the congruences via free rotations.

Also, the above system is a subsystem of~$\UNC3$, so we know that
it is satisfied by open subsets of~$S^2$ with dense union
if arbitrary rotations are allowed.

But one cannot combine the above:

\proclaim{Theorem 5.1} The system of congruences $A_1 \cong A_2 \cong A_3$,
$A_1 \cup A_2 \cong A_1 \cup A_3$ cannot be satisfied by
open subsets of~$S^2$ with dense union using free rotations. \endproclaim

\demo{Proof}
Suppose $\sigma$, $\rho$, and~$\tau$ are members of a free group~$F$ of
rotations of~$S^2$ and $A_1$,~$A_2$, and~$A_3$ are pairwise disjoint
open subsets of $S^2$ such that $\sigma(A_1 \cup A_3) = A_1 \cup A_2$,
$\tau(A_2) = A_3$, and $\rho(A_1) = A_3$.  We will show that
$A_1 \cup A_2 \cup A_3$ cannot be dense in~$S^2$.

As we saw in section~3, $\sigma$~must map connected components of $A_1$
or~$A_3$ to connected components of $A_1$ or~$A_2$. Of course, $\rho$
maps connected components of~$A_1$ to connected components of~$A_3$, and
similarly for~$\tau$. Hence, one can form a labeled directed
graph~$\setgraph$ whose vertices are the connected components of the
sets $A_1,A_2,A_3$ (labeled $1,2,3$, respectively), and whose edges are
given as follows: if $C$~is a component of $A_1$ or~$A_3$, then there is
an edge labeled~$\sigma$ from~$C$ to~$\sigma(C)$; if $C$~is a component
of $A_2$, then there is an edge labeled~$\tau$ from~$C$ to~$\tau(C)$; if
$C$~is a component of $A_1$, then there is an edge labeled~$\rho$
from~$C$ to~$\rho(C)$.

This digraph is related to a much larger digraph~$\Fgraph$, whose
vertices are {\sl all} nonempty open connected proper subsets of the
sphere, with edges given by: if $g$~is one of the generators of~$F$
(from a free generator set fixed in advance) and $C$~is a vertex
of~$\Fgraph$, then there is an edge of~$\Fgraph$ from~$C$ to~$g(C)$. So
the connected component of the vertex~$C$ in~$\Fgraph$ is just the orbit
of~$C$ under~$F$.  Each edge of~$\setgraph$ corresponds to a finite `path'
in~$\Fgraph$, given by the expression of the label of that edge
($\sigma$,~$\tau$, or~$\rho$) as a reduced word in the generators
of~$F$ (where an occurrence of an inverse generator means that an
$\Fgraph$\snug-edge is to be traversed backward); hence,
(the vertex set of) each component of~$\setgraph$ is
included in a component of~$\Fgraph$.

The structure of each component of~$\Fgraph$ is rather simple.
If $F$~acts freely on the vertices of the component, then the component
looks just like the Cayley graph of~$F$ --- a tree with $m$~edges leading
from each vertex and $m$~edges leading to each vertex (where $m$~is the
number of generators of~$F$), and no cycles even if the orientation of
edges is ignored.  If $F$~does not act freely on the vertices, let $v$~be
a vertex fixed by some nontrivial element of~$F$.  The elements of~$F$
which fix~$v$ form a subgroup of~$F$ which is abelian (these elements all
have to be rotations around the same axis, by the remark after Lemma~2.1)
and hence cyclic; let $w$~be a generator of this subgroup.  We may assume
$w$ is cyclically reduced (if $w = g^{-1}w'g$, then $w'$ generates the
subgroup fixing the vertex $g(v)$, so we can use that instead).  So the
word~$w$ describes a `cycle' in~$\Fgraph$ (where $w$~is read from right
to left, and inverse generators in~$w$ mean that edges of~$\Fgraph$
are traversed backward); since no elements of~$F$ other than powers
of~$w$ fix~$v$, there are no other cycles in this component of~$\Fgraph$
(ignoring the direction of edges), so the component resembles a single
ring (which we will call the {\it prime ring} of the component) with
copies of parts of the free $F$-tree attached to each vertex.  By Lemma~2.1,
the word~$w$ cannot be a proper power in~$F$ (if $w$~is a power of~$w'$,
then $w'$ is also a rotation around the same axis as~$w$, so $w'$ must
also fix~$v$, so $w'$ is in the subgroup generated by~$w$); this is
why the ring is called `prime.'

The components of~$\setgraph$ are rather different.  First, the vertices
of a component of~$\setgraph$ are pairwise disjoint subsets of~$S^2$
of the same positive measure, so the component must be finite.
Second, it is possible for a component of~$\setgraph$ to include
multiple cycles, if the words $\rho$, $\sigma$, and $\tau$
satisfy nontrivial relations in~$F$.

From any vertex of~$\setgraph$, we can follow $\sigma$\snug-edges
forward either forever or until we reach a 2-vertex, and we can follow
$\sigma$\snug-edges backward either forever or until we reach a 3-vertex.
(Any 1-vertex or 3-vertex has a unique $\sigma$\snug-edge leading from
it, and any 1-vertex or 2-vertex has a unique $\sigma$\snug-edge leading
to it.)  These vertices and edges form the {\it $\sigma$\snug-path}
containing the given vertex.  This path is included in a component
of~$\setgraph$, so it cannot contain infinitely many vertices, so it must
be a terminating path or a finite cycle.  If it is a terminating path,
it has the form $3 \to 1 \to 1 \to \dots \to 1 \to 2$ (0~or more~1's);
if it is a finite cycle, it must consist entirely of 1-vertices.  But then
Lemma~2.1 implies that the cycle must have length~1, because any component
fixed by~$\sigma^k$ ($k \ge 1$) must in fact be fixed by~$\sigma$.
Call a $\sigma$\snug-path of the latter type a {\it $1$\snug-loop}.

For each terminating $\sigma$\snug-path, there is a $\tau$\snug-edge
connecting the final vertex of the path to the initial vertex of
another (or perhaps the same) terminating $\sigma$\snug-path.
Again, since the components of~$\setgraph$ are finite, if one follows
these $\sigma$\snug- and $\tau$\snug-edges, one must eventually
repeat a vertex.  So these $\sigma$\snug-paths are joined together
into {\it $\sigma,\tau$-cycles}; each vertex of~$\setgraph$ which is
not (the unique vertex of) a 1-loop is in a unique $\sigma,\tau$-cycle.

Hence, for every component of one of the sets $A_1,A_2,A_3$, there is
a nontrivial word in $\sigma$ and~$\tau$ which fixes that component.
If this word does not collapse to the identity element of~$F$ (when
expressed in terms of the generators of~$F$), then it is a rotation of
infinite order around some axis~$\axis$, so, by Lemma~2.1, the component
is invariant under all rotations around~$\axis$, and hence must be a disk or
an annulus.  If we want to cover a dense part of the sphere with such
components, we will have to use a wide variety of them:

\proclaim{Lemma 5.2} Suppose that we have a collection of pairwise
disjoint nonempty connected open subsets of the sphere~$S^2$, with
union dense in~$S^2$.  Suppose that each of these subsets is completely
symmetric around some axis (so it is a disk or annulus around that axis).
Then either the sets are all symmetric around the same axis, or
infinitely many different axes are used.
\endproclaim

\demo{Proof}
Suppose the sets are not all symmetric around the same axis,
but only have finitely many axes of symmetry.
Let $A$ be one of the sets in the collection, let $\axis$
be the axis of symmetry of~$A$, and let $x$ be
one of the two points of~$S^2$ on~$\axis$.
Let $y$ be a point in~$A$, and let~$r$ be the distance from~$x$
to~$y$; then all points of~$S^2$ at distance~$r$ from~$x$ are
in~$A$.  Let $z$ be a point in another member of the collection
with a different axis of symmetry, and let $R$ be the distance
from~$x$ to~$z$.  We may assume $R>r$ (if not, replace~$x$ with
the other point of intersection of $S^2$ and~$\axis$).

Let $R_0$ be the greatest number above~$r$ such that no point
of~$S^2$ at distance between~$r$ and~$R_0$ from~$x$ is in
a set with axis of symmetry other than~$\axis$.  Then we
have $r < R_0 < R$.

Let $C$ be the circle in~$S^2$ with center~$x$ and radius~$R_0$.
Then the sets in the collection which have axis of symmetry~$\axis$
cover (at least a dense part of) the points just inside~$C$
(those at distance between $r$ and~$R_0$ from~$x$).

For each axis of symmetry~$\axis'$ other than~$\axis$, the plane through
$\axis$ and~$\axis'$ meets~$C$ in two points.  There are only finitely
many axes~$\axis'$, and hence only finitely many such points; let $w$~be
a point of~$C$ which is not one of these points.  So no circle of rotation
around such an axis~$\axis'$ is tangent to~$C$ at~$w$.

Now, for each~$\axis'$, there is a positive number~$\eps_{\axis'}$
such that no point~$w'$ within distance~$\eps_{\axis'}$ of~$w$ can
lie in a member of the collection with axis of symmetry~$\axis'$,
because the circle obtained by rotating~$w'$ around~$\axis'$
crosses over~$C$ and hence meets the members of the collection
with axis of symmetry~$\axis$.  Let $\eps$ be the least of these
numbers~$\eps_{\axis'}$.  The neighborhood of~$w$ with radius~$\eps$
cannot meet the members of the collection with axis other than~$\axis$,
so the members of the collection with axis~$\axis$ must cover a
dense part of this neighborhood.  By symmetry around~$\axis$, these
members actually cover a dense part of all of the points at distance
between $R_0$ and~$R_0+\eps$ from~$x$.  But this contradicts the
maximality of~$R_0$, so we are done.
\QED\enddemo

The proof of the theorem now proceeds by cases.
         
{\bf Case 1}: $\sigma$ and $\tau$ commute.
      
So $\sigma$ and~$\tau$ are both powers of some $\alpha\in F$.
         
A component of~$A_3$ cannot be fixed under~$\sigma$, so it cannot be
fixed under~$\alpha$, so (by Lemma~2.1) it cannot be fixed under any
nonzero power of~$\alpha$.  But any such component is a vertex
in a $\sigma,\tau$-cycle, so there is a word in $\sigma$ and~$\tau$
(a positive word, not using inverses) which fixes the component;
when expressed as a power of~$\alpha$, this word
must come out to $\alpha^0$.
            
Hence, $\sigma$ and~$\tau$ are powers of~$\alpha$ with exponents of
opposite sign; we may assume $\sigma = \alpha^n$ and $\tau = \alpha^{-m}$
with $n,m > 0$. And the ratio of $\sigma$\snug-edges to $\tau$\snug-edges
in each $\sigma,\tau$-cycle is $m:n$, so the ratio of 1-vertices to
2-vertices to 3-vertices in the cycle is $m-n:n:n$.
            
{\bf Subcase 1a}: $m-n > n$.

Since each component of~$\setgraph$ is a disjoint union of
$\sigma,\tau$\snug-cycles and 1-loops, we get that such a component
will contain more 1-vertices than 3-vertices.  This is impossible, because
$\rho$~gives a bijection between the 1-vertices and the 3-vertices
in the component.
            
{\bf Subcase 1b}: $m-n = n$.

In this case, there are exactly as many 1-vertices as 3-vertices in
each $\sigma,\tau$\snug-cycle; since the 1-vertices and 3-vertices must
balance in each component of~$\setgraph$, there cannot be any 1-loops in
any component.  But $\tau = \sigma^{-2}$, so a $\sigma$\snug-path cannot
end in $1 \to 1 \to 2$ ($\tau$ would send the 2-vertex back to a 1-vertex)
or $3 \to 2$ ($\sigma^{-1} = \sigma\tau$ would send the 2-vertex back to
the 3-vertex, which is impossible because $\sigma\circ\tau$~must send any
2-vertex to a 1-vertex or a 2-vertex).  So the only possible form for a
$\sigma$\snug-path is $3 \to 1 \to 2$.  This means that $\sigma(A_3) =
A_1$ and $\sigma(A_1) = A_2$.

If $A_1 \cup A_2 \cup A_3$ were dense, then $A_2$~would differ from the
complement of $A_1 \cup A_3$ by a meager set, and $A_3$ would differ from
the complement of $A_1 \cup A_2$ by a meager set; it would follow that
$\sigma(A_2)$ differs from~$A_3$ by a meager set.  So $\sigma^3(A_3)$
differs from~$A_3$ by a meager set.  By Lemma~2.1, $\sigma(A_3)$~differs
from~$A_3$ by a meager set; this is impossible, because $\sigma(A_3)$
is~$A_1$, which is an open set disjoint from~$A_3$.  So $A_1 \cup A_2
\cup A_3$ must not be dense.
            
{\bf Subcase 1c}: $m-n < n$.

Suppose we have a component of~$\setgraph$ in which the
$\sigma,\tau$\snug-cycles contain a total of $cm$~$\sigma$\snug-edges
and $cn$~$\tau$\snug-edges.  Then these cycles contain $c(2n-m)$
more 3-vertices than 1-vertices; the component must contain $c(2n-m)$
1-loops to make up the deficit.  (So the component has a total of
$3cn$ vertices.)  But there is a fixed limit on the number of 1-loops
one can have in a single component of~$\Fgraph$ (and hence in a single
component of~$\setgraph$).  Let $k$ be the length of the word~$\sigma$
in terms of the generators of~$F$.  If the $\Fgraph$\snug-component has
no prime ring, or if its prime ring has length greater than~$k$, then
the component cannot contain any 1-loops at all.  If the prime ring has
length at most~$k$, then any 1-loop gives a path determined by~$\sigma$
(some edges may be traced backward) of length~$k$ which must pass around
the prime ring at least once; it is easy to see that there are at most
$k$ such paths in the component.  (Fix a vertex~$v$ on the prime ring;
the path determined by~$\sigma$ is known completely once we know how many
steps it takes to reach~$v$.)  Therefore, we must have $c(2n-m) \le k$,
so each component of~$\setgraph$ has at most $N = 3nk/(2n-m)$ vertices.

Each 1-loop is a vertex fixed by~$\sigma$.  If $v$ is any vertex
of~$\setgraph$, there is a 1-loop in the same $\setgraph$\snug-component,
which must be reachable from~$v$ by following a path of fewer than~$N$
$\setgraph$\snug-edges (forward or backward).  So there is a group
element~$w$ which is a word of length less than~$N$ in $\rho,\sigma,\tau$
such that $\sigma(w(v)) = w(v)$ and hence $(w^{-1}\sigma w)(v) = v$.

This gives a finite list of non-identity elements of~$F$
such that each component of the sets $A_1,A_2,A_3$ is fixed
under one of these elements, and hence under all rotations around
the axis of this element.  Note that at least two axes are used;
the 1-loops are fixed under~$\sigma$ and the 3-vertices are not.
Therefore, by Lemma~5.2, $A_1 \cup A_2 \cup A_3$ is not dense.
This completes Case~1.

{\bf Case 2}: $\sigma$ and $\tau$ do not commute.
      
So $\langle\sigma,\tau\rangle$ (the subgroup of~$F$ generated by $\sigma$
and~$\tau$) is a free group with free generators $\sigma,\tau$.
         
To handle this case, we will need the following lemma:

\proclaim{Lemma 5.3} If $F$ is a free group and $G$
and~$H$ are free subgroups of~$F$ of rank\/~$2$ such that $H = \alpha G
\alpha^{-1}$ for some~$\alpha$ which is not in~$G$, then $G \cap H$
has rank at most\/~$1$ (i.e., is cyclic). \endproclaim
            
\demo{Proof} This is basically a case of Proposition 3.4 from
Nickolas~\cite{\Nickolas}.  The statement of that proposition says
``rank $m > 2$,'' but the proof works also for $m = 2$ in the case
where the two conjugate subgroups are distinct.  To see that the
groups $G$ and~$H$ here are indeed distinct, we need to know that
the normalizer of $G$ in~$F$ is just $G$~itself.  Equivalently, if
$F'$~is a free group and $G$~is a normal subgroup of~$F'$ which is
free of rank~$2$, then $G = F'$.  This follows from Theorem~2.10 in
Magnus-Karrass-Solitar~\cite{\Magnus}.
\QED\enddemo

Throughout Case~2 there will be no need to distinguish between
$\sigma,\tau$\snug-cycles and 1-loops, so from now on the term
`$\sigma,\tau$\snug-cycle' will include 1-loops as a special case.

Consider an arbitrary component~$C_\setgraph$ of~$\setgraph$,
whose vertices lie within the component~$C_\Fgraph$ of~$\Fgraph$.
If $v$~is a vertex of~$C_\setgraph$, then $v$~lies within some
$\sigma,\tau$\snug-cycle, so there is a nontrivial word in $\sigma$
and~$\tau$ which fixes~$v$.  Since $\langle\sigma,\tau\rangle$ is free,
this word is not the identity in~$F$.  So $C_\Fgraph$ cannot be
a free $F$\snug-tree; it must have a prime ring.

We will show:
\roster
\item"$\bullet$" there is a fixed limit~$L_1$ on the distance in~$C_\Fgraph$
from an arbitrary vertex of~$C_\setgraph$ to the prime ring; and
\item"$\bullet$" there is a fixed limit~$L_2$ on the length of the
prime ring.
\endroster
Both of these limits depend only on $\rho$,~$\sigma$, and~$\tau$,
not on the particular components being considered.

Given this, for every vertex~$v$ of~$C_\setgraph$, there is a path
in~$\Fgraph$ (ignoring direction of edges, as usual) from this vertex
to the prime ring, around the ring, and back to the vertex, with total
length at most $2L_1+L_2$.  This gives a nontrivial word in~$F$ of length
at most $2L_1+L_2$ which fixes~$v$, and this word is a rotation whose
axis is an axis of complete symmetry of~$v$, by Lemma~2.1.  But there
are only finitely many such words, so there are only finitely many
axes of symmetry for the vertices of~$\setgraph$ (i.e., the components
of the sets $A_1,A_2,A_3$).  There must be at least two such axes,
though.  (Let $v$ be a 3-vertex.  Then there is a word~$w$ in $\sigma$
and~$\tau$, ending in~$\sigma$, which fixes~$v$, and this word does
use~$\tau$.  The vertex~$\sigma(v)$ is fixed by~$\sigma w\sigma^{-1}$.
Since $\sigma w\sigma^{-1}w$ ends in more~$\sigma$\snug's than
$w\sigma w\sigma^{-1}$ does, $w$ and~$\sigma w\sigma^{-1}$ do not commute
in~$\langle\sigma,\tau\rangle$, so they must be rotations with different
axes.)  Therefore, by Lemma~5.2, $A_1 \cup A_2 \cup A_3$ is not dense.
So we will be done with the proof of the theorem once we have shown that
the limits $L_1$ and~$L_2$ exist.

Let $k$ be the maximum of the lengths of $\sigma$ and~$\tau$
when written in terms of the generators of~$F$.  There are only finitely
many words in the generators of~$F$ of length at most~$2k$; since each
such word can be written in at most one way as a word in $\sigma$
and~$\tau$ (because $\langle\sigma,\tau\rangle$ is free), there are
only finitely many words in $\sigma$ and~$\tau$ which collapse in~$F$
to a word of length at most~$2k$.  Let $N$ be the greatest of the
lengths of these finitely many $\sigma,\tau$\snug-words (where here we
compute length by counting $\sigma$\snug's and~$\tau$\snug's,
not $F$\snug-generators).
In other words, if $w$~is a word in $\sigma$ and~$\tau$ of length greater
than~$N$, and $w'$~is the reduced form of the expression of~$w$ in terms
of the generators of~$F$, then $w'$~has length greater than~$2k$.
We will see that $(N+2)k$ is a
suitable value for the limit~$L_1$.

Let $v$ be a vertex of~$C_\setgraph$, and let $p$ be the
$\sigma,\tau$\snug-cycle it lies on.  Since $\sigma$ and~$\tau$ can
be written as words (of length at most~$k$) in the generators of~$F$,
the cycle~$p$ induces a path~$p'$ in~$C_\Fgraph$ which also starts
and ends at~$v$.  Some of the vertices of~$p'$ are the vertices of~$p$;
a $p$\snug-vertex occurs at least once every $k$~steps in~$p'$.
(We do not count as a `$p$\snug-vertex' of~$p'$ an instance where
the path~$p'$, in the process of following $\sigma$ or~$\tau$ to
get from one $p$\snug-vertex to the next, passes through
some intermediate vertex which happens to lie on~$p$.)
Since the $\sigma,\tau$\snug-word given by~$p$ is a
nontrivial element of~$F$, the path~$p'$ must go around the prime ring
of~$C_\Fgraph$ at least once.

Whenever $p'$~moves away from the
prime ring, say at a vertex~$u$, it must eventually return to
the prime ring at the same vertex~$u$.  If this part of~$p'$ does
not contain any~$p$\snug-vertices, then it has length less than~$k$.
If the part does contain at least one $p$\snug-vertex, let $x$ and~$y$
be the first and last $p$\snug-vertices encountered along the $p'$\snug-path
from~$u$ to~$u$.  Then the most direct path
in~$C_\Fgraph$ from~$x$ to~$y$ is given
by a reduced $F$\snug-word~$w'$ of length at most~$2k$; while the
path along~$p'$ from~$x$ to~$y$ is given by a word~$w$, not necessarily
reduced.  Since neither of these paths uses the edges of the prime ring,
and since $C_\Fgraph$ has no other cycles, the words $w$ and~$w'$ must
be equal in~$F$.  But $w$ is given by a part of~$p$, so it is the
expansion of a word in $\sigma$ and~$\tau$.  Since this latter word is equivalent to~$w'$, it must have length at most~$N$, by the definition
of~$N$.  So the word~$w$ has length at most~$Nk$, and the entire
part of~$p'$ from~$u$ to~$u$ has length at most~$(N+2)k$.
This was true for every part of~$p'$ off the prime ring, so
every vertex of~$p'$ (in particular, the vertex~$v$ we started with)
is within distance $(N+2)k$ (actually, half that) of the prime ring.
Since $v$~was arbitrary, the value $(N+2)k$ works for~$L_1$.

It remains to find a limit~$L_2$ for the length of the prime ring.
For this, the following fact will be useful: if $w$~is a word in the
generators of~$F$ such that some non-identity power of~$w$ fixes a
vertex~$v$ of~$C_\Fgraph$, then the length of the prime ring is at
most the length of~$w$.  (To see this, note that $w$~is a rotation
of~$S^2$, and the power of~$w$ fixing the vertex~$v$ is a rotation around
the same axis having infinite order; by Lemma~2.1, $w$~itself fixes~$v$.
So $w$~induces a path in~$C_\Fgraph$ from~$v$ to~$v$; since $w$~is not
the identity, this path must traverse the prime ring at least once,
so its length is at least the length of the prime ring.)

Fix a vertex~$v_0$ of~$C_\Fgraph$.  Then we can define a function $h
\colon F \to C_\Fgraph$ by $h(x) = x(v_0)$.  This~$h$ maps the action
of~$F$ on~$F$ by left multiplication to the action of~$F$ on~$C_\Fgraph$;
that is, $g(h(x)) = h(gx)$ for all $g,x \in F$.  So, if $T_F$ is the
Cayley graph of~$F$ (the vertices are the elements of~$F$, and there
is an edge from~$x$ to~$gx$ whenever $g$~is one of the given generators
of~$F$; so $T_F$ is a free~$F$-tree), then $h$~gives a graph homomorphism
from~$T_F$ to~$C_\Fgraph$.

For each vertex~$x$ of~$T_F$ such that $h(x)$~is in~$C_\setgraph$,
give~$x$ the same label that $h(x)$~has (1,~2, or~3).  Then, just
as in~$C_\Fgraph$, if
$x$~has label 1 or~3, then $\sigma(x)$~will have label 1 or~2;
if $x$ has label~2, then $\tau(x)$~will have label~3; and,
if $x$~has label~1, then $\rho(x)$~will have label~3.
This means that each labeled vertex of~$T_F$ is in a well-defined
$\sigma,\tau$\snug-path, which $h$~maps to a $\sigma,\tau$\snug-cycle
in~$C_\setgraph$; however, the $\sigma,\tau$\snug-paths in~$T_F$ are
infinite in both directions.

Each vertex of~$C_\Fgraph$ has infinitely many preimages in~$T_F$.
However, we will now show that each $\sigma,\tau$\snug-cycle
in~$C_\setgraph$ give rise to only finitely many $\sigma,\tau$\snug-paths
in~$T_F$.  Let $x$ be a vertex of~$F$ such that $h(x)$~is in the
$\sigma,\tau$\snug-cycle in question.  There is a minimal word~$w$
in~$F$ which fixes~$h(x)$ (describing a path from~$h(x)$ to itself
which goes around the prime ring once); then the elements of~$F$
which fix~$h(x)$ are just the powers of~$w$, so the $h$\snug-preimages
of~$x$ are the vertices $w^j(x)$ for $j \in \Z$.

The $\sigma,\tau$\snug-cycle containing~$h(x)$ yields a nontrivial word~$u$
in $\sigma$ and~$\tau$ such that $u(h(x)) = h(x)$; hence, $u$~must
be~$w^n$ for some nonzero integer~$n$.  So the $\sigma,\tau$\snug-path
containing~$x$ also contains $u(x) = w^n(x)$; in fact, it
contains $u^j(x) = w^{nj}(x)$ for all $j \in Z$.  This shows
that the preimages of~$h(x)$ lie on at most $|n|$ $\sigma,\tau$\snug-paths,
so, as stated, the $\sigma,\tau$\snug-cycle in~$C_\setgraph$
yields only finitely many $\sigma,\tau$\snug-paths in~$T_F$.

Since $C_\setgraph$ only includes finitely many
$\sigma,\tau$\snug-cycles, there are only finitely many (labeled)
$\sigma,\tau$\snug-paths in~$T_F$.  Note that $\rho^{-1}$ maps the
3-vertices in these paths to 1-vertices in these paths.

We now break into subcases based on the form of~$\rho$.
         
{\bf Subcase 2a}: $\rho$ is in the subgroup $\langle\sigma,\tau\rangle$.

Here we will obtain a value for~$L_2$ depending on the exact form
of~$\rho$ as a word in $\sigma$ and~$\tau$.

Fix a 3-vertex~$x$ in~$T_F$.  If we start at~$x$ and apply all possible
words in $\sigma$ and~$\tau$ (and their inverses),
we get a subset~$T'$ of~$T_F$
(closed under~$\rho$ as well as $\sigma$ and~$\tau$) which can
be viewed as a graph by putting edges from~$y$ to $\sigma(y)$ and~$\tau(y)$
for all vertices~$y$.  Since $\langle\sigma,\tau\rangle$ is free and acts
freely on~$T_F$, $T'$ is isomorphic to the Cayley graph of
$\langle\sigma,\tau\rangle$, which is a free $\sigma,\tau$\snug-tree.
If a labeled vertex is in~$T'$, then its entire $\sigma,\tau$\snug-path
is included in~$T'$, and is a $\sigma,\tau$\snug-path
in~$T'$ (where its vertices will be consecutive, unlike in~$T_F$).

The vertex~$h(x)$ is in a $\sigma,\tau$\snug-cycle; by tracing around
this cycle, we get a nontrivial word~$u$ in $\sigma$ and~$\tau$ (not
using inverses) such that $u(h(x)) = h(x)$.  Let $P$ be
the $\sigma,\tau$\snug-path containing~$x$ in~$T_F$ (and in~$T'$).
So $u^j(x) \in P$ for all $j \in \Z$, and $P$~is a periodic path
with a period given by~$u$.

We next show that $\rho^{-1}(x) \in P$.  Suppose $\rho^{-1}(x)$ is
not in~$P$.  Then, by periodicity, we have $\rho^{-1}(u^j(x))\notin P$
for all integers~$j$.  (If $\rho^{-1}(u^j(x))$ is in $P$, then
it has the form $\alpha u^i(x)$ for some integer~$i$ and some
final segment~$\alpha$ of~$u$.  Since the group action is
free on~$T_F$, this gives $\rho^{-1}u^j = \alpha u^i$, so
$\rho^{-1} = \alpha u^{i-j}$, so $\rho^{-1}(x) \in P$.)
So these vertices $\rho^{-1}(u^j(x))$ must lie in other
$\sigma,\tau$\snug-paths within~$T'$.  But $T'\setminus P$ consists
of infinitely many separate components (each attached to
a single vertex of~$P$), and a $\sigma,\tau$\snug-path
other than~$P$ must lie within one of these components.  Since
the vertices $\rho^{-1}(u^j(x))$ are in separate components,
they must lie in separate $\sigma,\tau$\snug-paths.
This is impossible because there are only finitely many
$\sigma,\tau$\snug-paths in~$T_F $ arising from~$C_\setgraph$.

So $\rho^{-1}(x) \in P$.  The same argument shows that every
3-vertex in any $\sigma,\tau$\snug-path in~$T'$ is sent by~$\rho^{-1}$
to a 1-vertex in that same $\sigma,\tau$\snug-path.

Since $\rho^{-1}(x)$~is in~$P$, $\rho^{-1}$~must be either a final segment
of~$u^n$ for some $n>0$ (so $\rho^{-1}$ is a product of $\sigma$\snug's
and $\tau$\snug's, ending in~$\sigma$) or a final segment of~$(u^{-1})^n$
for some $n>0$ (so $\rho^{-1}$ is a product of $\sigma^{-1}$\snug's and
$\tau^{-1}$\snug's, ending in~$\sigma^{-1}\tau^{-1}$; it can't
be just~$\tau^{-1}$ because $\tau^{-1}(x)$ is a 2-vertex).

If $\rho^{-1}=\sigma$, then we have $\sigma(A_3) = A_1$ and hence
$\sigma(A_1)=A_2$.  This means that $u$~must be a power
of~$\tau\sigma^2$, and $u$ fixes $h(x) \in C_\Fgraph$,
so the length of the prime ring is at most the length of~$\tau\sigma^2$.
So we can let $L_2 = 3k$.

If $\rho^{-1} = \sigma^j$ for some $j>1$, then since $\rho^{-1}(x)$~is on the path~$P$, $\sigma(x)$ must be a 1-vertex rather than a 2-vertex (so we can
apply~$\sigma$ again).  But $\rho(\sigma(x)) = \sigma^{1-j}(x)$ is
not on the path~$P$, so it must be a 3-vertex on some other
$\sigma,\tau$\snug-path.  Hence, $\rho^{-1}$~maps a 3-vertex on that
other path to a 1-vertex on this path, which is impossible by
previous remarks.

If $\rho^{-1}$ is a product involving~$\tau$, let~$j$ be the unique
positive integer such that $\rho$~ends in $\tau\sigma^j$.
Since $\rho^{-1}(x)$ is on~$P$, the path from~$x$ to~$\rho^{-1}(x)$
within~$T'$ must be part of~$P$.  In particular,
$\tau\sigma^j(x)$ is on~$P$, and is the next 3-vertex on~$P$
after~$x$.  But now we can apply~$\rho^{-1}$ to~$\tau\sigma^j(x)$
to get another vertex on~$P$, and use this to conclude that
$(\tau\sigma^j)^2(x)$ is the next 3-vertex on~$P$ after~$\tau\sigma^j(x)$,
and so on.  Eventually the 3-vertex~$u(x)$ must be reached, so
$u$~must be a power of~$\tau\sigma^j$.  Therefore, the length of the
prime ring is at most the length of~$\tau\sigma^j$, so we can
let $L_2 = k(j+1)$.  (Actually, one can show that $j$~must be~2 here,
because $j>2$ would yield a contradiction as in the preceding paragraph
while $j=1$ would give no 1-vertices on~$P$ at all.)

Finally, if $\rho^{-1}$~is a product of $\sigma^{-1}$\snug's and
$\tau^{-1}$\snug's, we can proceed as in the positive cases above.  If
$\rho^{-1}=\sigma^{-1}\tau^{-1}$, then we get
$\sigma(A_1)=A_2$ and hence $\sigma(A_3)=A_1$, so $L_2=3k$ works.
If $\rho^{-1}=(\sigma^{-1})^j\tau^{-1}$ for some $j>1$, we
get a contradiction because $\rho(\sigma^{-1}\tau^{-1}(x))$
is a 3-vertex on some path other than~$P$.  If
$\rho^{-1}$~ends in $\tau^{-1}(\sigma^{-1})^j\tau^{-1}$, then
$u^{-1}$~must be a power of $(\sigma^{-1})^j\tau^{-1}$ and
we can let $L_2 = k(j+1)$.
            
{\bf Subcase 2b}: $\rho$ is not in $\langle\sigma,\tau\rangle$.

Hence, by Lemma~5.3, $\rho\langle\sigma,\tau\rangle\rho^{-1} \cap
\langle\sigma,\tau\rangle$ is a cyclic group.  Let $\beta$~be a generator
of this cyclic group, and let $L_2$ be the length of~$\beta$ as a
word in the generators of~$F$.  To see that this~$L_2$ works, it
will suffice to show that some non-identity power of~$\beta$
fixes a vertex of~$C_\Fgraph$.

As in Subcase~2a, fix a 3-vertex~$x$ in~$T_F$, let $u$ be the
word in $\sigma$ and~$\tau$ such that $u(h(x))=h(x)$ obtained
from the $\sigma,\tau$\snug-cycle containing~$h(x)$,
and let $P$ be the $\sigma,\tau$\snug-path in~$T_F$ containing~$x$.
Since there are only finitely many $\sigma,\tau$\snug-paths
obtained from~$C_\setgraph$, there must be positive integers $i<j$
such that $\rho^{-1}(u^i(x))$ and $\rho^{-1}(u^j(x))$ lie in the
same $\sigma,\tau$\snug-path.  This means that there is
$\gamma \in \langle\sigma,\tau\rangle$ such that
$\gamma(\rho^{-1}(u^i(x))) = \rho^{-1}(u^j(x))$, which implies
$\gamma\rho^{-1}u^i = \rho^{-1}u^j$ and hence
$u^{j-i} = \rho\gamma\rho^{-1}$.
So $u^{j-i}$ is in $\rho\langle\sigma,\tau\rangle\rho^{-1} \cap
\langle\sigma,\tau\rangle$, so it is a non-identity power of~$\beta$
which fixes the vertex~$h(x)$, as desired.  This completes
the proof of Theorem~5.1.
\QED\enddemo

Note that the system of congruences in Theorem~5.1
is actually satisfied by nonempty
finite subsets of a free group: let the group be~$\Z$, and let
$A_1 = \{2\}$, $A_2 = \{1\}$, and $A_3 = \{3\}$.  So we have another
example showing that the `connected and prime' restriction in Theorem~4.3
is needed.

Now let us consider the smaller system of congruences
$$A_1 \cong A_3, \qquad A_1 \cup A_2 \cong A_1 \cup A_3.$$
Since the system from Theorem~5.1 is satisfied by nonempty subsets of a
free group, this subsystem is also satisfied by such sets.  It turns out,
though, that there is basically only one way to get these sets:

\proclaim{Proposition 5.4}
The only finite subsets of a free group satisfying the congruences
$A_1 \cong A_3$, $A_1 \cup A_2 \cong A_1 \cup A_3$ are those for which
$\sigma(A_3) = A_1$ and $\sigma(A_1)=A_2$ for some group element~$\sigma$.
\endproclaim

\demo{Proof} Suppose we have pairwise disjoint finite subsets
$A_1,A_2,A_3$ of a free group, and $\rho$ and~$\sigma$ are elements
of the group such that $\rho(A_1)=A_3$ and $\sigma(A_1 \cup A_3) = A_1
\cup A_2$; we will show that $\sigma(A_3) = A_1$ and $\sigma(A_1)=A_2$.
We may assume that $A_1,A_2,A_3$ are not all empty (otherwise the
conclusion is trivial).

If $\sigma$ and~$\rho$ do not commute, then, as recalled before,
$\sigma$ and~$\rho$ are free generators for the subgroup they generate.
Now, if we start with an element~$x$ of $A_1 \cup A_2 \cup A_3$,
we can apply either $\sigma^{-1}$ (if $x \in A_1 \cup A_2$) or~$\rho^{-1}$
(if $x \in A_3$) to get another element of~$A_1\cup A_2 \cup A_3$.
By iterating this, we can get arbitrarily long words~$w$
in $\sigma^{-1}$ and~$\rho^{-1}$ such that $w(x) \in A_1 \cup A_2 \cup A_3$.
But these elements $w(x)$ are all distinct, because the words~$w$ are
distinct elements of the group, and the group acts freely on itself.
So $A_1 \cup A_2 \cup A_3$ is infinite, contradiction.

Therefore, $\sigma$ and~$\rho$ commute.  Now, since $\sigma(A_3) \subseteq
A_1 \cup A_2$, we have $\sigma(A_3) \cap A_3 = \nullset$.  Applying
$\rho^{-1}$ to this gives $\rho^{-1}(\sigma(A_3) \cap A_3) = \nullset$,
so $\rho^{-1}(\sigma(A_3)) \cap \rho^{-1}(A_3) = \nullset$ since the
mapping given by $\rho^{-1}$ is one-to-one, so $\sigma(\rho^{-1}(A_3))
\cap \rho^{-1}(A_3) = \nullset$ since $\sigma$ and~$\rho$ commute, so
$\sigma(A_1) \cap A_1 = \nullset$.  But $A_1 \subseteq \sigma(A_1 \cup
A_3)$, so $A_1 \subseteq \sigma(A_3)$.  Since $A_1$ and $A_3$ have the
same finite size (because of~$\rho$), we must have $A_1 = \sigma(A_3)$.
This and $\sigma(A_1 \cup A_3) = A_1 \cup A_2$ give $\sigma(A_1) = A_2$,
as desired.
\QED\enddemo

An argument similar to the last part of the proof of Proposition~5.4 shows
that, if we have pairwise disjoint open subsets $A_1,A_2,A_3$ of~$S^2$
and commuting rotations $\rho$ and~$\sigma$ such that $\rho(A_1) = A_3$
and $\sigma(A_1 \cup A_3) = A_1 \cup A_2$, then $\sigma(A_2) = A_1$ and
$\sigma(A_1) = A_3$.  To see this, note that we again have $\sigma(A_3)
\cap A_3 = \nullset$, and applying~$\rho^{-1}$ to both sides gives
$\sigma(A_1) \cap A_1 = \nullset$.  Now, if $\sigma(A_3) \cap A_2$ were
nonempty, then it would be an open subset of~$A_2$ (hence of positive
measure) disjoint from~$\sigma(A_1)$, so $\sigma(A_1)$ would have to
be a subset of~$A_2$ of measure smaller than that of~$A_2$.  But $\rho$
and~$\sigma$ preserve measure, so $A_1$, $A_2$, $A_3$, and $\sigma(A_2)$
all have the same measure.  Therefore, we must have $\sigma(A_3) \cap
A_2 = \nullset$, and this implies $\sigma(A_3) = A_1$ and $\sigma(A_1)
= A_2$, as desired.  (This argument is easier than the proof of Case~1
in Theorem~5.1; unfortunately, the argument for Case~2 is based on the
assumption that $\sigma$ and~$\tau$, rather than $\sigma$ and~$\rho$,
do not commute.)

So, just as in Subcase~1b of the proof of Theorem~5.1, the subsets
$A_1,A_2,A_3$ satisfying these congruences cannot have dense union if
the witnessing rotations $\rho$ and~$\sigma$ are commuting
members of a free group.

Now consider the case where the witnessing rotations $\rho$ and~$\sigma$
do not commute (and hence are free generators for a free group of
rotations of~$S^2$).  Starting with any component of one of the sets
$A_1,A_2,A_3$, we can repeatedly apply $\rho^{-1}$ or~$\sigma^{-1}$
as appropriate to get additional components of these sets of the same
measure as the original; eventually the same component must be repeated.
This means that some nontrivial word in $\rho$ and~$\sigma$ fixes the
original component.  So, by Lemma~2.1, every component of these three
sets must have an axis of symmetry (i.e., must be a disk or an annulus).
This means that a construction like that in Theorem~4.3 will not work
here without substantial modification, because the new components added
at each stage of that construction could be of arbitrary shape, as long
as they were small enough.

The problem in Theorem~5.1 was that there were not enough
ways to satisfy the system of congruences using {\it finite}
subsets of the sphere.  However, for the smaller system we are considering
now, there is a very wide variety of finite sets satisfying it.
Here are a few examples:
$$\gather
A_1 = \{x\},\quad A_2=\{\sigma\rho x\},\quad A_3 = \{\rho x\},
\qquad \text{where $\sigma x = x$;} \\
A_1 = \{y\},\quad A_2=\{\sigma y\},\quad A_3 = \{\rho y\},
\qquad \text{where $\sigma\rho y = y$;} \\
A_1 = \{z,\sigma z,\sigma\rho\sigma z\},\quad A_2=\{\sigma\rho z,
\sigma^2 z, \sigma^2\rho\sigma z\},\quad A_3 = \{\rho z,
\rho\sigma z, \rho\sigma\rho\sigma z\},
\qquad \text{where $\sigma\rho\sigma\rho\sigma z = z$;} \\
{\align
A_1 = \{x, \sigma\rho x, \sigma^2\rho x, \sigma\rho\sigma\rho
x\},\quad A_2=\{&\sigma^3\rho x,\sigma\rho\sigma^2\rho
x, \sigma^2\rho\sigma\rho x, \sigma\rho\sigma\rho\sigma\rho x\}, \\
&\quad A_3 = \{\rho x, \rho\sigma\rho x, \rho\sigma^2\rho
x, \rho\sigma\rho\sigma\rho x\},
\qquad \text{where $\sigma x = x$.}
\endalign}
\endgather$$
In each of these cases we have $\rho(A_1) = A_3$ and
$\sigma(A_1 \cup A_3) = A_1 \cup A_2$.

This gives hope that there are enough such finite sets that we can
use disks or annuli around the points in these sets to form open sets
$A_1,A_2,A_3$ with dense union satisfying the congruences.  (If we use
completely disjoint finite sets, such as in the first three examples
above, then we can try to use small disks around the points to form
suitable open sets.  On the other hand, if we use different finite
subsets which have points in common, such as the first and fourth
examples, then we will have to use annuli rather than disks so that the
point~$x$ itself is not used more than once.)  It turns out that such
a construction is indeed possible, but extreme care is needed to set up
the right inductive hypotheses.

\proclaim{Theorem 5.5} The system of congruences $A_1 \cong A_3$,
$A_1 \cup A_2 \cong A_1 \cup A_3$ can be satisfied by
open subsets of~$S^2$ with dense union using free rotations. \endproclaim

\demo{Proof}
For this proof it will be convenient to specify that distances between
points of~$S^2$ are measured along (minimal) great-circle paths; this means
that the metric (call it~$d$) will be additive along short great-circle arcs.
Let $B(x,\eps)$ denote the open disk $\{y \in S^2 \colon d(x,y) < \eps\}$.

Let $\sigma$ and~$\rho$ be free generators for a free group~$G$
of rotations of~$S^2$ (so $G$~is countable).  We will build pairwise
disjoint open subsets $A_1,A_2,A_3$ of~$S^2$ with dense union such that
$\rho(A_1)=A_3$ and $\sigma(A_1\cup A_3) = A_1\cup A_2$.

Let $p_0$ be a fixed point of~$\sigma$, and let $\orbit$ be the orbit
of~$p_0$ under~$G$.  As we have seen before,
all elements of~$G$ which fix~$p_0$ must
commute with~$\sigma$, and hence (because $G$~is free and $\sigma$~is
one of its generators) must be powers of~$\sigma$. So the action of~$G$
on~$\orbit$ is free except that $\sigma(p_0) = p_0$; more precisely, if
$w$ and~$w'$ are distinct words in~$\sigma,\rho$ neither of which ends
in $\sigma$ or~$\sigma^{-1}$, then $w(p_0) \ne w'(p_0)$.
So, if we view~$\orbit$ as a graph by putting edges from~$x$
to $\sigma(x)$ and~$\rho(x)$ for each $x \in \orbit$, then the graph
will be a free $\sigma,\tau$\snug-tree except for a single loop
from the vertex~$p_0$ to itself.

The general strategy of the proof will be the same as that of
Theorem~4.3 (or Theorem~2.1 of Dougherty~\cite{\Dougherty}, or
Theorem~3.1 of Dougherty and Foreman~\cite{\DoughertyForeman}).
We will construct open subsets $A_1^n,A_2^n,A_3^n,B_1^n,B_2^n,B_3^n$
of~$S^2$, increasing with~$n$, which satisfy a list of
inductive hypotheses.  Fix a list $\langle Z_n\colon n=0,1,2,\dotsc \rangle$
of the nonempty sets in some countable base for~$S^2$.  At stage~$n$,
we will enlarge the sets $A_i^n$ and~$B_i^n$ to sets $A_i^{n+1}$
and~$B_i^{n+1}$ so that the inductive hypotheses are true for the new
sets and at least one of the sets~$A_i^{n+1}$ meets~$Z_n$.
Hence, the sets $A_i = \bigcup_{n=0}^\infty A_i^n$ for $i=1,2,3$
will have dense union, and the inductive hypotheses will ensure
that they are open sets satisfying the system of congruences.
The sets~$B_i^n$ contain the points that are to be explicitly
excluded from~$A_i^N$ for all $N$.

For convenience, here is a list of all of the inductive hypotheses
to be used (some of which mention terms to be defined later).

\roster
\item $A_i^n$ and $B_i^n$ are open sets whose boundaries do not contain
any of the points in~$\orbit$.

\item $A_i^n \subseteq B_j^n$ for $j \ne i$.

\item $B_1^n\cap B_2^n \cap B_3^n = \nullset$.

\item $\rho(A_1^n) = A_3^n$ and $\sigma(A_1^n \cup A_3^n) = A_1^n \cup A_2^n$.

\item $\rho(B_1^n) = B_3^n$, $\sigma(B_1^n \cap B_3^n) = B_1^n \cap B_2^n$,
$\rho(B_2^n \cap B_3^n) = B_1^n \cap B_2^n$, and $\sigma(B_2^n) = B_3^n$.

\item
There exist a positive number~$\delta_0$ and an acceptable path which
starts at~$p_0$ and ends at a point in $B(p_0,\delta_0)\setminus\{p_0\}$
such that, for every point~$x$ on the path, $B(x,\delta_0)$ is disjoint
from $B_1^n \cup B_2^n \cup B_3^n$.

\item
For any point $p \ne p_0$ in~$\orbit$ but not in
$A_1^n \cup A_2^n \cup A_3^n$, and any~$i$ such that $p \notin B_i^n$,
there is an acceptable path from some point in $B(p_0,\delta_0)$ to~$p$
which gives~$p$ the label~$i$.

\item
There is a natural number~$M$ such that, for any allowed modification
to a point in~$\orbit$, the relevant modification propagation algorithm
will terminate and will only modify points at most $M$~steps from the
original point.
\endroster

Requirement~\ri1 ensures that the unions~$A_i$ are open sets, requirements
\ri2 and~\ri3 ensure that the sets~$A_i$ are pairwise disjoint (because,
for any~$n$, a point in two of the sets~$A_i^n$ must be in all three of
the sets~$B_j^n$), and requirement~\ri4 ensures that $\rho(A_1)=A_3$ and
$\sigma(A_1 \cup A_3) = A_1 \cup A_2$.  So carrying out the construction
of the sets $A_i^n$ and~$B_i^n$ as above will suffice to prove the
theorem.

As in previous proofs, the process of enlarging the sets $A_i^n$ and~$B_i^n$
to sets $A_i^{n+1}$ and~$B_i^{n+1}$ consists of two phases.  In the
first phase, we build {\it finite} sets which can be added to the current
sets to yield sets $\hat A_i$ and~$\hat B_i$ satisfying most of the
inductive hypotheses (specifically, \ri2--\ri5 and~\ri8).
Then, in the second phase, each new point is replaced with a small
open annulus to yield sets $A_i^{n+1}$ and~$B_i^{n+1}$ satisfying all of
the inductive hypotheses.

The first of these two phases involves difficulties that did not occur
in the previous proofs.  If we want to add a point~$y$ to $\hat A_1$
or~$\hat A_2$, we must add~$\sigma^{-1}(y)$ to $\hat A_1$ or~$\hat A_3$;
if we add~$y$ to~$\hat A_3$, we must add~$\rho^{-1}(y)$ to $\hat A_1$.
This leads to an apparently infinite sequence of points to be added
to the three sets~$\hat A_i$.  In order to have this process add only
finitely many points to these sets, we will have to arrange for the
sequence of new points to loop back on itself.  We will do this by
causing the sequence of points to terminate at the point~$p_0$.
(If we add~$p_0$ to~$\hat A_1$, we do not
have to go back any farther, because $\sigma^{-1}(p_0) = p_0$.)
Of course, we must ensure that
there is a sufficient supply of target points~$y$ from which one
can backtrack to~$p_0$ in this way; we will need to find such a~$y$
in the given open set~$Z_n$.

Let us define an {\it acceptable path}
from~$x$ to~$y$ (for the sets
$A_j^n$ and~$B_j^n$) to be a sequence $x_0,x_1,\dots,x_k$ with
$x_0=x$ and $x_k=y$, together with
a labeling which assigns either 1,~2, or~3 to each~$x_i$ so
that:
\roster
\item"$\bullet$" $x_0$~is labeled~1;
\item"$\bullet$" if $x_i$ is labeled~$j$, then $x_i \notin B_j^n$; and
\item"$\bullet$" for any consecutive points $x_i$ and~$x_{i+1}$ on the path,
either $x_{i+1}=\rho(x_i)$, $x_i$~is labeled~1, and $x_{i+1}$~is
labeled~3, or $x_{i+1}=\sigma(x_i)$, $x_i$~is labeled 1 or~3, and $x_{i+1}$~is
labeled 1 or~2.
\endroster
Let $\alpha_{i+1}$ be $\rho$ if $x_{i+1}$ is labeled~3, $\sigma$~otherwise,
so that we have $x_{i+1}=\alpha_{i+1}(x_i)$ for all~$i$.  The word
$w = \alpha_k \alpha_{k-1} \dots \alpha_1$ in $\sigma$ and~$\tau$ is said to
be {\it associated} with the given acceptable path; note that $w(x) = y$.

It would be natural to expect that the inductive hypothesis on
reachability would have the form ``the set of all~$y$ such that there
is an acceptable path from~$p_0$ to~$y$ assigning to~$y$ the label~$i$
is dense in~$S^2 \setminus B_i^n$.''  (We could also exclude $A_1^n \cup
A_2^n \cup A_3^n$ if necessary.)  However, it seems to be necessary
to make a stronger inductive hypothesis in order to prove the hypotheses
at stage~$n+1$ from the hypotheses at stage~$n$.  The stronger version
is given in hypotheses \ri6 and~\ri7.  It is not yet obvious that these
hypotheses imply the statement above; we will prove that next.

Given~$n$, fix~$\delta_0$ and and an acceptable path~$P_0$ as in~\ri6.
Let $Q_0$ be the set of points on this path, let $w_0$ be the word
associated with the path, and let $D_0 = B(p_0,\delta_0)$.
We may assume that the final point~$w_0(p_0)$ is not labeled~2,
because, if it were, we could just change that label to~1
(since $w_0(p_0)$~is known not to be in~$B_1^n$).
Let $\delta_1 = d(p_0,w_0(p_0))$.  (Later we will define further
$\delta$\snug's, with $\delta_0 > \delta_1 > \delta_2 > \delta_3 >
\dotsb$.)

\procl{Claim 1}
Let $\eps$ be a positive number less than $\delta_0 - \delta_1$.
Let $E$ be the set of all points~$x \in B(p_0,\delta_0-\eps)$
for which there is an acceptable path from~$p_0$ to~$x$ such that
$x$~is labeled\/~$1$ and, for every~$y$ on the path, $B(y,\eps)$ is
disjoint from $B_1^n \cup B_2^n \cup B_3^n$.  Then $E$~is dense in
$B(p_0,\delta_0-\eps)$.

\procl{Proof} Since the path~$P_0$ satisfies the conditions in~\ri6,
all of the disks $B(y,\delta_0)$ for $y \in Q_0$, including~$D_0$, are
disjoint from $B_1^n \cup B_2^n \cup B_3^n$.  Clearly $p_0$~is in~$E$
(via a trivial path).
Now, if $x \in E$, then $\sigma(x)$ is in $B(p_0,\delta_0-\eps)$,
so $B(\sigma(x),\eps) \subseteq D_0$.  So we can append~$\sigma(x)$
(with label~1) to the acceptable path from~$p_0$ to~$x$ to get such
a path from~$p_0$ to~$\sigma(x)$; this shows that $\sigma(x) \in E$
for all $x \in E$.

Next, suppose $x \in E$ and $d(w_0(x),p_0) < \delta_0 - \eps$.
The points of~$Q_0$ are the points of the form~$w(p_0)$ where
$w$~is a final segment of the word~$w_0$.  For each such~$w$,
we have $d(w(x),w(p_0)) < \delta_0 - \eps$, so
$B(w(x),\eps)$ is included in $B(w(p_0),\delta_0)$ and hence is disjoint
from $B_0^n \cup B_1^n \cup B_2^n$.  So we get an acceptable path
from~$p_0$ to~$w_0(x)$ by taking the given path from~$p_0$ to~$x$
and appending the points~$w(x)$ for $w$~a nonempty final segment
of~$w_0$ (where the label assigned to~$w(x)$ is the same as the
label assigned to~$w(p_0)$ in~$P_0$).  This path has all the properties
needed in the definition of~$E$, except that the final label may be
1 or~3; but in either case we can append~$\sigma(w_0(x))$ with
label~1.  So we have shown that, if $x \in E$ and
$d(w_0(x),p_0) < \delta_0 - \eps$, then $\sigma(w_0(x)) \in E$.

In particular, since $p_0 \in E$, we get $\sigma(w_0(p_0)) \in E$, so
$\sigma^j(w_0(p_0)) \in E$ for all $j>0$.  The points
$\sigma^j(w_0(p_0))$ are dense in the circle around~$p_0$ of radius
$\delta_1$.  Applying~$w_0$ to this circle gives another
circle of radius~$\delta_1$, which passes through the point~$p_0$, so there
are points of~$w_0(E)$ on this latter circle at a set of distances from~$p_0$
which is dense in the interval $[0,2\delta_1]$.  If such a point~$y$ is
within distance $\delta_0 - \eps$ of~$p_0$, then $\sigma(y)$~is
in~$E$; in fact, we have $\sigma^j(y) \in E$ for all $j>0$.
This shows that $E$~is dense in $B(p_0, \min(2\delta_1,\delta_0-\eps))$.

By applying $w_0$ and positive powers of~$\sigma$ again, we can show
that $E$ is dense in $B(p_0, \min(3\delta_1,\allowbreak
\delta_0-\nobreak\eps))$.
We can repeat this until we reach a multiple~$k\delta_1$ so large that
$k\delta_1 > \delta_0 - \eps$; this will show that $E$~is dense in
$B(p_0,\delta_0-\eps)$, as desired. \QEd\par\smallskip

So, given an open set~$Z$ which we want to reach by an acceptable path
from~$p_0$, we can start by finding a point $y \in \orbit \cap Z$,
fixing an~$i$ such that $y \notin B_i^n$ by~\ri3, and applying~\ri7
to get a point~$x \in D_0$ and an acceptable path from~$x$ to~$y$, with
associated word~$w$.  Let $\eps$ be a positive number less than $\delta_0
- d(p_0,x)$.  Then the set~$E$ in Claim~1 will contain points~$x'$
arbitrarily close to~$x$.  In particular, we will be able to make~$x'$
so close to~$x$ that the path from~$x'$ to $y'=w(x')$ given by the
word~$w$ is still acceptable (with the same labeling as the path from~$x$
to~$y$) and $y'$~is also in~$Z$.  Concatenating this path with the path
from~$p_0$ to~$x'$ as in the definition of~$E$ will yield an acceptable
path from~$p_0$ to a point $y' \in Z$.

The labeling of an acceptable path from~$p_0$ indicates which of
the points on the path are to be added to the given sets sets~$A_i^n$
(and to~$B_j^n$ for $j \ne i$) in the process of building the
new sets~$\hat A_i$ (and~$\hat B_j$).  However, just adding the
points on the path will not be enough to satisfy requirements
\ri4 and~\ri5; other points not on the path will have to be added
as well.  We now describe two modification propagation algorithms;
one shows how to add additional points to the sets~$\hat A_i$, and
the other shows how to add additional points to the sets~$\hat B_i$.

An {\it allowed modification} for a point~$x$ means either adding~$x$
to a set~$\hat A_i$, where $x$~is not in~$B_i^n$, or adding~$x$ to one
or more of the sets~$\hat B_j$, making sure that $x$~does not end up in
all three of the sets~$\hat B_j$.

Once an allowed modification has been made at~$x$, we may have to make
modifications at other points $\sigma^{\pm 1}(x)$ or~$\rho^{\pm 1}(x)$
in order to make \ri4 and~\ri5 true for the sets $\hat A_j$
and~$\hat B_j$.  These new modifications may entail further
modifications, and so on.  We will now give a more precise description
of the algorithms for propagating these modifications, one for the
sets~$\hat A_j$ and one for the sets~$\hat B_j$.

For the purposes of these algorithms, a point in~$\orbit$ is said to be
{\it eligible} if it is not~$p_0$ and has not yet been modified for the
sets in question ($\hat A_j$ or~$\hat B_j$).  The reason for excluding~$p_0$
is to avoid the loop in the graph~$\orbit$; this means that
all propagation will move away from the initial modification in this graph,
and no point will be reached more than once.  The exclusion does not
matter, because $p_0$~will always be part of the initial modification
before the propagation takes place; it just makes requirement~\ri8
easier to state.

{\narrower\narrower
\noindent Modification propagation algorithm A: Suppose that
the point $x \in \orbit$ has been added to~$\hat A_i$.  Then do all
of the following:
\roster
\item"$\bullet$" If $i=1$ and $\rho(x)$~is eligible, then
add~$\rho(x)$ to~$\hat A_3$ and apply algorithm~A to $\rho(x)$.
\item"$\bullet$" If $i=1$ or $i=3$, and $\sigma(x)$~is eligible,
then add~$\sigma(x)$ to either $\hat A_1$ or~$\hat A_2$ (specifically,
add~$\sigma(x)$ to~$\hat A_2$ if $\sigma(x) \notin B_2^n$, and to~$\hat
A_1$ otherwise), and apply algorithm~A to $\sigma(x)$.
\endroster

}

It is important to note that the additional modifications produced
in algorithm~A are allowed modifications.  If $x$~has been added
to~$\hat A_1$ by an allowed modification, then $x \notin B_1^n$,
so $\rho(x) \notin B_3^n$ by~\ri5, so adding~$\rho(x)$
to~$\hat A_3$ is allowed.  Similarly, if $x$ has been added to
$\hat A_1$ or~$\hat A_3$, then $x \notin B_1^n \cap B_3^n$,
so $\sigma(x) \notin B_1^n \cap B_2^n$, so adding~$\sigma(x)$
to $\hat A_1$ or~$\hat A_2$ is allowed.

Algorithm~A only proceeds forward, from~$x$ to $\rho(x)$ and~$\sigma(x)$;
it does not try to propagate modifications backward to $\rho^{-1}(x)$
and~$\sigma^{-1}(x)$.  As described earlier, such backward propagation
would normally be endless; the acceptable paths are specifically designed
to handle this.

Another fact we will need later is that the modifications produced
by algorithm~A do not reach points that were already in one of the
sets~$A_i^n$ (assuming the initially modified point was not already
in one of these sets).  Of course, the algorithm would not add a point
to~$\hat A_i$ if it were already in~$A_j^n$ for some $j \ne i$, because
the point would be in~$B_i^n$ and the modification would not be allowed.
But the algorithm also will not do a redundant addition (adding a point
to~$\hat A_i$ when it is already in~$A_i^n$).  For instance, if $x$~has
just been added nonredundantly to~$\hat A_1$, then $x$~is not in~$A_1^n$,
so $\rho(x)$~is not in~$A_3^n$, so the resulting addition of~$\rho(x)$
to~$\hat A_3$ is nonredundant.

Once the construction of the sets~$\hat A_i$ is complete, we will need
to build the sets~$\hat B_i$.  We start with $\hat B_i = B_i^n$.
Next come the initial modifications: for each point~$x$ which has been
added to one of the sets~$\hat A_i$, add~$x$ to~$\hat B_j$ for
all $j \ne i$.  This is an allowed modification because $x$~could
not have been added to~$\hat A_i$ if it were in~$B_i^n$, so $x$~will not
be in $\hat B_i$.  These initial modifications now propagate according
to the following algorithm:

{\narrower\narrower
\noindent Modification propagation algorithm B: Suppose that
the point $x \in \orbit$ has been added to one or more of
the sets~$\hat B_j$.  Then do all of the following:
\roster
\item"$\bullet$" If $\rho(x)$~is eligible, then: if
$x$~is now in~$\hat B_1$, add~$\rho(x)$ to~$\hat B_3$;
and if $x$~is now in $\hat B_2 \cap \hat B_3$, add~$\rho(x)$
to~$\hat B_1$ and to~$\hat B_2$.
If this actually modifies~$\rho(x)$ (i.e., $\rho(x)$~has been added
to a set~$\hat B_i$ that it was not already in), then
apply algorithm~B to~$\rho(x)$.
\item"$\bullet$" If $\rho^{-1}(x)$~is eligible, then: if
$x$~is now in~$\hat B_3$, add~$\rho^{-1}(x)$ to~$\hat B_1$;
and if $x$~is now in $\hat B_1 \cap \hat B_2$, add~$\rho^{-1}(x)$
to~$\hat B_2$ and to~$\hat B_3$.
If this actually modifies~$\rho^{-1}(x)$, then
apply algorithm~B to~$\rho^{-1}(x)$.
\item"$\bullet$" If $\sigma(x)$~is eligible, then: if
$x$~is now in~$\hat B_2$, add~$\sigma(x)$ to~$\hat B_3$;
and if $x$~is now in $\hat B_1 \cap \hat B_3$, add~$\sigma(x)$
to~$\hat B_1$ and to~$\hat B_2$.
If this actually modifies~$\sigma(x)$, then
apply algorithm~B to~$\sigma(x)$.
\item"$\bullet$" If $\sigma^{-1}(x)$~is eligible, then: if
$x$~is now in~$\hat B_3$, add~$\sigma^{-1}(x)$ to~$\hat B_2$;
and if $x$~is now in $\hat B_1 \cap \hat B_2$, add~$\sigma^{-1}(x)$
to~$\hat B_1$ and to~$\hat B_3$.
If this actually modifies~$\sigma^{-1}(x)$, then
apply algorithm~B to~$\sigma^{-1}(x)$.
\endroster

}

Note that algorithm~B propagates in all directions, not just forward;
this does not lead to an infinite regress here.  Note also that
the propagated modifications are allowed, assuming the original one was.
For instance, consider~$\rho(x)$.  Since the modification at~$x$
is allowed, either $x$~does not end up in~$\hat B_1$ or it does
not end up in $\hat B_2 \cap \hat B_3$.  If $x$~does not end up
in~$\hat B_1$, then $\rho(x)$~will not be added to~$\hat B_3$,
and $\rho(x)$~cannot have been in~$B_3^n$ to start with
(because that and~\ri5 would give $x \in B_1^n$ and hence
$x \in \hat B_1$), so $\rho(x)$~does not end up in~$\hat B_3$.
Similarly, if $x$~does not end up in
$\hat B_2 \cap \hat B_3$, then $\rho(x)$~is added to neither~$\hat B_1$
nor~$\hat B_2$, and it cannot have been in~$B_1^n \cap B_2^n$
to start with (by~\ri5 again), so it does not end up
in $\hat B_1 \cap \hat B_2$.  So, in any case,
$\rho(x)$~does not end up in all three of the sets~$\hat B_i$.
The same applies to propagation in the other three directions.

Unlike algorithm~A, algorithm~B can add a point to sets it was already in.
But if a point~$x$ is already in two of the sets~$\hat B_i$,
then it cannot be added to the third such set, so
$x$~will not be added to any {\sl new} sets and algorithm~B
will not be applied to~$x$.

We can now describe the full construction of the intermediate sets
$\hat A_i$ and~$\hat B_i$.  Start with $\hat A_i = A_i^n$ and
$\hat B_i = B_i^n$.  Find an acceptable path from~$p_0$ to a point in
the target open set, and add the points on the path to the
sets~$\hat A_i$ as specified by the labeling of the path.  Now apply
algorithm~A to all of the points on the path (with adjacent points on the
path being ineligible because they have already been modified).
For each point that has been added to one of the sets~$\hat A_i$ this
way (the points on the path and the points modified by algorithm~A),
add that point to the sets~$\hat B_j$ for all $j \ne i$.  Now apply
algorithm~B to all of these points to complete the construction
of the sets~$\hat B_i$.

This construction adds only finitely many points to the open sets
we started with ($A_i^n$ and~$B_i^n$).  The initial acceptable path
is finite, and algorithm~A terminates after finitely many steps
for each point on the path, by~\ri8.  (Hypothesis~\ri8 refers to running the
algorithm with a single starting point.  But running it starting
from an entire path of initial modifications simply means that
more points will be declared ineligible for each individual
run of the algorithm; this can only make the algorithm terminate sooner.)
This gives finitely many initial modifications for the sets~$\hat B_i$,
and then \ri8 ensures that all of the required executions of algorithm~B
will terminate after finitely many steps as well.  So the whole process
is finite.

\procl{Claim 2}
Suppose we have an acceptable path (for the sets $A_i^n$ and~$B_i^n$)
starting at~$p_0$.
If we follow the procedure above to construct sets $\hat A_i$ and~$\hat
B_i$, then these sets will satisfy \ri2--\ri5 and~\ri8.

\procl{Proof}
We took care of~\ri2 by adding the new points in~$\hat A_i$ to~$\hat B_j$
for all $j \ne i$.  The fact that all modifications to the sets~$\hat B_i$
were allowed implies that \ri3~holds for the resulting sets.

For~\ri4, we will show that $x \in \hat A_1$ if and only if $\rho(x) \in
\hat A_3$ for any point~$x$; the proof that $x \in \hat A_1 \cup \hat A_3$
iff $\sigma(x) \in \hat A_1 \cup \hat A_2$ is similar.
First, suppose $x \in \hat A_1$.  If $x \in A_1^n$, then
$\rho(x) \in A_3^n$ by the old~\ri4.  If $x$~is on the acceptable path
and $\rho(x)$~is also on the path, then $\rho(x)$~must be the next point
after~$x$ on the path and must be labeled~3, so $\rho(x) \in \hat A_3$.
If $x$~is on the path but $\rho(x)$~is not, or if $x$~was added by
algorithm~A, then $\rho(x)$~must have been added to~$\hat A_3$
when algorithm~A was applied to~$x$ (note that $\rho(x)$~must have been
eligible, because algorithm~A always moves farther away from the
path, never toward it, and the tree structure of~$\orbit$ guarantees
that $\rho(x)$~could not have been reached from the path by
any other route).  So $\rho(x) \in \hat A_3$ in any case.

Now suppose $\rho(x) \in \hat A_3$.  If $\rho(x) \in A_3^n$, then
$x \in A_1^n$ by the old~\ri4.  If $\rho(x)$~is on the acceptable path
and is labeled~3, then its predecessor on the path must be~$x$, and
$x$~must be labeled~1.  If $\rho(x)$~was added to~$\hat A_3$ by
algorithm~A, this must be because $x$~had previously been added
to~$\hat A_1$.  So, in any case, $x \in \hat A_1$.

For~\ri5, we again handle the $\rho$~case; the $\sigma$~case is similar.
Let $x \in S^2$ be arbitrary; we must show that $x \in \hat B_1$ iff
$\rho(x) \in \hat B_3$ and $x \in \hat B_2 \cap \hat B_3$ iff $\rho(x)
\in \hat B_1 \cap \hat B_2$.  If neither $x$ nor~$\rho(x)$ was ever
modified, then this follows from the old~\ri5.  Another case is when
each of $x$
and~$\rho(x)$ was added to a set~$\hat A_i$ (initially or by
algorithm~A); but the only way in which this can happen is when $x$~is
added to~$\hat A_1$ and $\rho(x)$~is added to~$\hat A_3$.  In this case,
$x$~is in $\hat B_2$ and~$\hat B_3$ but not~$\hat B_1$, while
$\rho(x)$~is in $\hat B_1$ and~$\hat B_2$ but not~$\hat B_3$, so the
desired relationships hold.

The remaining possibility is that $x$ or~$\rho(x)$ or both was modified by
algorithm~B.  Because of the tree structure of~$\orbit$ (other than
at~$p_0$, which was initially added to $\hat B_2$ and~$\hat B_3$),
we must have that either $\rho(x)$~was an eligible point when
algorithm~B was applied to~$x$, or vice versa.  If it is the former,
then $\rho(x)$~was added to~$\hat B_3$ iff $x$~was in~$\hat B_1$; and
if $x$~was not in~$\hat B_1$, then $\rho(x)$~could not have been
in~$B_3^n$ to start with (if it were, then $x$~would be in~$B_1^n$
by the old~\ri5), so $\rho(x)$~did not end up in~$\hat B_3$.
Therefore, we do get $x \in \hat B_1$ iff $\rho(x) \in \hat B_3$.
Similar reasoning shows that $x \in \hat B_2 \cap \hat B_3$ iff
$\rho(x) \in \hat B_1 \cap \hat B_2$.  The same argument works
if $x$~was eligible when algorithm~B was applied to~$\rho(x)$.
So \ri5~holds for the new sets~$\hat B_i$.

It remains to show that the new sets satisfy~\ri8.  Let $M$ be the
number from~\ri8 for the sets $A_i^n$ and~$B_i^n$.  Let $Q_1$ be the
set of points modified at any time during the construction (the points on
the acceptable path and the points modified by algorithms A and~B).
Then $Q_1$~is a finite connected subset of~$\orbit$ containing~$p_0$;
let $N$~be the largest number of edges for a non-self-intersecting path
within~$Q_1$ (ignoring orientation of edges as usual).  If we were to
make a new allowed modification and then run algorithm A or~B, then
the algorithm would only reach points within $M$~steps of the starting
point unless it reached a point in~$Q_1$.  In this case, the algorithm
can proceed at most $N$~steps farther before leaving~$Q_1$; after
leaving~$Q_1$, it can proceed at most $M$~steps farther before halting
(it cannot reenter~$Q_1$ because of the tree structure of~$\orbit$).
So, in all, the algorithm cannot go farther than $2M+N+1$ steps from
the starting point.
\QEd\par\smallskip

Once we have suitable intermediate sets $\hat A_i$ and~$\hat B_i$,
we will build new open sets $A_i^{n+1}$ and~$B_i^{n+1}$ by replacing
each point in $\hat A_i \setminus A_i^n$ or $\hat B_i \setminus B_i^n$
with a very small new annulus.  (So the new points in the intermediate
sets will actually not be put in the new open sets.  In particular,
$p_0$~will not be in any of the new open sets; this will leave $p_0$
free to be used again in the construction at the next stage.)  As in
previous proofs, most of the induction hypotheses will hold for the new
open sets because they hold for the intermediate sets; we will have to
argue directly that \ri6 and~\ri7 hold for the new open sets.

The plan of the proof has now been presented; it remains to fill in
the rest of the details.

To start with, let $A_i^0 = B_i^0 = \nullset$ for $i=1,2,3$.  It is
obvious that hypotheses \ri1--\ri5 hold for these sets.  For~\ri6, we
can let $\delta_0$ be so large that $B(p_0,\delta_0)$ is the entire
sphere~$S^2$; the required acceptable path is just the single step
from~$p_0$ to~$\rho(p_0)$.  For~\ri7, the required acceptable path will
have zero, one, or two steps, depending on whether $i$ is 1,~3, or~2.

Finally, for~\ri8, we can use the value $M=4$.  In fact, as exhaustive
checking of the possibilities will verify, modification propagation
algorithm A always terminates within two steps of the starting point
(worst case: adding $x$ to~$\hat A_1$ will cause $\sigma\rho(x)$ to be
added to~$\hat A_2$), while modification propagation algorithm B can
go up to four steps from the starting point (worst case: adding $x$ to
$\hat B_1 \cap \hat B_3$ will cause $\rho^{-1}\sigma\rho^{-1}\sigma(x)$
to be added to~$\hat B_1$).

This completes the initialization of the construction.  Now, suppose
we have already constructed sets $A_i^n$ and~$B_i^n$ satisfying the
inductive hypotheses.  Let $Z_n$ be a nonempty open subset of~$S^2$.
We must show how to enlarge $A_i^n$ and~$B_i^n$ to sets $A_i^{n+1}$
and~$B_i^{n+1}$ so that the inductive hypotheses are true for the new
sets and at least one of the sets~$A_i^{n+1}$ meets~$Z_n$.

If $Z_n$~intersects $A_1^n \cup A_2^n \cup A_3^n$, then we do not have
to do anything at stage~$n$; just let $A_i^{n+1}=A_i^n$ and $B_i^{n+1} =
B_i^n$.  So suppose $Z_n$~is disjoint from $A_1^n \cup A_2^n \cup A_3^n$.

Let us say that a disk $B(x,\eps)$ {\it lies on one side of} a set~$A$
if and only if $B(x,\eps) \subseteq A$ or $B(x,\eps) \subseteq S^2 \setminus
A$.  If $x$~is not a boundary point of~$A$, then, for any sufficiently small
positive number~$\eps$, $B(x,\eps)$ lies on one side of~$A$.
In particular, by~\ri1, this holds when $x \in \orbit$ and
$A$~is one of the sets $A_i^n$ or~$B_i^n$.

As described in the paragraph after Claim~1, let $y$ be a point in~$\orbit
\cap Z_n$ other than~$p_0$.  By~\ri3, there is an~$i$ such that $y
\notin B_i^n$.  Fix $\delta_0,P_0,Q_0,w_0,D_0,\delta_1$
as described just before Claim~1.
Apply~\ri7 to get an acceptable path from a point $x
\in D_0$ to~$y$, with associated word~$w$.  So the points on the path
are the points~$v(x)$ where $v$~is a final segment of~$w$.  By~\ri1,
none of these points is a boundary point of any of the sets $A_j^n$
or~$B_j^n$, so we can find a number $\eps > 0$ so small that each of
the disks $B(v(x),\eps)$ (where $v$~is a final segment of~$w$) lies on
one side of each of the sets $A_j^n$ and~$B_j^n$.  We may
also assume that $\eps < \delta_0 - d(p_0,x)$ and
$\eps < \delta_0 - \delta_1$, and that $B(y,\eps) \subseteq Z_n$.

By Claim~1, the neighborhood $B(x,\eps)$ meets the set~$E$ defined in
that claim; let $x'$ be a point in $E \cap B(x,\eps)$.  Let~$P'$
be an acceptable path from~$p_0$ to~$x'$ such that $x'$~is labeled~1.
Then we can extend~$P'$ to an acceptable path from~$p_0$ to $y' = w(x')$
by following the word~$w$ and using the same labeling as in the
acceptable path from~$x$ to~$y$.  This is a suitable labeling
because, for any final segment~$v$ of~$w$, $v(x')$~is in $B(v(x),\eps)$,
which lies on one side of each set~$A_i^n$ or~$B_i^n$, so
$v(x)$ and~$v(x')$ are in the same such sets, so a label which
is allowed for~$v(x)$ is also allowed for~$v(x')$.  So we have an
acceptable path~$P$ from~$p_0$ to a point $y' \in Z_n$.

Note that none of the points on the path~$P$ are in any of the
sets~$A_i^n$.  (If a point~$z$ on the path were in~$A_i^n$, then it would
also be in~$B_j^n$ for $j \ne i$, so the only possible label for~$z$
would be~$i$.  Now, using \ri4~and the definition of an acceptable path,
we see that the point following~$z$ on the path is also in one of the
sets~$A_i^n$.  Repeating this, we eventually conclude that $y'$~is
in one of the sets~$A_i^n$; this is impossible because $y'$~was
chosen from~$Z_n$.)

Using the acceptable path~$P$, construct the sets $\hat A_i$ and~$\hat B_i$
as described before Claim~2.  By Claim~2, these sets satisfy \ri2--\ri5
and \ri8; let~$\hat M$ be the bound obtained from~\ri8 for these sets
(while~$M$ is the bound for the sets~$A_i^n$ and~$B_i^n$).
Also, let~$Q_1$ be the set of points modified during the construction
of $\hat A_i$ and~$\hat B_i$.

Since the sets $Q_0$ and~$Q_1$ are finite, we can
find a positive number $\delta_2 < \delta_1$ so small that there do not exist
$x \in Q_0$ and $y \in Q_1$ such that either $0 < d(x,y) \le \delta_2$
or $0 < |\delta_1 - d(x,y)| \le \delta_2$.

Since $d(w_0^{-1}(p_0),p_0) = \delta_1$, the point~$w_0^{-1}(p_0)$
lies on the circle with center~$p_0$ and radius~$\delta_1$.
The points $\sigma^n(w_0(p_0))$ for positive integers~$n$ are dense in
this circle, since $\sigma$~is a rotation of infinite order around~$p_0$.
Therefore, we can fix a number $n_1 > 0$ such that
$d(\sigma^{n_1}w_0(p_0),w_0^{-1}(p_0)) < \delta_2$.

Let $\delta_3 = d(\sigma^{n_1}w_0(p_0),w_0^{-1}(p_0))$ and
$D_3 = B(p_0,\delta_3)$.

\procl{Claim 3} For any point $y \ne p_0$ in~$\orbit$ but not in $A_1^n
\cup A_2^n \cup A_3^n$, and any~$i$ such that $y \notin B_i^n$, there
is an acceptable (for the sets $A_j^n$ and~$B_j^n$) path from some $x
\in D_3$ to~$y$ which gives~$y$ the label~$i$.  Similarly, for any point
$y \in \orbit$ which is not in $\hat A_1 \cup \hat A_2 \cup \hat A_3$,
and any~$i$ such that $y \notin \hat B_i$, there is an acceptable (for
the sets $\hat A_j$ and~$\hat B_j$) path from some $x \in D_3$ to~$y$
which gives~$y$ the label~$i$.

\procl{Proof}
For the first part, begin by applying~\ri7 to get an acceptable path~$P_1$
from some $x \in D_0$ to~$y$ which gives~$y$ the label~$i$.  Choose a
positive number~$\eps$ less than $\delta_3$ and also less than $\delta_0 -
\delta_1$.  Since $x \in B(p_0,\delta_0)$, the open set $B(x,\eps)
\cap B(p_0,\delta_0-\eps)$ is nonempty.  Apply Claim~1 to get a point~$x'$
in this open set which is in the set~$E$ (defined in Claim~1).  Let $P_2$
be the acceptable path from~$p_0$ to~$x'$ as specified in the definition
of~$E$, and let $w$ be the associated word for $P_2$.  Then we can get a
path~$P_3$ from~$w^{-1}(x)$ to~$x$ using the same associated word~$w$
and the same labeling as for~$P_2$.  Each point on~$P_3$ is within
distance~$\eps$ of the corresponding point on~$P_2$, and hence is not
in $B_1^n \cup B_2^n \cup B_3^n$.  Hence, $P_3$~is an acceptable path,
and appending~$P_1$ to~$P_3$ gives the desired acceptable path from a
point in~$D_3$ (since $d(x',x) < \eps$, we have
$d(p_0,w^{-1}(x)) < \eps$) to~$y$.

For the second part, we start constructing the acceptable path backward
from~$y$.  If $y \notin \hat B_3$, then $y_1 = \rho^{-1}(y)$ is not
in~$\hat B_1$ by~\ri5, and $y_1$~is also not in~$\hat A_1$ by~\ri4
(since $y \notin \hat A_3$), so, by~\ri2, $y_1$~is not in any of the
sets~$\hat A_i$.  Similarly, if $y\notin \hat B_1$ or $y \notin \hat
B_2$, and $y_1 = \sigma^{-1}(y)$, then we get $y_1 \notin \hat B_1$ or
$y_1 \notin \hat B_3$, as well as $y_1 \notin \hat A_1 \cup \hat A_2
\cup \hat A_3$.  In either case, $y_1$~satisfies the hypotheses of the
claim.  We can now apply the same reasoning to~$y_1$ to get a new
point~$y_2$ (either $\rho^{-1}(y_1)$ or~$\sigma^{-1}(y_1)$), and so on
as many times as desired.  (Note that we never reach the point~$p_0$,
because $p_0 \in \hat A_1$.)

Each member of~$\orbit$ is of the form $w(p_0)$ for some word~$w$ (which
is unique if we require that $w$~not end in $\sigma$ or~$\sigma^{-1}$).
Since $Q_1$~is finite, there is a number~$N$ such that each $x \in Q_1$ is
of the form~$w(p_0)$ with $w$~a word of length at most~$N$.  So, starting
with~$y$, we can step backward repeatedly as in the preceding paragraph
until we reach a point $y_* = w_*(p_0)$ where the word~$w_*$ begins
with more than~$N$ inverse generators ($\rho^{-1}$ or~$\sigma^{-1}$).
Then apply the first part of the claim to get an acceptable path (for the
sets $A_i^n$ and~$B_i^n$) from sone~$x \in D_3$ to~$y_*$.  All points on
this path are of the form~$w(p_0)$ where $w$~begins with more than~$N$
inverse generators, so none of them are in~$Q_1$.  Hence, this path is
also acceptable for the sets $\hat A_i$ and~$\hat B_i$.  By reversing
the backward steps taken from~$y$ to~$y_*$, we get an extension of this
path to an acceptable path from~$x$ to~$y$, as desired.
\QEd\par\smallskip

Let~$y$ be a point which is in~$Q_1$ or adjacent (via an edge of the
graph~$\orbit$) to a point of~$Q_1$, but is not in~$\hat A_1
\cup \hat A_2 \cup \hat A_3$.  Then, for any~$i$ such that
$y \notin \hat B_i$, we can apply the second part of Claim~3
to get an acceptable path from a point in~$D_3$ to~$y$ so that
$y$ is labeled~$i$.  None of the points on this path are in
$\hat A_1 \cup \hat A_2 \cup \hat A_3$ (by the same argument used earlier
to show that none of the points on the path~$P$ are in
$A_1^n \cup A_2^n \cup A_3^n$); in particular, $p_0$~is not
on the path.

Do this for each such point~$y$ and each~$i$ to get a finite collection
of acceptable paths.  Let $Q_2$ be the set of all points on these paths,
and let $\delta_4$ be the minimum distance from a point of~$Q_2 \cap D_3$
to the boundary of~$D_3$.  (So $\delta_3 - \delta_4$ is the
maximum distance from $p_0$ to a point in~$Q_2
\cap D_3$.)

By following the edges specified by the word
$w_0\sigma^{n_1}w_0w_0\sigma^{n_1}w_0$, we get a path in~$\orbit$ from
$w_0^{-1}\sigma^{-n_1}w_0^{-1}(p_0)$ to $w_0\sigma^{n_1}w_0(p_0)$.
Let $Q_3$ be the set of points on this path; note that $Q_0 \subset Q_3$.

Let $Q_*$ be the set of all points within at most $\hat M+1$ steps (in
the graph~$\orbit$) of a point in~$Q_1 \cup Q_2 \cup Q_3$.
So $Q_*$~is a finite subset of~$\orbit$.  Now choose $\delta_5>0$ so
small that:
\roster
\item"$\bullet$" if $x$ and~$y$ are distinct points of~$Q_*$,
then $d(x,y) > 2\delta_5$;
\item"$\bullet$" if $x \in Q_*$, then $B(x,\delta_5)$ lies on one side of
each of the sets $A_i^n$ and~$B_i^n$; and
\item"$\bullet$" $\delta_5 < \min(\delta_0 - \delta_1, \delta_2 -
\delta_3, \delta_4)/2$.
\endroster

We can now find $\delta_7$ and~$\delta_6$ such that $0 < \delta_7 <
\delta_6 < \delta_5$, $\delta_7$ and~$\delta_6$ are not in the (countable)
set $\{d(x,y) \colon x,y \in \orbit\}$, and $B(y',\delta_6) \subseteq
Z_n$, where $y'$~is the point in~$Z_n$ at the end of the acceptable
path~$P$.  We can also require~$\delta_6$ to be so small that the
following geometrical condition holds: if $C$~is a circle on~$S^2$
of radius~$\delta_1$, $a$~is a point on~$C$, $C'$~is the circle with
center~$a$ and radius~$\delta_3$, $b$~is an intersection point of
$C$ and~$C'$, and $x$~is a point on~$C'$ such that
$2\delta_5 - 2\delta_6 \le d(b,x) \le \delta_3$, then
the distance from~$x$ to~$C$ is at least~$2\delta_6$.

Construct the sets $A_i^{n+1}$ and~$B_i^{n+1}$ by replacing
each of the new points in $\hat A_i$ and~$\hat B_i$ with an open annulus
of inner radius~$\delta_7$ and outer radius~$\delta_6$.  That is, let
$$\align
A_i^{n+1} &= A_i^n \cup \bigcup \{ An(x,\delta_7,\delta_6) \colon
x \in \hat A_i \setminus A_i^n \} \qquad\text{and}\\
B_i^{n+1} &= B_i^n \cup \bigcup \{ An(x,\delta_7,\delta_6) \colon
x \in \hat B_i \setminus B_i^n \},
\endalign$$
where $An(x,\delta_7,\delta_6)$ is the open annulus $B(x,\delta_6)
\setminus \overline{B(x,\delta_7)}$.
(Note that if $x$~is a point in~$Q_1$ which is already in $A_i^n$,
then $An(x,\delta_7,\delta_6) \subseteq A_i^n$ because $\delta_6 <
\delta_5$; hence, we could have written ``$x \in Q_1 \cap \hat A_i$''
instead of ``$x \in \hat A_i \setminus A_i^n$'' above.  The same applies
to the $B$~sets.)
Clearly we have $A_i^n
\subseteq A_i^{n+1}$ and $B_i^n \subseteq B_i^{n+1}$.  The annulus
$An(y',\delta_7,\delta_6)$ will be a subset of~$Z_n$, so one of the
sets~$A_i^{n+1}$ meets $Z_n$.  It remains to verify that $A_i^{n+1}$
and~$B_i^{n+1}$ satisfy \ri1--\ri8.

The choice of $\delta_7$ and~$\delta_6$ ensures that $A_i^{n+1}$
and~$B_i^{n+1}$ satisfy~\ri1.

Since $\delta_6 < \delta_5$, the annuli around the points in~$Q_*$
are disjoint from each other.  Also, for each point $z \in Q_*$, the
disk $B(z,\delta_5)$ lies on one side of each set $A_i^n$ or~$B_i^n$.
Using this, we see that, if $x \in S^2$ is in one of the annuli
$An(z,\delta_7,\delta_6)$ for $z \in Q_*$, then $x \in A_i^{n+1}$ iff
$z \in \hat A_i$, and $x \in B_i^{n+1}$ iff $z \in \hat B_i$.  On the
other hand, if $x$ is not in any of the annuli $An(z,\delta_7,\delta_6)$
for $z \in Q_1$ [{\it sic}], then $x \in A_i^{n+1}$ iff $x \in A_i^n$
and then $x \in B_i^{n+1}$ iff $x \in B_i^n$.

It is now straightforward to verify \ri2--\ri5 for the sets $A_i^{n+1}$
and~$B_i^{n+1}$.  For instance, here is the proof that $\rho(A_1^{n+1})
\subseteq A_3^{n+1}$.  Let $x$ be a point in~$A_1^{n+1}$.  If $x \in
An(z,\delta_7,\delta_6)$ for some $z \in Q_*$ such that $\rho(z)$~is also
in~$Q_*$, then $x \in A_1^{n+1}$ implies $z \in \hat A_1$, which implies
$\rho(z) \in \hat A_3$, which implies $\rho(x) \in A_3^{n+1}$, because
$\rho(x) \in An(\rho(z),\delta_7,\delta_6)$.  On the other hand, if $x$~is
not in $An(z,\delta_7,\delta_6)$ for any such~$z$, then in particular $x
\notin An(z,\delta_7,\delta_6)$ for all $z \in Q_1$ (because $z \in Q_1$
implies $\rho(z) \in Q_*$), so $x \in A_1^{n+1}$ implies $x \in A_1^n$,
which implies $\rho(x) \in A_3^n$, which implies $\rho(x) \in A_3^{n+1}$.

To show that \ri8~holds for the new sets, suppose that an allowed
modification is made to a point $x \in \orbit$, and the corresponding
modification propagation algorithm is applied.  As long as the algorithm
does not reach a point in an annulus $An(z,\delta_7,\delta_6)$ for some
$z \in Q_1$, the points it reaches will be in~$A_i^{n+1}$ iff they are
in~$A_i^n$, and similarly for~$B_i^{n+1}$, so the algorithm will work
for the new sets exactly as it does for the old sets; hence, it will not
go more than $M$~steps from~$x$.  So suppose the algorithm does reach
a point~$x'$ (within $M+1$ steps of~$x$) in such an
annulus $An(z,\delta_7,\delta_6)$.  Consider what happens if we start
with the sets $\hat A_i$ and $\hat B_i$, modify the point~$z$ in the
way $x'$~was modified above, and apply the algorithm.  This algorithm
will terminate and will only modify points at most $\hat M$~steps from~$x'$
(although it may examine points one step farther away to see whether
they need modification).  For any word~$w$ of length at most $\hat M+1$, we
have $w(z) \in Q_*$, so $w(x') \in A_i^{n+1}$ iff $w(z) \in \hat A_i$
annd similarly for~$B_i^{n+1}$.  Therefore, the algorithm execution
starting at~$x'$ for $A_i^{n+1}$ and~$B_i^{n+1}$ will behave exactly
like the algorithm execution starting at~$z$ for $\hat A_i$ and~$\hat
B_i$, so it will go at most $\hat M$~steps from~$x'$.  (An exception
occurs when the execution for~$z$ reaches the ineligible point~$p_0$.
In this case, the corresponding point from the execution for~$x'$ is in
$An(p_0,\delta_7,\delta_6)$, which is included in~$A_1^{n+1}$.  So this
point is not reached by algorithm~A; it may be reached by algorithm~B,
but since it is in $B_2^{n+1} \cap B_3^{n+1}$ the algorithm will not
proceed any farther.  Hence, the algorithm execution for~$x'$ may go one
step farther than the algorithm execution for~$z$.)  Therefore, the full
algorithm starting at~$x$ will go no further than $M+\hat M+2$ steps from~$x$.

We next prove~\ri6 for the new sets.  Since $w_0\sigma^{n_1}w_0$
is an isometry of~$S^2$, we have
$$d((w_0\sigma^{n_1}w_0)^{-1}(p_0),p_0) = d(p_0,w_0\sigma^{n_1}w_0(p_0)).$$
So $w_0\sigma^{n_1}w_0(p_0)$ and $(w_0\sigma^{n_1}w_0)^{-1}(p_0)$ lie on
the same circle centered at~$p_0$; hence, we can find $n_2 > 0$ such
that $\sigma^{n_2}(w_0\sigma^{n_1}w_0(p_0))$ is at distance less
than~$\delta_7/2$ from $(w_0\sigma^{n_1}w_0)^{-1}(p_0)$.
Let $$\delta_8 = \delta_7 - d(\sigma^{n_2}w_0\sigma^{n_1}w_0(p_0),
(w_0\sigma^{n_1}w_0)^{-1}(p_0))$$ and $D_8 = B(p_0,\delta_8)$.
Let $P_8$ be the path starting from~$p_0$ with associated word
$w_0\sigma^{n_1}w_0\sigma^{n_2}w_0\sigma^{n_1}w_0$, labeled so that
the $w_0$~parts have labeling corresponding to~$P_0$ while all of
the extra $\sigma$~steps lead to points labeled~1.

We will show that, for each point~$x$ on~$P_8$, $B(x,\delta_8)$ is
disjoint from $B_1^{n+1} \cup B_2^{n+1} \cup B_3^{n+1}$.
It then follows immediately that $P_8$ is acceptable.
Since the last point on~$P_8$ is at distance $\delta_7-\delta_8$
from~$p_0$, and $\delta_7 - \delta_8 < \delta_7/2 < \delta_8$,
this last point is in~$D_8$; so we will have~\ri6 for the new sets.

We first show that each point~$x$ on~$P_8$ is within distance
$\delta_0-\delta_8$ of some point~$x'$ on~$P_0$; it will then follow
that $B(x,\delta_8) \subseteq B(x',\delta_0)$, so, by the choice
of~$P_0$, $B(x,\delta_8)$ must be disjoint from $B_1^n \cup B_2^n \cup
B_3^n$.  The initial part of~$P_8$ given by the first~$w_0$ is
$P_0$~itself.  The point $w_0(p_0)$ is at distance~$\delta_1$ (which is
less than $\delta_0-\delta_8$ because $\delta_8 < \delta_5 < \delta_0 -
\delta_1$) from $p_0$; the same applies to the following points on~$P_8$
up to $\sigma^{n_1}w_0(p_0)$.  The points on~$P_8$ coming from the
second~$w_0$ are at distance~$\delta_1$ from the corresponding points
on~$P_0$.  The last of these is at distance~$\delta_3$ (again less than
$\delta_0 - \delta_8$) from~$p_0$ and so are the following points given
by the~$\sigma^{n_2}$, so the third~$w_0$ gives points at
distance~$\delta_3$ from the corresponding points of~$P_0$.  The last of
these (call it~$x$) is at distance $\delta_7-\delta_8$ from
$(w_0\sigma^{n_1})^{-1}(p_0)$, which is at distance~$\delta_1$
from~$p_0$; so $d(x,p_0) \le \delta_1 + \delta_7 - \delta_8$, which is
less than $\delta_0 - \delta_8$ because $\delta_7 < \delta_5 < \delta_0
- \delta_1$.  The same applies to $\sigma^m(x)$ for $m \le n_1$, and now
the final $w_0$~segment gives points within $\delta_0-\delta_8$ of
points on~$P_0$, as desired.

We must now show that the sets $B(x,\delta_8)$ for $x$ on~$P_8$ do not
contain any of the new points added to~$B_i^n$ to get $B_i^{n+1}$;
to do this, it will suffice to show that $B(x,\delta_8)$ is
disjoint from all annuli $An(x',\delta_7,\delta_6)$ for $x' \in Q_1$.  If
$x$~is on the initial part of~$P_0$ given by the first
$w_0\sigma^{n_1}w_0$, then $x$~is in~$Q_3$.  Hence, any $x' \ne x$
in~$Q_1$ is at distance at least~$2\delta_5$ from~$x$, so $B(x,\delta_8)
\cap An(x',\delta_7,\delta_6) = \nullset$.  The point $x$~itself could
be in~$Q_1$, but we would have $B(x,\delta_8) \cap
An(x,\delta_7,\delta_6) = \nullset$ because $\delta_8 < \delta_7$.
So $B(x,\delta_8)$ does not contain any of the new points.

The next part of~$P_8$ (given by the middle~$\sigma^{n_2}$) consists of
points~$x$ at distance~$\delta_3$ from~$p_0$.  Since $2\delta_5$~is less
than~$\delta_3$ and also less than~$\delta_2-\delta_3$, and any point
of~$Q_1$ must either be~$p_0$ or at distance greater than~$\delta_2$
from~$p_0$ (by the definition of~$\delta_2$), it must be that any
point~$x$ on this middle part of~$P_8$ must be at distance at
least~$2\delta_5$ from any point $x' \in Q_1$.  We therefore get
$B(x,\delta_8) \cap An(x',\delta_7,\delta_6) = \nullset$ again.

Each point~$x$ on the final part of~$P_8$ (given by the remaining
$w_0\sigma^{n_1}w_0$) is at distance exactly $\delta_7 - \delta_8$ from
a point $y \in Q_3$.  Any point $x' \ne y$ in~$Q_1$ is at distance
at least~$2\delta_5$ from~$y$, so $B(x,\delta_8)
\cap An(x',\delta_7,\delta_6) = \nullset$ because
$(\delta_7-\delta_8)+\delta_8+\delta_6 < 2\delta_5$.  And even if
$y$~itself is in~$Q_1$, we have $B(x,\delta_8)
\cap An(y,\delta_7,\delta_6) = \nullset$ because
$d(x,y) = \delta_7-\delta_8$.  This completes the proof of~\ri6 for the
new sets.

For~\ri7 we will use one more claim:

\procl{Claim 4} Let~$x$ be a point in~$D_3$ which is not in or on the
boundary of the annulus $An(p_0,\delta_7,\delta_6)$.  Then there is an
acceptable (for the sets $A_i^{n+1}$ and~$B_i^{n+1}$) path from~$p_0$ to
a point~$x' \in D_3$ such that $d(x',x) < \delta_8$, $x'$ is labeled~1,
and, for every point~$a$ on the path, $B(a,\delta_8)$ is disjoint from
$B_1^{n+1} \cup B_2^{n+1} \cup B_3^{n+1}$.

\procl{Proof}
First, suppose $x$~is inside the annulus, so $d(x,p_0) < \delta_7$.
Let $z$ be the last point on~$P_8$.  Since $\delta_7/2 < \delta_8 <
\delta_7$ and $d(z,p_0) = \delta_7 - \delta_8$, we can cover
the entire disk $B(p_0,\delta_7)$ by rotating the disk~$B(z,\delta_8)$
around the point~$p_0$.  But the points $\sigma^m(z)$ for
$m=1,2,3,\dotsc$ are dense in the circle with center~$p_0$ and
radius $\delta_8 - \delta_7$; it follows that the open disks
$B(\sigma^m(z),\delta_8)$ for $m>0$ cover the open disk
$B(p_0,\delta_7)$.  Therefore, we can find a positive integer~$m$ such
that $d(\sigma^m(z),x) < \delta_8$.  Let $x' = \sigma^m(z)$; the
path~$P_8$, followed by the $m$~steps from $z$ to $\sigma^m(z)$, gives
the desired path from~$p_0$ to~$x'$.

Now suppose $x$~is outside the annulus, so $\delta_6 < d(x,p_0) <
\delta_3$.  The points $\sigma^kw_0\sigma^{n_1}w_0(p_0)$ for $k>0$ are
dense in the circle of radius~$\delta_3$ around~$p_0$, and the points
$w_0\sigma^kw_0\sigma^{n_1}w_0(p_0)$ are dense in the corresponding
circle around $w_0(p_0)$.  Hence, we can choose $k > 0$ so that, if $y =
w_0\sigma^kw_0\sigma^{n_1}w_0(p_0)$, then $y$ lies inside the circle with
center~$p_0$ and radius~$\delta_1$, and
$$\max(\delta_6+\delta_8,d(x,p_0)-\delta_8) < d(y,\sigma^{-n_1}w_0^{-1}(p_0))
< \min(\delta_3,d(x,p_0)+\delta_8).$$
This will imply $d(x,p_0)-\delta_8 < d(w_0\sigma^{n_1}(y),p_0) <
d(x,p_0)+\delta_8$; it follows that there is a positive number~$m$ such
that $d(\sigma^mw_0\sigma^{n_1}(y),x) < \delta_8$.  Fix such an~$m$,
and let $x' = \sigma^mw_0\sigma^{n_1}(y)$; we will see that the path
from~$p_0$ to~$x'$ given by the word
$\sigma^mw_0\sigma^{n_1}w_0\sigma^kw_0\sigma^{n_1}w_0$
has the desired properties.

We must see that, for every point~$a$ on the path, $B(a,\delta_8)$ is
disjoint from $B_1^{n+1} \cup B_2^{n+1} \cup B_3^{n+1}$.  (Given this,
the labeling where the $w_0$~segments are labeled like~$P_0$ and the
extra $\sigma$~steps lead to points labeled~1 will make this path
acceptable.)  The argument is similar to that for~$P_8$.
Each of the four $w_0$~segments ends up at a point within distance
$\delta_0-\delta_8$ of~$p_0$ (the distances are respectively $\delta_1$,
$\delta_3$, less than~$\delta_1$, and less than~$\delta_3$), so, as
for~$P_8$, each point on this path is within distance
$\delta_0-\delta_8$ of a point in~$P_0$; it follows that the disks
$B(a,\delta_8)$ do not meet $B_1^n \cup B_2^n \cup B_3^n$.
It remains to show that these disks do not meet any of the annuli
$An(z,\delta_7,\delta_6)$ for $z \in Q_1$.

The points on the initial segment of the path given by the word
$w_0\sigma^{n_1}w_0$ are also on the path~$P_8$, so they have already
been taken care of.  The next segment (given by~$\sigma^k$) consists of
points at distance~$\delta_3$ from~$p_0$; these are handled by the same
argument as for the middle segment of~$P_8$.  Then comes the third
$w_0$~segment of the path; each point of this segment is at
distance~$\delta_3$ from the corresponding point of~$P_0$, so the same
argument using the definition of~$\delta_2$ applies to handle these
points.

We have now reached the point~$y$.
Let $s = d(y,\sigma^{-n_1}w_0^{-1}(p_0))$.  Each point~$a$ on the path
from~$y$ to the endpoint~$x'$ is at distance~$s$ from a point $a'\in
Q_3$.  Since $s > \delta_6+\delta_8$, $B(a,\delta_8)$ cannot intersect
$An(a',\delta_7,\delta_6)$.  If $z$~is a point of~$Q_1$ other than~$a'$,
then $d(z,a') \ge 2\delta_5$.  Hence, if $s < 2\delta_5 - 2\delta_6$,
then $B(a,\delta_8)$ cannot meet $An(z,\delta_7,\delta_6)$, so we are
done.

So assume $s \ge 2\delta_5 - 2\delta_6$.  Then the extra geometrical
condition imposed on~$\delta_6$ implies that the distance from~$y$ to
the circle with center~$p_0$ and radius~$\delta_1$ is at
least~$2\delta_6$.  In other words, $\delta_1 - \delta_3 \le d(y,p_0)
\le \delta_1 - 2\delta_6$.  The same applies to~$\sigma^j(y)$, since
$d(\sigma^j(y),p_0) = d(y,p_0)$.  Now, the definition of~$\delta_2$ implies
that any point $z\in Q_1$ must satisfy either $d(z,p_0) \ge \delta_1$
or $d(z,p_0) \le \delta_1-\delta_2$.  Since $\delta_6 < \delta_5 <
(\delta_2-\delta_3)/2$, we must have $d(z,\sigma^j(y)) \ge 2\delta_6$
for any~$j$ and any~$z \in Q_1$; it follows that
$B(\sigma^j(y),\delta_8)$ and $An(z,\delta_7,\delta_6)$ are disjoint.

This takes care of the segment from~$y$ to~$\sigma^{n_1}(y)$.  For the
next segment (the fourth $w_0$~segment), each point~$a$ is at distance
$d(\sigma^{n_1}(y),p_0)$ from a point on~$P_0$; so the same argument
used for~$\sigma^{n_1}(y)$ will handle~$a$.

The final segment of the path consists of points~$a$ at distance~$s$
from~$p_0$.  Since $s > \delta_6 + \delta_8$, $B(a,\delta_8)$ cannot
intersect $An(p_0,\delta_7,\delta_6)$.  Any point $z \in Q_1$ other
than~$p_0$ is at distance at least~$\delta_2$ from~$p_0$; since $s \le
\delta_3 < \delta_2 - 2\delta_5$, $B(a,\delta_8)$ cannot intersect
$An(z,\delta_7,\delta_6)$ either.  So the path has the desired
properties.
\QEd\par\smallskip

Now, to prove~\ri7, let~$p\ne p_0$ be in~$\orbit$ but not in $A_1^{n+1}
\cup A_2^{n+1} \cup A_3^{n+1}$, and let~$i$ be such that $p \notin
B_i^{n+1}$.  By Claim~3, there is an acceptable (for the sets $A_j^n$
and~$B_j^n$) path from some point in~$D_3$ to~$p$ which gives~$p$ the
label~$i$.  If there is no point on this path which is in any of the
annuli $An(z,\delta_7,\delta_6)$ for $z \in Q_1$, then the path is also
acceptable for the sets $A_i^{n+1}$ and~$B_i^{n+1}$. On the other hand,
if there is such a point, then let $p'$ be the last point on the path
which is in an annulus $An(z,\delta_7,\delta_6)$ where $z$~is in {\sl or
adjacent to}~$Q_1$.  Then either $p' = p$ or $z$ is adjacent to rather
than in~$Q_1$; in either case, if $i'$~is the label of~$p'$, then we get
$p' \notin B_{i'}^{n+1}$.  (This is given if $p' = p$; if $p' \ne p$,
then since $z \notin Q_1$ we get $p' \in B_{i'}^{n+1}$ iff $p' \in
B_{i'}^n$, and the latter does not hold because the path is acceptable.)
It follows that $z \notin \hat B_{i'}$.  So, as described just after
Claim~3, we already selected an acceptable path (for the sets $\hat A_j$
and~$\hat B_j$) from some point in~$D_3$ to~$z$ so that $z$ got
label~$i'$, and the set~$Q_2$
includes the points on this path.  Let~$w$ be the associated word for
this latter path; then there is a corresponding path from $w^{-1}(p')$
to~$p'$ with the same labeling, and this path will be acceptable for the
sets $A_i^{n+1}$ and~$B_i^{n+1}$ (because, if $b'$~is a point on this
path and $b$~is the corresponding point on the path to~$z$, then
$b \in Q_*$, so we have $b' \in B_j^{n+1}$ iff $b \in \hat B_j$).
The continuation of this from~$p'$ to~$p$ as on the original path is
also acceptable for these sets, because no point from~$p'$ to~$p$ is
in~$Q_1$ (unless $p'=p$).  Note that $w^{-1}(p')$ is within
distance~$\delta_6$ of~$w^{-1}(z)$, which is within distance
$\delta_3-\delta_4$
of~$p_0$; since $\delta_6 < \delta_4$, we have $d(w^{-1}(p'),p_0) <
\delta_3$, so $w^{-1}(p') \in D_3$.

Hence, in any case, there is an acceptable (for the sets $A_j^{n+1}$
and~$B_j^{n+1}$) path from some point $x \in D_3$ to~$p$ so that $p$~is
labeled~$i$.  Note that $x$~is not in any of the sets $A_j^{n+1}$.
(As argued previously, if a point on the acceptable path were in one of
these sets, then the next point on the path would be also, and so on all
the way to~$p$; but $p$~is not in any of these sets.)  So $x$~cannot be
in the annulus $An(p_0,\delta_7,\delta_6)$, which is included in
$A_1^{n+1}$.  Also, $x$~cannot lie on the boundary of this annulus,
since this boundary contains no points in~$\orbit$.  So we can apply
Claim~4 to get an acceptable path from~$p_0$ to~$x'$ with the properties
listed in that Claim.  Let $w_4$ be the associated word for this new
path.  Then there is a corresponding path from~$w_4^{-1}(x)$ to~$x$.
Each point on this path is at distance less than~$\delta_8$ from the
corresponding point on the path from~$p_0$ to~$x'$, and hence is not in
$B_1^{n+1} \cup B_2^{n+1} \cup B_3^{n+1}$.  Therefore, this new path is
acceptable, and combining it with the path from~$x$ to~$p$ gives an
acceptable path from~$w_4^{-1}(x)$ (which is in~$D_8$) to~$p$.
So \ri7~holds for the sets $A_i^{n+1}$ and~$B_i^{n+1}$.

Therefore, the new sets $A_i^{n+1}$ and~$B_i^{n+1}$ satisfy \ri1--\ri8.
This completes the construction and the proof of the theorem.
\QED\enddemo

The proof of Theorem~5.5 can be modified to yield pairwise disjoint
open subsets $A_1,A_2,A_3,A_4$ of the sphere with dense union such that
$\rho(A_4) = A_1$, $\rho(A_1) = A_3$, and $\sigma(A_1 \cup A_3) = A_1 \cup
A_2$, where $\rho$ and~$\sigma$ are given free rotations.  (There will
now be four sets $A_i^n$ and four sets~$B_i^n$.  Hypotheses \ri1,~\ri2,
and~\ri8 are unchanged, and \ri3,~\ri4, and~\ri6 have the obvious changes.
Hypothesis~\ri5 now states that $\rho(B_1^n) = B_3^n$, $\rho(B_4^n) =
B_1^n$, $\sigma(B_1^n \cap B_3^n) = B_1^n \cap B_2^n$, $\rho(B_2^n \cap
B_3^n) = B_2^n \cap B_4^n$, and $\sigma(B_2^n \cap B_4^n) = B_3^n \cap
B_4^n$.  Hypothesis~\ri7 has the obvious ``${}\cup A_4^n$'' added, and
also restricts~$i$ to the values $1,2,3$.  Acceptable paths will still
only use the labels 1,~2, and~3; if we need to meet an open set~$Z_n$
by adding a point of it to~$\hat A_4$, we will do so by adding a point
of~$\rho(Z_n)$ to~$\hat A_1$.  The needed changes to modification
propagation algorithm~B are straightforward, read off directly
from the new~\ri5.  For algorithm~A, we add a clause that, if
$x$~has been added to~$\hat A_1$, then we should add~$\rho^{-1}(x)$
to~$\hat A_4$.  The rest of the proof goes through as before.)

Hence, one gets open subsets of the sphere with dense union satisfying
(via free rotations) the system $$A_1 \cong A_3, \qquad A_1 \cup A_2
\cong A_1 \cup A_3 \cong A_1 \cup A_4.$$ (The congruence $A_1 \cong
A_4$ is also satisfied.)  This is of interest because one can show
that this system cannot be satisfied by finite subsets of a free group
(unless they are all empty).  If we had such finite subsets, then all
four of them would have to have the same number of elements (hence,
they would all be nonempty).  Now, by Proposition~5.4, there would be
group elements $\sigma$ and~$\sigma'$ such that $\sigma(A_3) = A_1$,
$\sigma(A_1) = A_2$, $\sigma'(A_3) = A_1$, and $\sigma'(A_1) = A_4$.
Since $\sigma(A_1) \ne \sigma'(A_1)$, $\sigma \ne \sigma'$.  But now
$\sigma'(\sigma^{-1}(A_1)) = A_1$, which is impossible for a nonempty
finite set $A_1$ (given any element of it, we could apply $\sigma'
\circ \sigma^{-1}$ repeatedly to get infinitely many elements of it).

If we allow arbitrary isometries rather than just free rotations,
then all of the systems of congruences seen so far in this section
are very simply satisfiable by open subsets of the sphere with dense
union, because they are all subsystems of $\UNC s$ for some $s \le 4$.
However, another modification of the proof of Theorem~5.5 yields pairwise
disjoint open subsets $A_1,\dots,A_6$ of~$S^2$ with dense union such
that $\rho(A_6)=A_5$, $\rho(A_5)=A_4$, $\rho(A_4)=A_1$, $\rho(A_1)=A_3$,
and $\sigma(A_1 \cup A_3) = A_1 \cup A_2$, where $\rho$ and~$\sigma$
are given free rotations.  These sets satisfy the system
$$\gather
A_1 \cong A_3 \cong A_4 \cong A_5 \cong A_6,\\
A_2 \cup A_1 \cong A_1 \cup A_3 \cong A_1 \cup A_4 \cong A_4 \cup A_5 \cong
A_5 \cup A_6,\\
A_3 \cup A_1 \cup A_4 \cong A_1 \cup A_4 \cup A_5 \cong A_4 \cup A_5 \cup
A_6,\\
A_3 \cup A_1 \cup A_4 \cup A_5 \cong A_1 \cup A_4 \cup A_5 \cup A_6.
\endgather$$
There is no obvious simpler proof that this system is satisfiable
by open sets with dense union even in the arbitrary-isometries case,
since $\UNC6$ is not known to be satisfiable on the sphere by such sets.

\head 6. The various solvability properties \endhead

In this paper, we have considered a number of variations of the question of
whether a system of congruences can be satisfied nontrivially
(i.e., by sets which are not all empty), depending
on what kind of subsets we are allowing (open or finite),
what space they are subsets of, and which isometries can be used
to witness the congruences.  Here is a list of these variations:

\roster
\item"OSF:" open subsets of the sphere, using free rotations
\item"FSF:" finite subsets of the sphere, using free rotations
\item"DSF:" open subsets of the sphere with dense union, using free rotations
\item"OSI:" open subsets of the sphere, using any isometries
\item"FSI:" finite subsets of the sphere, using any isometries
\item"DSI:" open subsets of the sphere with dense union, using any isometries
\item"OPS:" open subsets of any suitable Polish space
\item"FPS:" finite subsets of any suitable Polish space
\item"DPS:" open subsets of any suitable Polish space with dense union
\item"FFG:" finite subsets of a free group
\item"PFG:" finite subsets of a free group with connected prime union
\item"FFQ:" finite subsets of the cosets of a pure cyclic subgroup in
a free group [Theorem 3.2 (III)]
\endroster

And here are a few more properties a system of congruences can have
that are relevant in characterizing the satisfiablility of the system:
\roster
\item"w:" the system is weak
\item"nc:" the system is numerically consistent
\item"c:" the system is consistent
\endroster
Recall that the weak systems are those that can be satisfied by
a partition of a sphere into arbitrary pieces, using a free group of
rotations to witness the congruences; if one also wants the pieces
to be nonmeager sets with the property of Baire, then it is precisely
the weak consistent systems that have solutions.  (In both cases the
requirement of weakness can be dropped if one allows arbitrary isometries
to witness the congruences.)

\plotfigscale 7.2
\midinsert
\centerline{%
\plotfigbegin
   \plotvskip 42
   \plothphant 44
   \plottext{OSF}
   \plot 22 21
   \plottext{FSF}
   \plot 32 21
   \plottext{DSF}
   \plot 12 21
   \plottext{OSI}
   \plot 22 11
   \plottext{FSI}
   \plot 32 11
   \plottext{DSI}
   \plot 12 11
   \plottext{OPS}
   \plot 22 31
   \plottext{FPS}
   \plot 32 31
   \plottext{DPS}
   \plot 12 31
   \plottext{FFG}
   \plot 32 41
   \plottext{PFG}
   \plot 12 41
   \plottext{FFQ}
   \plot 42 21
   \plottext{w}
   \plot 2 21
   \plottext{nc}
   \plot 22 1
   \plottext{c}
   \plot 12 1
   \plotcentered{3.2}
   \plot 27 22.3
   \plot 27 19.7
   \plot 37 22.3
   \plot 37 19.7
   \plotcentered{4.2}
   \plot 27 32
   \plotcentered{6.1}
   \plot 27 12.3
   \plot 27 9.7
   \plotcentered{\S2}
   \plot 7 22
   \plot 17 2
   {\plotrotate 270
      \plot 23 6
      \plotcentered{4.1}
      \plot 33.3 36
      \plot 30.7 36
      \plotcentered{4.3}
      \plot 13 36 }
   \plotlinewidth 0.12\plotunity
   \plotmove 14 41
   \plotarrow 30 41
   \plotmove 12 40
   \plotarrow 12 32
   \plotmove 31.7 40
   \plotarrow 31.7 32
   \plotmove 32.3 32
   \plotarrow 32.3 40
   \plotmove 33.5 40
   \plotarrow 41.5 22
   \plotmove 14 31
   \plotarrow 20 31
   \plotmove 30 31
   \plotarrow 24 31
   \plotmove 12 30
   \plotarrow 12 22
   \plotmove 22 30
   \plotarrow 22 22
   \plotmove 32 30
   \plotarrow 32 22
   \plotmove 10 21
   \plotarrow 3 21
   \plotmove 14 21
   \plotarrow 20 21
   \plotmove 30 21.3
   \plotarrow 24 21.3
   \plotmove 40 21.3
   \plotarrow 34 21.3
   \plotmove 12 20
   \plotarrow 12 12
   \plotmove 22 20
   \plotarrow 22 12
   \plotmove 32 20
   \plotarrow 32 12
   \plotmove 14 11
   \plotarrow 20 11
   \plotmove 30 11.3
   \plotarrow 24 11.3
   \plotmove 22 10
   \plotarrow 22 1.8
   \plotmove 20.5 1
   \plotarrow 13 1
   \plotlinewidth 0.01\plotunity
   \plotmove 24 20.7
   \plotarrow 30 20.7
   \plotmove 34 20.7
   \plotarrow 40 20.7
   \plotmove 24 10.7
   \plotarrow 30 10.7
\plotfigend}
\botcaption{Figure 6.1}
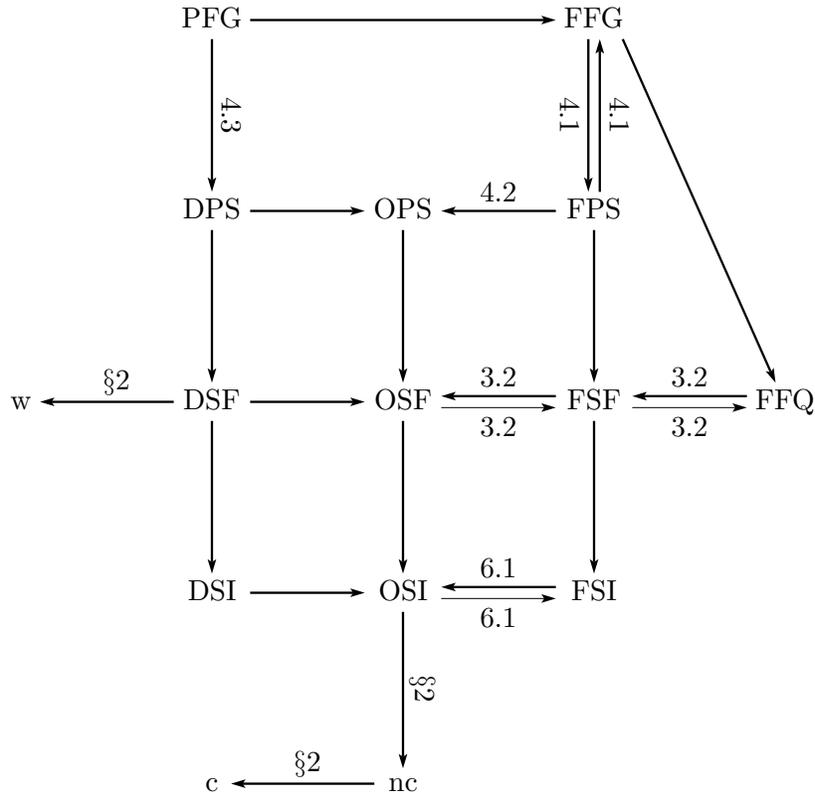
Implications between the various satisfiability properties.
\endcaption
\endinsert

Figure~6.1 shows the known implications between these properties.
Most of them were proved earlier in the paper (as indicated in the
figure), or are trivial; the remaining ones are given by:

\proclaim{Proposition 6.1} If arbitrary isometries can be used to witness
the congruences, then a system of congruences is satisfiable by open
subsets of the sphere (not all empty) if and only if the system is
satisfiable by finite subsets of the sphere (not all empty). \endproclaim

\demo{Proof}
The right-to-left implication is proved in the same way as in
Theorem~3.2 --- replace the points with identical small open disks.

For the other direction, suppose we have open sets~$A_j$ satisfying the
congruences.  Choose a connected component~$C$ of one of the sets~$A_k$,
and let $G$ be the stabilizer group of~$C$ (i.e., the set of all
isometries~$g$ of the sphere such that $g(C) = C$).  If $C$~is the
entire sphere~$S^2$, then only one of the sets~$A_j$ is nonempty, so we
can get finite sets satisfying the congruences by making just that
one of the finite sets nonempty.  So assume $C$~is not all of~$S^2$.

The isometry group~$\isomgroup$ of~$S^2$ is a compact group
under the maximum-distance metric $d_\isomgroup(g,h) = \max_{x \in S^2}
d(g(x),h(x))$.  Since $C$~is open, the group~$G$ must be a closed subgroup
of~$\isomgroup$.  To see this, suppose $g \in \isomgroup \setminus G$;
then $g(C) \ne C$, so there must be a point~$x$ such that either $x
\in g(C) \setminus C$ or $x \in C \setminus g(C)$.  In the former case,
for all $g'$ sufficiently close to~$g$ we have $x \in g'(C) \setminus
C$; in the latter case, for all $g'$ sufficiently close to~$g$ we have
$g'(g^{-1}(x)) \in C \setminus g'(C)$.  So the complement of~$G$ is
open in~$\isomgroup$.

If $G$~is finite, choose a point $z \in C$ and let $Z$ be the $G$-orbit
of~$z$.  Then $Z$ is a finite subset of~$C$ which is fixed under
any isometry which fixes~$C$.

If $G$~is infinite, then, since $G$~is closed in a compact group and hence
compact, we can choose a sequence of distinct members $h_n$ of~$G$
converging to some~$h \in G$.  Let $g_n = h^{-1} h_n$; then we
have $g_n \in G$ and the isometries~$g_n$ are distinct and converge to the
identity isometry.  Any isometry close to the identity must be
orientation-perserving, so we may assume that all of the isometries~$g_n$
are non-identity rotations.  Let $\axis_n$ be the axis of~$g_n$; by moving to
a subsequence if necessary, we may assume that the axes~$\axis_n$
converge to an axis~$\axis$.  Now, for large~$n$, the rotation~$g_n$ is
close to the identity, so its order is large if not infinite; hence, the
powers of~$g_n$ come close to all rotations around axis~$\axis_n$.
Therefore, if $\rho$~is any rotation around the limiting axis~$\axis$,
then $\rho$~can be approximated arbitrarily well by a power of~$g_n$ for
a sufficiently large~$n$, so $\rho$~is in the closure of~$G$, which is~$G$.

Thus, if~$G$ is infinite, then there is an axis~$\axis$ such that
all rotations around~$\axis$ are in~$G$; this means that $C$ must be
a disk or annulus centered on~$\axis$.  Note that there can only be one
such axis, since $C$~is a nonempty proper open subset of~$S^2$.
Also, $C$~must be symmetric under reflections of~$S^2$ which
leave the points of~$\axis$ fixed.  If $C$~is not symmetric under
reflections which reverse~$\axis$, let $Z$~be a set containing just
one point, one of the two intersections of~$\axis$ with~$S^2$;
if $C$ is symmetric under such reflections, let $Z$~be the set
comprising both of these intersections.  So $Z$ is a finite subset
of~$S^2$ (which need not be included in~$C$ in this case), and
the stabilizer group of~$Z$ is exactly the same as that of~$C$,
namely~$G$.

We now define a function~$F$ whose domain is the set of all
components of the sets~$A_j$ which are congruent to~$C$.
(Since these components all have the same positive measure and are
disjoint from each other, there are only finitely many of them.)
Given such a component~$C'$, let $g$~be an isometry such that
$g(C) = C'$, and define~$F(C')$ to be $g(Z)$.  Then $F(C')$~is
well-defined, because if $h$~is another isometry such that
$h(C) = C'$, then $h^{-1}(g(C)) = C$, so $h^{-1} \circ g \in G$,
so $h^{-1}(g(Z)) = Z$, so $g(Z) = h(Z)$.

If $C_1$ and~$C_2$ are in the domain of~$F$ and $h(C_1) = C_2$,
then we have $h(F(C_1)) = F(C_2)$.  To see this, fix isometries
$g_1$ and~$g_2$ such that $g_1(C) = C_1$ and $g_2(C) = C_2$.
Then $g_2^{-1}(h(g_1(C))) = C$, so $g_2^{-1}(h(g_1(Z))) = Z$,
so $h(F(C_1)) = h(g_1(Z)) = g_2(Z) = F(C_2)$.

Next, we note that, if $C_1$ and~$C_2$ are distinct members of the
domain of~$F$, then $F(C_1)$ and~$F(C_2)$ are disjoint.
In the case that $G$~is finite, this follows from the fact that
$C_1$ and~$C_2$ are disjoint and $F(C_i) \subseteq C_i$, $i=1,2$.
If $G$~is infinite, then $F(C_1)$ and~$F(C_2)$ are either both single
points or both pairs of antipodal points, so, if they are not disjoint,
then they coincide.  But if we have $F(C_1) = F(C_2)$ where
$C_i = g_i(C)$ for $i=1,2$, then we get $g_1(Z) = g_2(Z)$,
so $g_2^{-1}(g_1(Z)) = Z$, so $g_2^{-1}\circ g_1$ is in the
stabilizer group of~$Z$, which is~$G$ in this case; hence,
$g_2^{-1}(g_1(C)) = C$, so $C_1 = g_1(C) = g_2(C) = C_2$.
Hence, if $C_1 \ne C_2$, then $F(C_1)$ and~$F(C_2)$ must be disjoint.

Now we can define finite sets~$B_j$ as follows: for each~$j$,
let $B_j$ be the union of all of the sets~$F(C')$ where $C'$~is
a component of~$A_j$ congruent to~$C$.  Then the sets~$B_j$ are
pairwise disjoint and not all empty (one of them includes~$F(C)$).
And the fact that $h(F(C_1)) = F(C_2)$ whenever $h(C_1) = C_2$
implies that any congruence satisfied by the sets~$A_j$ is also satisfied
by the sets~$B_j$, using the same isometry.  So the sets~$B_j$
are finite sets satisfying the given system of congruences.
\QED\enddemo

In most cases, the proofs of the satisfiability implications in Figure~6.1
actually give stronger implications: if a system of congruences
is satisfiable in the first context by sets that are all nonempty
(not just ``not all empty''), then it is satisfiable in the second
context by sets that are all nonempty.  The three exceptions are
shown in the figure using lighter arrows.  The open-to-finite parts
of the proofs of Theorem~3.2 and Proposition~6.1 only ensure that
some of the finite sets satisfying the congruences are nonempty, even if
all of the given open sets were nonempty; the same thing happens in
the proof that (II) implies (III) in Theorem~3.2.  It is not known whether
one can give modified proofs that would yield finite sets satisfying
the congruences that are all nonempty.
(The implication $\text{DSF}\to\text{w}$ in Figure~6.1 is to be read
in the usual way: if a system is satisfiable in case~DSF using sets
which are not all empty, then the system must be weak.  However,
the implication $\text{OSI}\to\text{nc}$ is not quite that strong:
if a system is satisfiable in case~OSI using sets which are all
nonempty, then the system is numerically consistent.  If only some of the
sets are nonempty, then all one can conclude is that the given
system can be made numerically consistent by deleting zero or more
of the sets mentioned in it.)

We have seen a number of examples of systems of congruences which can
be used to show that various implications in Figure~6.1 are
not reversible.  Here is a summary of these examples:

The system $A_1 \cong A_1 \cup A_3 \cup A_4$, $A_3 \cong A_1 \cup A_2
\cup A_3$ used in Wagon's presentation~\cite{\Wagon} of Robinson's
results (Robinson~\cite{\Robinson} actually used a different system)
is weak but not consistent, and hence not satisfiable by open or finite
sets in any of the cases listed here.  The system $A_1 \cong A_2 \cong
A_3 \cong A_4 \cong\nobreak A_5$, $A_1 \cup A_2 \cong A_1 \cup A_3 \cup A_4$
from section~2 is weak and consistent, but not numerically consistent,
and hence also not satisfiable in any of these cases.

The trivial system $A_1 \cong A_2$ is not weak, but it is satisfiable
in all the cases not shown in Figure~6.1 as implying weakness
(i.e., it is satisfiable in cases DSI and FFG).

The system $A_1 \cup A_2 \cong A_1 \cup A_3 \cong A_2 \cup A_3$
from Theorem~3.1 is weak and is satisfiable in cases DSI and~FSI,
but is not satisfiable in case~OSF.

The system $A_1 \cong A_2 \cong A_3$, $A_1 \cup A_2 \cong A_1 \cup A_3$
from Theorem~5.1 is weak and is satisfiable in cases DSI and~FFG,
but is not satisfiable in case~DSF.

The system $A_1 \cong A_3$, $A_1 \cup A_2 \cong A_1 \cup A_3$ from Theorem~5.5
is satisfiable in cases DSF and~FFG, but not in case~PFG;
it is not known whether this system is satisfiable in case~DPS.

The system $A_1 \cong A_3$,
$A_1 \cup A_2 \cong A_1 \cup A_3 \cong A_1 \cup A_4$
given after Theorem~5.5 is satisfiable in case~DSF, but not in case~FFG;
it is not known whether this system is satisfiable in case~DPS or
case~OPS.

This leaves a few implications in Figure~6.1 which may or may not be
reversible: $\text{OSI}\to\text{nc}$ (the system~$\UNC6$ may be a
counterexample here), $\text{DSI}\to\text{OSI}$, $\text{FPS}\to\text{OPS}$,
and $\text{PFG}\to\text{DPS}\to\text{DSF}$ (the system from
Theorem~5.5 shows that these last two implications cannot both be
reversible).

\head 7. Completeness of congruence deduction rules \endhead

As noted early in section~2, a given system of congruences
on sets $A_1,A_2,\dots,A_r$ can imply
other congruences, because congruence must be an equivalence relation
(reflexive, symmetric, and transitive).  Also, if we are considering
the case where the sets~$A_i$ are required to
form a partition of the space in question, then one can also
use the complementation rule to deduce new congruences from old ones.
One can ask whether this set of rules is complete, in the sense
that any congruence which necessarily follows from a given system
of congruences is in fact deducible by these rules alone.
(I thank Harvey Friedman for bringing up this question.)

If we allow improper congruences in the system, then the answer
is no.  For instance, if the improper congruence
$A_1 \cong \nullset$ is satisfied, then the congruence
$A_1 \cup A_2 \cong A_2$ must also be satisfied, but this cannot
be deduced from the above rules (if $r>2$).  Similarly, if
if $r=3$ and the sets are required to form a partition, then
the improper congruence $A_1 \cup A_2 \cup A_3 \cong A_2 \cup A_3$ implies
the congruence $A_1 \cup A_2 \cong A_2$ (because it forces~$A_1$ to
be empty), and again one cannot deduce this by the given rules.

However, if we restrict ourselves to proper congruences, then the
answer is yes:

\proclaim{Theorem 7.1}
If one has a system of proper congruences and an additional congruence
which is not deducible from the system by the equivalence relation rules,
then one can find a suitable space and nonempty open subsets of that space
which satisfy the system of congruences but not the additional congruence.
If the additional congruence is not deducible from the system using
the equivalence relation rules and the complementation rule, then the
open subsets of the suitable space can be taken to form a partition of
the space.
\endproclaim

\demo{Proof}
The suitable space we will use is the discrete space $F \times \N$,
where $\N$ is the set of natural numbers and $F$ is a free group on $m$
generators $f_1,\dots,f_m$ (here $m$ is at least 2 and at least the
number of congruences in the given system).  The group~$F$ acts on this
space by left multiplication on the first coordinate: $g((h,n)) = (gh,n)$.
Let $[r]$ denote the set $\{1,2,\dots,r\}$.

In order to prove the second part of the theorem, we will randomly construct
a partition of~$F\times\nobreak\N$ into sets $A_1,A_2,\dots,A_r$ which
satisfies the given system but, with probability~$1$, satisfies
no congruence other than those deducible from the system by the
equivalence relation rules and the complementation rule.

The assignment of each pair~$(g,n)$ to one of the sets $A_1,\dots,A_r$
is done recursively on the reduced form of the group element~$g$.
For the identity element~$e$, assign $(e,n)$ to one of the sets~$A_k$ at
random with equal probability for each~$k$, and independently for all $n
\in \N$.  If $g\ne e$, then $g$~has a unique expression
as~$\rho\circ g'$ where $g'$~has a shorter reduced form than $g$~does,
and $\rho=f_i$ or $\rho=f_i^{-1}$ for some~$i \le m$.  Suppose that the
$i$\snug'th congruence in the given system is $\bigcup_{k\in L_i}A_k
\cong \bigcup_{k\in R_i}A_k$, where $L_i$ and~$R_i$ are nonempty proper
subsets of~$[r]$.  (If there is no $i$\snug'th congruence,
then we can just add a trivial and deducible $i$\snug'th congruence $A_1
\cong A_1$ to the system, so let $L_i = R_i = \{1\}$.)
If $\rho=f_i$, and we have already assigned $(g',n)$ to one of the
sets~$A_{k'}$, then put $(g,n)$ in~$A_k$, where: if $k' \in L_i$,
then $k$ is chosen randomly from~$R_i$; if $k' \notin L_i$, then
$k$ is chosen randomly from~$[r]\setminus R_i$.
If $\rho = f_i^{-1}$, then do the same thing, but with $L_i$
and~$R_i$ interchanged.  All random choices are to be made
uniformly from the options available and independently of each other.

It is easy to see that the sets constructed this way satisfy the
given congruences, with $f_i$~witnessing congruence number~$i$.
It remains to show that (with probability~$1$) no congruences not
deducible from this system are satisfied.

The sets~$A_k$ will (almost certainly) be nonempty; with probability~$1$,
each of the sets~$A_k$ will contain infinitely many points~$(e,n)$.
So the only congruences witnessed by
the identity element are those given by the reflexive law.

Define a nonempty set~$P_j(g)\subseteq\{1,2,\dots,r\}$ for each
$j\in[r]$ and $g \in F$ as follows.  If $g=e$,
then $P_j(g) = \{j\}$.  If $g = f_i \circ g'$ for some shorter~$g'$, then
$P_j(g)$ is $R_i$ if $P_j(g')\subseteq L_i$, $[r]\setminus R_i$
if $P_j(g') \cap L_i = \nullset$, and $[r]$ otherwise.
If $g = f_i^{-1} \circ g'$, do the same with $L_i$ and~$R_i$ interchanged.

The set $P_j(g)$ gives the possible values of~$k$ for which we can
have $(g,n) \in A_k$, given that $(e,n) \in A_j$.  We easily verify
by induction on~$g$ that, if $(e,n)\in A_j$, then $(g,n)$~must
be in $A_k$ for some~$P_j(g)$.  Furthermore, if $k \in P_j(g)$,
then the conditional probability that $(g,n) \in A_k$, given that
$(e,n) \in A_j$, is nonzero.  Since there are infinitely many~$n$\snug's
treated independently, with probability~$1$ there will be at least
one~$n$ such that $(e,n) \in A_j$ and $(g,n) \in A_k$.

It is straightforward to prove the following by induction on the length
of the reduced form of~$g$.  For each non-identity $g \in F$, there are
nonempty proper subsets $L(g)$ and~$R(g)$ of $[r]$ such that: if $j \in
L(g)$, then $P_j(g)$ is either $R(g)$ or $[r]$, and is the same for all
such~$j$; if $j \notin L(g)$, then $P_j(g)$ is either $[r]\setminus R(g)$
or $[r]$, and is the same for all such~$j$.  Furthermore, if $P_j(g)
\ne [r]$ for all~$j$, then the congruence $\bigcup_{k\in L(g)}A_k \cong
\bigcup_{k\in R(g)}A_k$ is deducible from the given system.

If $P_j(g) = [r]$ for some~$j$, then $g$~almost certainly cannot witness
any nontrivial congruence on the sets~$A_1,...,A_r$, because $g$~will
send points in~$A_j$ to all of the sets~$A_k$.  If $P_j(g) \ne [r]$
for all~$r$, then with probability~$1$ the only nontrivial congruences
witnessed by~$g$ are $\bigcup_{k\in L(g)}A_k \cong
\bigcup_{k\in R(g)}A_k$ and its complementary form, and both of these
are deducible from the given system.  So we have shown that (with
probability~$1$) no congruence holds between the sets~$A_k$ except
those deducible from the given system.  This completes the second part
ot the theorem.

For the first part of the theorem, we can use the same construction,
except that we will produce sets $A_1,A_2,\dots,A_{r+1}$ (so the sets
$A_1,\dots,A_r$ will no longer be a partition of the whole space).
We proceed exactly as above, except that $[r]$ is replaced by $[r+1]
= \{1,2,\dots,r+1\}$ throughout.  (This is why we were careful to
use $[r]\setminus R$ instead of writing $R^c$ in the above argument.)
We may assume $r+1 \notin L(g)$ for all~$g$ (otherwise, just replace
$L(g)$ and~$R(g)$ with their complements in $[r+1]$).  Since the given
congruences only involve sets $A_1,\dots,A_r$, it is easy to see that
$r+1 \in P_{r+1}(g)$ for all~$g$.  We now find that, if $g\in F$~is such
that $g \ne e$ and $P_j(g) \ne [r+1]$ for all~$j$, then $r+1 \notin R(g)$
and the congruence $\bigcup_{k\in L(g)}A_k \cong \bigcup_{k\in R(g)}A_k$
(which is a congruence among the sets $A_1,\dots,A_r$) is deducible from
the given system using the equivalence relation rules alone.  This is
(with probability~$1$) the only case in which a non-identity~$g$ can
witness a nontrivial congruence among the sets $A_1,\dots,A_r$, so no
such congruence holds except those deducible from the given system by
the equivalence relation rules.
\QED\enddemo

Actually, the argument for the first part of Theorem~7.1 works even if
improper congruences involving $A_1 \cup A_2 \cup \dots \cup A_r$ are
allowed in the system; it is only the congruences involving~$\nullset$
that must be excluded in this case.

We also considered subcongruences $\bigcup_{k\in L}A_k \scong
\bigcup_{k\in R}A_k$ in section~2, and gave the following deduction rules:
the inclusion rule (if $L \subseteq R$, then $\bigcup_{k\in L}A_k \scong
\bigcup_{k\in R}A_k$); transitivity; and, from $B \cong C$, one can deduce
$B \scong C$ and $C \scong B$.  Again there is a complementation rule
($\bigcup_{k\in L}A_k \scong \bigcup_{k\in R}A_k$ implies $\bigcup_{k\in
R^c}A_k \scong \bigcup_{k\in L^c}A_k$) in the case where the sets~$A_k$
form a partition of the space.  And again it is natural to ask whether
this set of rules is complete.

Just as for congruences, we run into difficulties if we allow improper
subcongruences (or improper congruences) in our assumptions.  For
instance, if the subcongruence $A_1 \scong \nullset$ is true, then
the subcongruence $A_1 \cup A_2 \scong A_2$ (and even the congruence
$A_1 \cup A_2 \cong A_2$) must also be true, but we cannot deduce this
from the given rules.  There are similar difficulties if we assume
an improper subcongruence of the form $A_1 \cup \dots \cup A_r \scong B$
in the partition case.

However, again as before, if we restrict ourselves to proper congruences
and subcongruences, then the answer is yes:

\proclaim{Theorem 7.2}
If one has a system of proper congruences and proper subcongruences, and an
additional subcongruence which is not deducible from the system by the
subcongruence rules (excluding complementation), then one can find a
suitable space and nonempty open subsets of that space which satisfy
the system of congruences and subcongruences but not the additional
subcongruence.  If the additional subcongruence is not deducible from
the system using the subcongruence rules including the complementation
rule, then the open subsets of the suitable space can be taken to form
a partition of the space.
\endproclaim

\demo{Proof}
The proof is very similar to that of Theorem~7.1.  Again use the suitable
space $F\times\N$, where $F$~is free on $m$~generators and
$m$~is at least the number of given congruences and subcongruences.
(In fact, we may assume $m$~is exactly this number, since we can add
trivial congruences $A_1 \cong A_1$
or subcongruences $A_1 \scong A_1$ to the given system.)

For the second part of the theorem, we randomly generate a partition
of $F\times\N$ into pieces $A_1,A_2,\dots,A_r$ as before.  The difference is
that we need to handle the case $g = \rho \circ g'$ where $\rho=f_i^{\pm 1}$
and the $i$\snug'th given congruence or subcongruence is a subcongruence.
Suppose this subcongruence is $\bigcup_{k\in L_i}A_k \scong
\bigcup_{k\in R_i}A_k$.  Then, if $\rho = f_i$ and
$(g',n)$ has been assigned to~$k'$
where $k' \in L_i$, we choose~$k$ randomly from $R_i$ and assign
$(g,n)$ to~$A_k$; if $k' \notin L_i$, we choose~$k$ randomly from $[r]$.
If $\rho = f_i^{-1}$, then if $k' \in R_i$, we choose $k$ randomly
from~$[r]$, while if $k' \notin R_i$, we choose $k$ randomly from
$[r]\setminus L_i$.
Again the resulting sets~$A_k$ must satisfy the given congruences and
subcongruences, with the $i$\snug'th of them being witnessed by~$f_i$.

Define $P_j(g)$ as before, but with new clauses: If the $i$\snug'th
member of the given system is the subcongruence $\bigcup_{k\in L_i}A_k \scong
\bigcup_{k\in R_i}A_k$, then, if $g = f_i \circ g'$, let $P_j(g)$ be
$R_i$ if $P_j(g') \subseteq L_i$, $[r]$~otherwise.  If $g = f_i^{-1} \circ
g'$, let $P_j(g)$ be $[r]\setminus L_i$ if $P_j(g') \cap R_i = \nullset$,
$[r]$~otherwise.  Again we get that (with probability~$1$) there exists
$n\in\N$ such that $(e,n) \in A_j$ and $(g,n) \in A_k$ if and only if
$k \in P_j(g)$.

Again, for each non-identity~$g$, there are nonempty proper subsets $L(g)$
and~$R(g)$ of~$[r]$ such that: if $j \in L(g)$, then $P_j(g)$ is either
$R(g)$ or $[r]$, and is the same for all such~$j$; if $j \notin L(g)$,
then $P_j(g)$ is either $[r]\setminus R(g)$ or $[r]$, and is the same
for all such~$j$.  Furthermore, if $P_j(g)=R(g)$ for $j \in L(g)$, then
the subcongruence $\bigcup_{k\in L(g)}A_k \scong \bigcup_{k\in R(g)}A_k$
is deducible from the given system; if $P_j(g)=[r]\setminus R(g)$ for
$j \notin L(g)$, then the reverse subcongruence $\bigcup_{k\in R(g)}A_k
\scong \bigcup_{k\in L(g)}A_k$ is deducible from the given system.

Now, the only cases in which a non-identity group element~$g$ witnesses
a nontrivial subcongruence $\bigcup_{k\in L}A_k \scong \bigcup_{k\in
R}A_k$ (here `nontrivial' means $L \ne \nullset$ and $R \ne [r]$) are
when $P_j(g) = R(g)$ for $j \in L(g)$, $L \subseteq L(g)$, and $R(g)
\subseteq R$, or when $P_j(g) = [r]\setminus R(g)$ for $j \notin L(g)$,
$L \subseteq [r]\setminus L(G)$, and $[r]\setminus R(g) \subseteq R$.
In either of these cases, the subcongruence $\bigcup_{k\in L}A_k
\scong \bigcup_{k\in R}A_k$ is deducible from the given system by the
subcongruence rules.  Therefore, the subcongruence rules are complete
for the second part of the theorem.

For the first part of the theorem, we again produce sets
$A_1,\dots,A_{r+1}$ instead of $A_1,\dots,A_r$ and replace $[r]$ with
$[r+1]$ throughout.  We may assume that $r+1 \notin L(g)$ for all~$g$.
Now, if a subcongruence $\bigcup_{k\in L}A_k \scong \bigcup_{k\in
R}A_k$ among the first $r$~sets is witnessed by the non-identity group
element~$g$, then we must have $L \subseteq L(g)$ and $R(g) \subseteq
R$, where $P_j(g) = R(g)$ for $j \in L(g)$ and $r+1 \notin R(g)$; in
this case, induction on~$g$ shows that $\bigcup_{k\in L(g)}A_k \scong
\bigcup_{k\in R(g)}A_k$ is deducible from the given system by the
subcongruence rules without using complementation, so the same holds
for $\bigcup_{k\in L}A_k \scong \bigcup_{k\in R}A_k$.  So again the
subcongruence rules are complete.
\QED\enddemo

One can also note in the proof of Theorem~7.2 that (with probability~$1$)
the only case in which the sets $A_1,\dots,A_r$ satisfy a congruence
is when this congruence is deducible from the congruences in the given
system by the congruence rules.  In other words, there are no useful
rules for using subcongruences (alone or in conjunction with congruences)
to deduce congruences; any congruence which follows from given proper
congruences and subcongruences must follow from the given congruences
alone.

\head 8. Open questions \endhead

A number of the theorems in this paper give specific examples
rather than general results.  Regarding general results, many of the
main questions remain open.  A few open questions have been mentioned
already (the satisfiability of $\UNC6$, and the converses of some
implications in section~6); here we list some more.

The main question remaining open is: can one give an explicit
characterization (in whatever form) of the satisfiable congruences,
in any of the cases listed in section~6?  Such characterizations have
been given for solutions to systems of congruences using arbitrary sets
(Robinson~\cite{\Robinson}, Adams~\cite{\Adams}) or using sets with
the property of Baire (Dougherty~\cite{\Dougherty}), but none has yet been
found for the open-sets cases.

In particular, is it even recursively decidable whether a given system
of congruences is (nontrivially) satisfiable, in any of these cases?
The possibility that this is undecidable is not entirely implausible;
since arbitrary computations can be coded in cellular automata
and related systems, it is conceivable that they could be encoded
in systems of congruences, so that, say, the system is satisfiable
by finite subsets of a free group (not all empty) if and only if
the computation terminates.

However, in this particular case, there is a partial decidability result.
If the group elements that are to witness the congruences are fixed in
advance, then the satisfiability question is decidable:

\proclaim{Proposition~8.1} There is an algorithm which, when given
a natural number~$m$, a system of~$k$ congruences, and
elements $g_1,\dots,g_k$ of the free group~$F_m$ on $m$ generators,
will decide whether there are finite subsets (not all empty) of~$F_m$
which satisfy the given congruences, where $g_i$~is the witness
for the $i$\snug'th congruence, $i=1,\dots,k$. \endproclaim

\demo{Proof} Let $L$ be the maximum of the lengths of the
group elements $g_1,\dots,g_k$ expressed as words in the generators
of~$F_m$, and let $N = 1 + 2m + (2m)^2 + \dots + (2m)^L$.
We will show that, if there exist finite subsets $A_1,\dots,A_r$
of~$F_m$ (not all empty) satisfying the congruences, with $g_i$ witnessing
the $i$\snug'th congruence for all $i \le k$, then there exist
such subsets consisting entirely of words of length less than
$(r+1)^N$.  This reduces the existence problem to a finite search,
so the problem is decidable.

Assume that there exist sets $A_1,\dots,A_r$ satisfying the congruences
as above.  We may assume that the identity element is in one of the
sets~$A_j$, because, given any element~$h$ of one of the sets, we can
multiply all elements of all of the sets by~$h^{-1}$ on the right
to get a new sequence of sets satisfying the congruences as before.
Now, among such $r$\snug-sequences of sets satisfying the congruences
(as witnessed by~$g_i$) and containing the identity element,
take $A_1,\dots,A_r$ to be one such that the sum of the lengths
of the words in $A_1 \cup \dots \cup A_r$ is as small as possible.
We will see that these sets cannot contain any word of length
as large as $(r+1)^N$.

Suppose $w$~is a reduced word of length at least $(r+1)^N$ which is in
one of the sets~$A_j$.  For each of the final segments~$v$ of~$w$,
let $p_v$ be the function whose domain is the set of words of
length at most~$L$ (note that there are $N$ of these), such that
$p_v(z)=j$ if $z \circ v \in A_j$, and
$p_v(z)=0$ if $z \circ v \notin A_1 \cup \dots \cup A_r$.
The number of possible functions~$p_v$ is $(r+1)^N$; since the
number of final segments~$v$ of~$w$
(counting the identity element and the word~$w$ itself) is greater
than $(r+1)^N$, there must exist final segments $v$ and~$v'$ with
$v$~shorter than~$v'$ (so $w = x \circ v'$ and $v' = y \circ v$
for sone words $x$ and~$y$, with no cancellation) such that
$p_v = p_{v'}$.

Now construct new subsets $A'_1,\dots,A'_r$ of $F_m$ as follows.
If the reduced word~$h$ does not end in~$v$, then put
$h \in A'_j$ iff $h \in A_j$ for all~$j$.  If $h$~does
end in~$v$, say $h = h' \circ v$, then put $h \in A'_j$ iff
$h' \circ v' \in A_j$.  This `cut-and-splice' operation
does not alter the relevant properties of the sets except near
the cut points $v$ and~$v'$.  Using the fact that
$p_v = p_{v'}$ (i.e., the sets $A_1,\dots,A_r$ ``look the same near~$v'$
as they do near~$v$''), it is not hard to show that the
sets $A'_1,\dots,A'_r$ satisfy the congruences as witnessed by
the group elements~$g_i$, since the sets $A_1,\dots,A_r$ do.
But the sum of the lengths of the words in~$A'_1 \cup \dots \cup A'_r$
is less than the sum of the lengths of the words in~$A_1 \cup \dots \cup A_r$.
This contradicts the minimality assumed earlier.  Therefore, the
word~$w$ cannot exist, and we are done.
\QED\enddemo

One possible form of a characterization of the satisfiable systems
in some context would be a list of systems which is universal in the
sense that any system is satisfiable if and only if it is reducible to
a system on the list.  We saw such a characterization of the numerically
consistent systems in section~2.  (Of course, the numerically consistent
systems can be characterized directly from the definition; it is a simple
linear programming problem to determine whether a system is numerically
consistent.)  Can such a universal list be given in any of the other
cases from section~6?  Note that such a list would not immediately imply
decidability of the satisfiability problem, even if the list were decidable.

Even if one is more interested in general results applying to
arbitrary suitable spaces or the like, the specific case
of the sphere with free rotations is useful as a source of limitative
results (showing that certain systems {\it cannot} be satisfied
nontrivially in general).  It would be helpful to have other
specific suitable spaces where systems of congruences can be shown
to be unsatisfiable.  The discrete free groups are of no use for this
purpose; any system of congruences has solutions there.  One
possible such space which deserves further study is the Cantor
space acted on freely by a free group of Lipschitz homeomorphisms.

In all of the cases we have examined involving free rotations of the
sphere, the arguments worked for arbitrary free rotations; it did not
matter which ones were used.  Is this always the case, or could it be
that there is a system of congruences satisfiable on the sphere under
one free group of rotations but not under a different free group?

The open sets produced by some of the constructions in this paper
are highly pathological (having infinitely many connected components,
boundaries of positive measure, etc.); one can consider what happens
if one is restricted to `nicer' open sets.  In particular, for what
systems of congruences can we find solutions using dissections of
the sphere?  Of course, one must define the term `dissection'; one
way to do this would be as the complement of a finite graph embedded
in the sphere.  (Is this significantly more restrictive than just
requiring the open sets in question to have finitely many
connected components?  What if the open sets actually have to
be connected?)  If we ask whether one finite union of pieces in
a dissection is congruent to another such finite union, should we
`erase' (i.e., add in) the boundary lines between adjacent pieces
in the same union?  This apparently gives a whole family of new
satisfiability questions, and one can ask whether the satisfiable
congruences can be characterized, or what the implications are
between these cases and those listed in section~6.

Finally, we should recall that a number of questions about
solutions to systems of congruences using Borel sets, or using
Lebesgue measurable sets, have been open for a long time.  For
instance, there is Question~4.15 from Wagon~\cite{\Wagon} (due to
Mycielski), which asks whether the system $A_1 \cong
A_2 \cong A_3$ has a solution using measurable subsets of~$S^2$.  So a
characterization of the solvable systems of congruences in these cases
appears to be a long way off.

\enddocument

%% file: plot.tex
%
%
%
%
%
\count255=\catcode`\!
\catcode`\!=11
\ifx\plot!loaded\relax
   \catcode`\!=\count255 \else\let\plot!loaded=\relax\fi
\chardef\plot!savecc=\count255
%
%
\def\plot!zero{0}
\def\plot!one{1}
\newdimen\plotunitx
\newdimen\plot!unitxu
\def\plot!figscalex{1}
\plotunitx=1truebp
\plot!unitxu=1truebp
\newdimen\plotunity
\newdimen\plot!unityu
\def\plot!figscaley{1}
\plotunity=1truebp
\plot!unityu=1truebp
\def\plot!sepscfac{1}
\newtoks\plot!symbol
\newbox\plot!figurebox
\newdimen\plotcurrx
\newdimen\plotcurry
\newdimen\plotlinewidth
\plotlinewidth=.6bp
\def\plot!savemem{0}
\let\plot!global=\relax
\newif\ifplot
\plottrue
\newif\ifplotPS
\newif\ifplotseparate
\newif\ifplot!infig
\newif\ifplot!trans
\newif\ifplot!loctrans
\newif\ifplot!sepscaled
\newif\ifplot!nowraw
\newif\ifplot!globdef
\newif\ifplot!rawbounds
\newif\ifplot!spseq
\newif\ifplot!stacklock
\newif\ifplot!lhcs
\plot!spseqtrue
%
%
\mathchardef\plot!ci="220E
\mathchardef\plot!bu="220F
\mathchardef\plot!cd="2201
\def\plot!sp{ }
\def\plot!em{}
\let\plot!bg={
\let\plot!eg=}
\def\plot!ni{\prevdepth=-1000pt }
\def\plot!sm#1{{\setbox0=\hbox{#1}\ht0=0pt \dp0=0pt \box0 }}
\def\plot!lo#1\plot!re{\def\plot!bo{#1}\plot!it}
\def\plot!it{\plot!bo \let\plot!ne=\plot!it \else\let\plot!ne=\relax\fi
   \plot!ne}
\let\plot!re=\fi
%
%
\long\def\plot!firstlet#1#2\plot!endofarg{#1}
\let\plot!PStrueold=\plotPStrue
\def\plotPStype#1{
   \plot!PStrueold
   \ifcase#1
      \plotPSfalse
   \or
      \def\plot!local{" }
      \def\plot!global{! }
      \def\plot!rawbegin##1{ps:SDict begin ##1\plot!sp end}
      \def\plot!raw##1{ps:SDict begin ##1\plot!sp end}
      \def\plot!rawend##1{ps:SDict begin ##1\plot!sp end}
      \def\plot!trans{currentpoint /p!s1 2 index def /p!s2 1 index def
         translate exec 0 0 moveto p!s1 neg p!s2 neg translate }
      \def\plot!setorig{}
      \def\plot!rawsetcurr{}
      \def\plot!rawsetorig{currentpoint translate }
      \def\plot!PSfile##1{\includegraphics{##1}}
      \plot!globdeftrue
      \plot!stacklocktrue
      \plot!lhcstrue
   \or
      \def\plot!local{ps:: }
      \def\plot!global{ps::[global] }
      \def\plot!rawbegin##1{ps::[inline,begin] ##1}
      \def\plot!raw##1{ps::[inline] ##1}
      \def\plot!rawend##1{ps::[inline,end] ##1}
      \def\plot!trans{Xpos Ypos translate
         exec Xpos neg Ypos neg translate }
      \def\plot!setorig{Xpos Ypos translate }
      \def\plot!rawsetcurr{Xpos Ypos moveto }
      \def\plot!rawsetorig{Xpos Ypos translate }
      \def\plot!PSfile##1{\plot!rawstart{gsave Xpos Ypos translate}%
         \special{ps: plotfile ##1 inline}\plot!rawfinish{grestore}}
      \plot!rawboundstrue
   \or
      \def\plot!local{empty.ps }
      \def\plot!rawbegin##1{empty.ps ##1}
      \def\plot!raw##1{empty.ps ##1}
      \def\plot!rawend##1{empty.ps ##1}
      \def\plot!trans{pop }
      \def\plot!setorig{}
      \def\plot!rawsetcurr{0 0 moveto }
      \def\plot!rawsetorig{}
      \def\plot!PSfile##1{\special{##1}}
      \plot!spseqfalse
   \else
      \immediate\write16{plotPStype: Illegal type specified -- type set to 0}
      \plotPSfalse
   \fi
   \ifplotPS\ifplot!globdef\else\ifx\plot!global\relax\else
      \ifcase\plot!savemem\or
         \plot!globdeftrue\plot!gwarning
      \or
         \immediate\write16{ }
         \immediate\write16{Do you want me to try to conserve TeX memory?}
         \immediate\write16{(If so, the resulting pages must be printed in %
            order from the beginning.)}
         \message{? }
         \read-16 to\plot!inputstr
         \edef\plot!inputstr{\plot!firstlet\plot!inputstr xx\plot!endofarg}
         \if y\plot!inputstr\plot!globdeftrue\fi
         \if Y\plot!inputstr\plot!globdeftrue\fi
      \fi
   \fi\fi\else\global\let\plot!gwarning=\relax\fi
   \ifplot!globdef\global\let\plot!gwarning=\relax\fi
   \plot!PSinit   
}
\def\plotPSask{
   \immediate\write16{ }
   \immediate\write16{Which DVI-to-PS converter are you using?}
   \immediate\write16{  1. Rokicki dvips (Radical Eye)}
   \immediate\write16{  2. ArborText dvips (DVILASER/PS)}
   \immediate\write16{  3. OzTeX [partial functionality]}
   \immediate\write16{  0. None of the above %
      [PostScript insertions replaced or suppressed]}
   \message{Enter a number: }
   \read-16 to\plot!inputstr
   \edef\plot!inputstr{\plot!firstlet\plot!inputstr xx\plot!endofarg}
   \ifcat0\plot!inputstr\else\def\plot!inputstr{Z}\fi
   \ifnum\expandafter`\plot!inputstr<`0\def\plot!inputstr{99}\fi
   \ifnum\expandafter`\plot!inputstr>`9\def\plot!inputstr{99}\fi
   \plotPStype{\plot!inputstr}
}
\def\plotPStrue{\plotPSask}

\let\plot!septrueold=\plotseparatetrue
\def\plotseparatetrue{\plot!septrueold\plottrue}
\let\plot!sepfalseold=\plotseparatefalse
\def\plotseparatefalse{\plot!sepfalseold\plotnosepscale}
\let\plot!falseold=\plotfalse
\def\plotfalse{\plot!falseold\plotseparatefalse}
\def\plotask{
   \immediate\write16{ }
   \immediate\write16{How do you want figures to appear?}
   \immediate\write16{  1. Normally within text}
   \immediate\write16{  2. On separate pages (blank spaces left in text)}
   \immediate\write16{  0. Not at all (blank spaces left in text)}
   \message{Enter a number: }
   \read-16 to\plot!inputstr
   \if 0\plot!firstlet\plot!inputstr xx\plot!endofarg
      \plotfalse
   \else\if 1\plot!firstlet\plot!inputstr xx\plot!endofarg
      \plottrue \plotseparatefalse
   \else\if 2\plot!firstlet\plot!inputstr xx\plot!endofarg
      \plotseparatetrue
   \else
      \immediate\write16%
         {plotask: Illegal option specified -- figures suppressed}
      \plotfalse
   \fi\fi\fi
}
\gdef\plot!gwarning{
   \immediate\write16{ }
   \immediate\write16{plot WARNING: This DVI-to-PS converter does not %
      fully support global}
   \immediate\write16{ \plot!sp\plot!sp definitions -- printing pages out %
      of order may cause problems.}
   \global\let\plot!gwarning=\relax}
%
%
\def\plot!zerocp{\plotcurrx=0pt\plotcurry=0pt\relax}
{\toksdef\abc=7 \catcode`p=12 \catcode`t=12 \global\abc={pt}}
\expandafter\def\expandafter\plot!unpt\expandafter#\expandafter1\the\toks7{#1 }
\def\plotthedim#1{\expandafter\plot!unpt\the#1 }

\def\plotpdim#1{\plotthedim{#1}\ifnum1000=\mag 1.00375 \else\ifx\plot!unmag
   \plot!em 1003.75 \the\mag\plot!sp div \else 1.00375 \fi\fi div }
\def\plotfdims#1#2{\plotthedim{#1}\plotthedim{\plotunitx}div
   \plotthedim{#2}\plotthedim{\plotunity}div }
\def\plot!recip#1{{\dimen6=1in \dimen8=#1in \dimen0=0sp \global\toks7={}
   \plot!lo
      \count255=\dimen6
      \divide\count255 by\dimen8
      \global\toks7=\expandafter\expandafter\expandafter{\expandafter\the
         \expandafter\toks\expandafter7\expandafter\plot!sp\the\count255}
      \ifdim\dimen0=0sp\global\toks7=\expandafter{\the\toks7.}\fi
      \dimen4=\dimen8
      \multiply\dimen4 by\count255
      \advance\dimen6 by -\dimen4
   \ifdim\dimen0<6sp
      \advance\dimen0 by 1sp
      \multiply\dimen6 by 10
   \plot!re}}
\def\plot!applyssf{\plotunitx=\plot!sepscfac\plot!unitxu
   \plotunity=\plot!sepscfac\plot!unityu}
\def\plot!PSssf{\ifplot!sepscaled\plot!sepscfac\plot!sp dup scale \fi}
\def\plot!errmsg{\errmessage
   {plot error: command allowed only within \string\plotfigbegin
   ...\string\plotfigend}}
\def\plot!rawstart#1{\plot!bg\ifplot!nowraw\special{\plot!raw{#1}}\else
   \special{\plot!rawbegin{#1}}\plot!nowrawtrue\plot!rawboundsfalse\fi}
\def\plot!rawfinish#1{\plot!eg\ifplot!nowraw\special{\plot!raw{#1}}\else
   \special{\plot!rawend{#1}}\fi}
\def\plot!rawstartnull{\ifplot!rawbounds\plot!rawstart{}\else
   \plot!bg\plot!nowrawtrue\fi}
\def\plot!rawfinishnull{\plot!eg\ifplot!rawbounds\ifplot!nowraw
   \else\special{\plot!rawend{}}\fi\fi}
%
%
\def\plotfigscales#1 #2 {\def\plot!figscalex{#1} \plot!unitxu=#1truebp
   \def\plot!figscaley{#2} \plot!unityu=#2truebp \plot!applyssf}
\def\plotfigscale#1 {\plotfigscales{#1} {#1} }
\def\plotrotate#1 {\def\plot!symangle{#1}\plot!loctranstrue\plot!transtrue}
\def\plotscale#1 {\def\plot!symscale{#1}\plot!loctranstrue\plot!transtrue}
\def\plotcancelmag{\def\plot!unmag{\ifnum1000=\mag \else1000 \the\mag\plot!sp
   div dup scale \fi}\plot!loctranstrue\plot!transtrue}
\def\plotresumemag{\def\plot!unmag{}}
\long\def\plotgentrans#1\plotPSend{\plot!loctranstrue\plot!transtrue
   \long\def\plot!transcode{\ifplot!lhcs 1 -1 scale \fi
   #1\ifplot!lhcs\plot!sp 1 -1 scale\fi}}
\def\plotnotrans{\def\plot!transcode{\ifx\plot!symangle\plot!zero\else
   \plot!symangle\plot!sp\ifplot!lhcs neg \fi rotate \fi
   \ifx\plot!symscale\plot!one\else\plot!symscale\plot!sp dup scale \fi}
   \plotresumemag\plot!loctransfalse\ifplot!sepscaled
   \else\plot!transfalse\fi\def\plot!symscale{1}\def\plot!symangle{0}}
\def\plotsepscale#1 {\ifplotseparate\def\plot!sepscfac{#1}
   \ifx\plot!sepscfac\plot!one\plotnosepscale\else
   \plot!sepscaledtrue\plot!transtrue\plot!applyssf \fi\fi}
\def\plotnosepscale{\plot!sepscaledfalse\ifplot!loctrans\else\plot!transfalse
   \fi\def\plot!sepscfac{1}\plot!applyssf}
%
%
\def\plotfigbegin{\setbox\plot!figurebox=\vbox\plot!bg
   \plot!infigtrue\plotnotrans\plot!zerocp}
\def\plotfigend{\xdef\plot!gtemp{\plot!sepscfac}\plot!eg
   {\dimen2=\ht\plot!figurebox
   \dimen4=\wd\plot!figurebox
   \dimen6=\dp\plot!figurebox
   \ifplotseparate
      \shipout\vbox to9truein{\vss\hbox to6.5truein{\hss\box\plot!figurebox
         \hss}\vss}
      \plot!recip{\plot!gtemp}
      \setbox\plot!figurebox=\hbox to \the\toks7\dimen4{\hfil
         \vrule width0pt height\the\toks7\dimen2 depth\the\toks7\dimen6}
   \else\ifplot\else
      \setbox\plot!figurebox=\hbox to \dimen4{\hfil
         \vrule width0pt height\dimen2 depth\dimen6}
   \fi\fi
   \box\plot!figurebox}}
\def\plot!plot{\ifplot
   \vbox to0pt{\kern-\dimen8
   \hbox{%
      \kern\dimen6
      \let\plot!temp=\plot!em
      \ifplotPS\ifplot!trans\ifplot!spseq
         \edef\plot!temp{\plot!unmag\plot!PSssf\plot!transcode}\fi\fi\fi
      \ifx\plot!temp\plot!em\else\plot!rawstart{{\plot!temp} \plot!pla}\fi
      \plot!sm{\hbox to0pt{\hss\the\plot!symbol\hss}}%
      \ifx\plot!temp\plot!em\else\plot!rawfinish{grestore}\fi
   }
   \kern\dimen8}
   \plot!ni\fi}
\def\plot#1 #2 {\ifplot!infig
   {\dimen6=#1\plotunitx\dimen8=#2\plotunity\plot!plot}
   \else\plot!errmsg\fi}
\def\plothere{\ifplot!infig
   {\dimen6=\plotcurrx\dimen8=\plotcurry\plot!plot}
   \else\plot!errmsg\fi}
\def\plotmove#1 #2 {\ifplot!infig
   \plotcurrx=#1\plotunitx\plotcurry=#2\plotunity
   \else\plot!errmsg\fi}
\def\plotrmove#1 #2 {\ifplot!infig
   \advance\plotcurrx by #1\plotunitx \advance\plotcurry by #2\plotunity
   \else\plot!errmsg\fi}
%
%
\def\plottext#1{\plot!symbol={\lower\fontdimen22\textfont2\hbox{#1}}}

\def\plotsfmla#1{\plot!symbol={\lower\fontdimen22\scriptfont2\hbox
   {$\scriptstyle #1$}}}
\def\plotssfmla#1{\plot!symbol={\lower\fontdimen22\scriptscriptfont2\hbox
   {$\scriptscriptstyle #1$}}}
\def\plotcentered#1{\plot!symbol={\vbox to0pt{\vss\hbox{#1}\vss}}}

%
%
\long\def\plotPSbegin#1\plotPSend{
   \ifplot!infig
   \ifplot\ifplotPS
      \let\plot!temp=\plot!em
      \ifx\plot!one\plot!figscalex\else\let\plot!temp=\plot!one\fi
      \ifx\plot!one\plot!figscaley\else\let\plot!temp=\plot!one\fi
      \ifx\plot!temp\plot!em\else\def\plot!temp{\plot!figscalex\plot!sp
         \ifdim\plot!unitxu=\plot!unityu dup dup scale \else
         \plot!figscaley\plot!sp 2 copy scale dup mul exch
         dup mul 2 div sqrt \fi div }\fi
      \special{\plot!local
         gsave \plot!setorig
         \ifnum1000=\mag\else1000 \the\mag\plot!sp div dup scale \fi
         \ifplot!sepscaled\plot!sepscfac\plot!sp dup scale \fi
         \plotpdim{\plotlinewidth}\plot!temp setlinewidth
         \plotfdims{\plotcurrx}{\plotcurry}moveto #1\plot!sp grestore}
   \fi\fi
   \else\plot!errmsg\fi}
\long\def\plotPSglobal#1\plotPSend{%
   \ifplotPS
      \ifx\plot!global\relax
         \errmessage{plotPSglobal: This PS type does not support global %
            insertions}%
      \else
         \special{\plot!global #1}%
         \plot!gwarning
      \fi
   \fi }
\long\def\plotPSmath#1\plotPSend{{%
   \ifplotPS
      \mathchoice
         {\special{\plot!raw{1.00001}}}%
         {\special{\plot!raw{1.0}}}%
         {\special{\plot!raw{0.7}}}%
         {\special{\plot!raw{0.5}}}
      \special{\plot!local gsave \plot!setorig #1\plot!sp
         \ifplot!stacklock 0 \fi grestore}
      \ifplot!stacklock\special{\plot!raw{count 0 ne {pop} if}}\fi
   \fi}}
\long\def\plotPSgdef#1#2\plotPSend{{%
   \ifplotPS
      \escapechar=-1 %
      \ifplot!globdef
         \special{\plot!global /p!0\string#1 {#2} def}%
         \xdef#1{p!0\string#1 }%
      \else
         \xdef#1{#2\plot!sp}%
      \fi
   \else\global\let#1=\relax\fi}}
\def\plot!befsetup{\setbox0=\hbox{%
   \hskip 1in\raise1in\hbox{\special{\plot!raw{\plot!bfa}}}}%
   \ht0=0pt \wd0=0pt \box0}

\long\def\plotPSbefore#1#2\plotPSend{{%
   \setbox8=\hbox{#1}%
   \ifplotPS
      \hbox{\plot!befsetup\special{\plot!raw{\plot!bfb #2\plot!sp
      grestore}}\box8}%
   \else \box8\fi }}
\long\def\plotPSduring#1#2\plotPSend{{%
   \toks0={{#1}}%
   \ifplotPS\ifplot!spseq
      \hbox{\plot!rawstart{gsave #2}\the\toks0 \plot!rawfinish{grestore}}%
   \else \the\toks0\fi \else \the\toks0\fi }}
\long\def\plotPSbefaft#1#2\plotPSend#3\plotPSend{{%
   \toks0={{#1}}\setbox8=\hbox{\the\toks0}\def\plot!bef{#2}\def\plot!aft{#3}%
   \ifplotPS\ifplot!spseq
      \setbox8=\hbox{\vrule width\wd8 height\ht8 depth\dp8}%
      \hbox to\wd8{\plot!befsetup
         \plot!rawstart{\plot!bab
            \ifx\plot!bef\plot!em\else
            {\plot!bac \plot!bef\plot!sp \plot!bad} \plot!tra \fi}%
         \hbox to0pt{\the\toks0 \hss}%
         \plot!rawfinish{\ifx\plot!aft\plot!em\else {\plot!bac
            \plot!aft} \plot!tra \fi grestore}\hfil}%
   \else \box8\fi \else \box8\fi }}
\long\def\plotPSraw#1\plotPSend{%
   \ifplotPS\ifplot!spseq\special{\plot!raw{#1}}\fi\fi }
%
%
\newdimen\plotdotspc
\plotdotspc=2bp
\def\plotDRline{\ifplot!infig
   \ifplot\ifplotPS
      {\plotfigscale 1
      \edef\plot!litemp{\plotfdims{\dimen0}{\dimen2}}%
      \plotPSbegin\plot!litemp\plot!lia \plotPSend}
   \else\plot!TeXline\fi\fi
   \else\plot!errmsg\fi}
\def\plot!TeXline{{
   \plotlinewidth=\plot!sepscfac\plotlinewidth
   \plotdotspc=\plot!sepscfac\plotdotspc
   \ifdim\dimen0=\plotcurrx
      \dimen6=\dimen0
      \dimen4=\dimen2
      \advance\dimen4 by -\plotcurry
      \ifdim\dimen4<0sp\dimen4=-\dimen4\fi
      \advance\dimen4 by \plotlinewidth
      \plotcentered{\vrule width\plotlinewidth height\dimen4}
      \dimen8=\dimen2
      \advance\dimen8 by \plotcurry
      \divide\dimen8 by 2
      \plot!plot
   \else\ifdim\dimen2=\plotcurry
      \dimen8=\dimen2
      \dimen4=\dimen0
      \advance\dimen4 by -\plotcurrx
      \ifdim\dimen4<0sp\dimen4=-\dimen4\fi
      \advance\dimen4 by \plotlinewidth
      \plotcentered{\vrule width\dimen4 height\plotlinewidth}
      \dimen6=\dimen0
      \advance\dimen6 by \plotcurrx
      \divide\dimen6 by 2
      \plot!plot
   \else\plot!TeXdotline\fi\fi}}
\def\plot!TeXdotline{      
   \plotcentered{\vrule width\plotlinewidth height\plotlinewidth}
   \ifdim\dimen0<\plotcurrx
      \dimen6=\dimen0
      \dimen8=\dimen2
      \dimen0=\plotcurrx
      \dimen2=\plotcurry
   \else
      \dimen6=\plotcurrx
      \dimen8=\plotcurry
   \fi
   \advance\dimen2 by -\dimen8
   \dimen4=\dimen0
   \advance\dimen4 by -\dimen6
   \global\dimen1=\dimen2
   \ifdim\dimen1<0sp\global\dimen1=-\dimen1\fi
   \global\advance\dimen1 by \dimen4
   \global\divide\dimen1 by \plotdotspc
   \global\advance\dimen1 by 1sp
   \global\dimen7=\dimen4
   \global\divide\dimen7 by \dimen1
   \global\dimen3=\dimen7
   \global\multiply\dimen3 by -\dimen1
   \global\advance\dimen3 by \dimen4
   \global\dimen9=\dimen2
   \global\divide\dimen9 by \dimen1
   \global\dimen5=\dimen9
   \global\multiply\dimen5 by -\dimen1
   \global\advance\dimen5 by \dimen2
   \ifdim\dimen5<0sp
      \global\advance\dimen5 by \dimen1
      \global\advance\dimen9 by -1sp
   \fi
   \dimen2=0sp
   \dimen4=0sp
   \plot!lo
      \plot!plot
   \ifdim\dimen6<\dimen0
      \advance\dimen2 by \dimen3
      \advance\dimen4 by \dimen5
      \advance\dimen6 by \dimen7
      \advance\dimen8 by \dimen9
      \ifdim\dimen2<\dimen1\else
         \advance\dimen6 by 1sp
         \advance\dimen2 by -\dimen1
      \fi
      \ifdim\dimen4<\dimen1\else
         \advance\dimen8 by 1sp
         \advance\dimen4 by -\dimen1
      \fi
   \plot!re
}
\def\plotline#1 #2 {
   {\dimen0=#1\plotunitx \dimen2=#2\plotunity
   \plotDRline}
   \plotmove {#1} {#2}
}
\def\plotrline#1 #2 {
   {\dimen0=\plotcurrx \dimen2=\plotcurry
   \advance\dimen0 by #1\plotunitx \advance\dimen2 by #2\plotunity
   \plotDRline}
   \plotrmove {#1} {#2}
}
%
%
\def\plotvskip#1 {\vskip #1\plotunity\relax}
\def\plothphant#1 {\vbox to0pt{\hbox to#1\plotunitx{\hfil}}\plot!ni}
%
%
\def\plot!PSinit{
   \plotPSgdef\plot!circ
      gsave
      \plot!setorig
      newpath
      0 0 2.16 0 360 arc fill
      1 setgray
      newpath
      0 0 1.44 0 360 arc fill
      grestore
   \plotPSend
   \plotPSgdef\plot!bullet
      gsave
      \plot!setorig
      newpath
      0 0 2.16 0 360 arc fill
      grestore
   \plotPSend
   \plotPSgdef\plot!cdot
      gsave
      \plot!setorig
      newpath
      0 0 0.44 0 360 arc fill
      grestore
   \plotPSend
   \plotPSgdef\plot!tra \plot!trans \plotPSend
   \plotPSgdef\plot!pla gsave \plot!tra \plotPSend
   \plotPSgdef\plot!bfa \plot!rawsetcurr currentpoint transform\plotPSend
   \plotPSgdef\plot!bfb gsave \plot!rawsetorig itransform
      72 div exch 72 div exch scale\plotPSend
   \plotPSgdef\plot!bab gsave \plot!rawsetorig itransform
      72 div /p!s4 exch def 72 div /p!s3 exch def grestore gsave\plotPSend
   \plotPSgdef\plot!bac p!s3 p!s4 scale\plotPSend
   \plotPSgdef\plot!bad 1 p!s3 div 1 p!s4 div scale\plotPSend
   \plotPSgdef\plot!lia currentpoint newpath moveto lineto stroke\plotPSend
   \plot!reclaim
}
\def\plot!reclaim{
   \def\plotPSask{\immediate\write16{plot: PStype has already been set}}
   \def\plotPSfalse{\plotPSask}
   \def\plotPStype##1{\plotPSask}
   \let\plot!PSinit=\relax
   \ifplotPS
      \let\plot!TeXline=\relax
      \let\plot!TeXdotline=\relax
   \fi
   \let\plot!reclaim=\relax
}
\catcode`\!=\plot!savecc

%% file: plotsupp.tex
%
%
%
%
%
\def\PSsymb#1#2#3#4\plotPSend{{%
   \ifplotPS
      \escapechar=-1 %
      \expandafter\plotPSgdef\csname plot0\string#1\endcsname#4\plotPSend
      \gdef#1{{\escapechar=-1
         \plotPSmath\csname plot0\string#1\endcsname\plotPSend
         \dimen0=#2\mathchoice{\hskip\dimen0}{\hskip\dimen0}{\hskip.7\dimen0}%
            {\hskip.5\dimen0}\mathstrut}}%
   \else \gdef#1{{#3}} \fi}}
\plotPSgdef\enlargebox
   dup 3 1 roll add 5 1 roll
   dup 3 1 roll add 5 1 roll
   dup 3 1 roll sub 5 1 roll
   sub 4 1 roll
\plotPSend
\plotPSgdef\allcorners
   4 copy 4 -1 roll 5 -1 roll
\plotPSend

\def\rotateCWorCCW#1#2{{%
   \setbox0=\hbox{#2}%
   \dimen0=\ht0 \advance\dimen0 by\dp0 \dimen2=\wd0%
   \setbox0=\hbox{}%
   \hbox to\dimen0{\hfil
   \vrule width0pt height\dimen2 \plottrue\plotseparatefalse
   \plotfigbegin
   \plotcurrx=-.5\dimen0 \plotcurry=.5\dimen2 %
   \plotcentered{#2}\plotrotate #1 \plothere
   \plotfigend}}}

\def\setslope#1 #2 {\plotrotate {#2 #1 atan} }
\def\PSDRto#1\plotPSend{
   \ifplotPS
      {\plotfigscale 1
      \edef\dimsNTRNL{\plotfdims{\dimen0}{\dimen2}}
      \plotPSbegin currentpoint \dimsNTRNL moveto translate currentpoint
         exch atan rotate #1\plotPSend}
   \fi
}
\def\PSto#1 #2 #3\plotPSend{
   {\dimen0=#1\plotunitx \dimen2=#2\plotunity \PSDRto #3\plotPSend}
   \plotmove {#1} {#2}
}
\def\PSrto#1 #2 #3\plotPSend{
   {\dimen0=\plotcurrx \dimen2=\plotcurry
   \advance\dimen0 by #1\plotunitx \advance\dimen2 by #2\plotunity
   \PSDRto #3\plotPSend}
   \plotrmove {#1} {#2}
}
\plotPSgdef\PSarrow
   dup -4 mul 0 rmoveto gsave currentpoint newpath 0 0 moveto lineto stroke
   grestore dup scale currentpoint newpath moveto -1 -1.5 rlineto
   5 1.5 rlineto -5 1.5 rlineto closepath fill
\plotPSend
\def\plotDRarrow{
   \ifplotPS
      {\dimen8=1bp \edef\argNTRNL{\plotpdim{\dimen8}\PSarrow}
      \expandafter\PSDRto\argNTRNL\plotPSend}
   \else\plotDRline\fi}
\def\plotarrow#1 #2 {
   {\dimen0=#1\plotunitx \dimen2=#2\plotunity \plotDRarrow}
   \plotmove {#1} {#2}
}
\def\plotrarrow#1 #2 {
   {\dimen0=\plotcurrx \dimen2=\plotcurry
   \advance\dimen0 by #1\plotunitx \advance\dimen2 by #2\plotunity
   \plotDRarrow}
   \plotrmove {#1} {#2}
}
\def\labeledDRline#1#2#3{{
   \dimen4=\dimen0
   \dimen6=\dimen2
   \advance\dimen4 by-\plotcurrx
   \advance\dimen6 by-\plotcurry
   \advance\plotcurrx by.5\dimen4
   \advance\plotcurry by.5\dimen6
   #3
   \advance\dimen0 by-\dimen4
   \advance\dimen2 by-\dimen6
   #1
   \def\emptyNTRNL{}
   \def\argNTRNL{#2}
   \ifx\emptyNTRNL\argNTRNL\else
      \edef\dimsNTRNL{\plotthedim{\dimen4}\plotthedim{\dimen6}}
      \expandafter\setslope\dimsNTRNL
      \plottext{#2}
      \plothere
   \fi
}}
\def\labeledline#1#2#3 #4 #5 {
   {\dimen0=#4\plotunitx \dimen2=#5\plotunity \labeledDRline{#1}{#2}{#3}}
   \plotmove {#4} {#5}
}
\def\labeledrline#1#2#3 #4 #5 {
   {\dimen0=\plotcurrx \dimen2=\plotcurry
   \advance\dimen0 by #4\plotunitx \advance\dimen2 by #5\plotunity
   \labeledDRline{#1}{#2}{#3}}
   \plotrmove {#4} {#5}
}